\documentclass[12pt]{amsart}

\usepackage{mathtools, amsthm, amsfonts, enumerate, comment, fancyhdr, mathrsfs, tikz-cd, xparse, amscd,bm}
\usepackage{tensor}

\usepackage[margin=3cm]{geometry}

\usepackage{tikz, adjustbox, mathtools}
\usepackage[hidelinks]{hyperref}
\usetikzlibrary{matrix,calc}

\usepackage{cite}
\renewcommand{\MR}[1]{}

\newcommand{\owedge}{\mathrel{\begin{picture}(8,8)\put(4,2){\circle{8}}\put(0.3,-0.6){\scalebox{0.92}{$\wedge$}}\end{picture}}}
\newcommand{\loplus}{\mathrel{\begin{picture}(8,8)
\qbezier(4.03,8)(2.34,8)(1.17,6.83)
\qbezier(1.17,6.83)(0,5.66)(0,4)
\qbezier(0,4)(0,2.34)(1.17,1.17)
\qbezier(1.17,1.17)(2.34,0)(4.03,0)
\put(-0.55,1){$+$}\end{picture}}}

\newcounter{teller}

\newenvironment{tabel}{\begin{list}
{\rm  (\alph{teller})\hfill}{\usecounter{teller} \leftmargin=1.1cm
\labelwidth=1.1cm \labelsep=0cm \parsep=0cm}
                      }{\end{list}}

\newcounter{tellerr}

\newcounter{tellerrr}

\newcounter{thms} \numberwithin{thms}{section}
\numberwithin{equation}{section}

\newtheorem{theorem}[thms]{Theorem}
\newtheorem{proposition}[thms]{Proposition}
\newtheorem{lemma}[thms]{Lemma}
\newtheorem{corollary}[thms]{Corollary}

\newtheorem{theoremA}{Theorem}

\newtheorem*{previousresult}{Theorem}

\theoremstyle{remark}

\theoremstyle{definition}

\allowdisplaybreaks 

\author{Simon C. W. Goodwin}
\address[Goodwin]{Department of Mathematics, The University of Auckland, Private Bag 92019, Auckland 1142, New Zealand}
\email{simon.goodwin@auckland.ac.nz}
\author{A. Rod Gover} 
\address[Gover]{Department of Mathematics, The University of Auckland, Private Bag 92019, Auckland 1142, New Zealand}
\email{r.gover@auckland.ac.nz}
\author{Thomas Leistner}
\address[Leistner]{School of Mathematical Sciences, Adelaide University, SA 5005, Australia}
\email{thomas.leistner@adelaide.edu.au}
\thanks{
This work was supported in part by the Marsden Fund of the Royal Society Te Apārangi (grants 19-UOA-008 and 24-UOA-005) and by the Australian Research Council (Discovery Projects DP190102360 and DP260102671). This work was also supported in part by the Simons Foundation (grant SFI-MPS-T-Institutes-00010825) and by the State Treasury funds under the project ``Organization of the Simons Semesters at the Banach Center -- New Energies in 2026--2028'' (MNiSW/2025/DAP/491).  Part of this work was carried out during visits by the first named author to Adelaide University and by the third named author to the University of Auckland; they thank their respective hosts for their hospitality.}
\begin{document}
\title{Killing Scales for Conformal Killing Tensors of Rank Two and Higher}
\begin{abstract}
	We study which metrics in a conformal class admit a Killing tensor with prescribed conformal Killing trace-free part. We call such metrics Killing scales. We also consider whether a given conformal Killing tensor is itself a Killing tensor for some metric in the conformal class, as well as the analogous problem for conformal Killing–Yano forms. Our approach to this problem uses conformal and projective tractors, and the prolongations provided by the BGG splitting operators of the projective and conformal Killing operators. 
	For rank two conformal Killing tensors, we characterise local Killing scales for arbitrary scales on arbitrary conformal manifolds, extending recent results that were restricted to Einstein scales on conformally flat manifolds. The criterion is expressed in terms of the scale tractor and the BGG prolongation of the conformal Killing tensor. We also obtain Killing-scale criteria for higher-rank conformal Killing tensors and conformal Killing–Yano forms, under additional hypotheses involving Einstein scales and normal solutions. We describe the tractor geometry of special Einstein products, including Ricci-flat products, and use these constructions to give examples beyond the conformally flat setting.
\end{abstract}
\maketitle

\section{Introduction}

A symmetric tensor $r_{a_{1}\cdots a_\ell}\in \Gamma(S^\ell(T^*M))$ on a pseudo-Riemannian manifold $(M,g)$ is a Killing tensor if
\[
	\nabla{}_{(a_{0}} r_{a_{1}\cdots a_\ell)}=0,
\] 
where $\nabla$ is the Levi-Civita connection of $g$. Killing tensors give constants of motion for the geodesics of $g$: If $\gamma\colon [a,b]\to M $ has velocity $u = \dot{\gamma}$ and $\nabla_{u}u=0$, then
\[
	\frac{d}{dt}(r(u,\dots,u)) = (\nabla r)(u,\dots ,u) = 0
\]
since the symmetric part of $\nabla r$ vanishes.

The trace-free part $k^g{}_{a_{1}\cdots a_\ell} $ of a Killing tensor $r_{a_{1}\cdots a_\ell}$ is a solution of the \textit{conformal Killing equation}
\begin{equation}\label{Eqn:CKTE}
	\nabla{}_{(a_{0}} k^g{}_{a_{1}\cdots a_\ell)_0}=0,
\end{equation} 
where the $0$ subscript indicates the trace-free part (see Section \ref{Sec:MLA}). The equation \eqref{Eqn:CKTE} is \textit{conformally covariant}, that is, when $\hat{g} = \Omega^2 g$ for some positive function $\Omega$, the tensor $k^{\hat{g}} = \Omega^{2\ell}k^{g}$ will satisfy \eqref{Eqn:CKTE} for the metric $\hat{g}$. 

Let $(M,\boldsymbol{c})$ be a \textit{conformal manifold} of dimension $n\ge 3$. The conformal structure is described by the conformally weighted \textit{conformal metric} $\boldsymbol{g}{}_{ab} \in \Gamma(\mathcal{E}{}_{(ab)} [2])$, where the $\mathcal{E}[w]$ are density bundles, see Section \ref{Sec:ConGeo}. Each metric $g\in \boldsymbol{c}$ is determined by a positive scale density $\sigma\in \Gamma(\mathcal{E}_+[1])$ via $g = \sigma^{-2}\boldsymbol{g}$. The \textit{conformal Killing operator}, acting on sections $k_{a_{1}\cdots a_\ell}\in \mathcal{E}_{(a_{1}\cdots a_\ell)_0}[2\ell]$,
\begin{equation}\label{Eqn:CKOp}
	k\mapsto \nabla{}_{(a_{0}}k_{ a_{1}\cdots a_\ell)_0}
\end{equation} 
is conformally invariant, that is, the choice of Levi-Civita connection of a metric in the conformal class above does not change the output of the operator. A $k$ in the kernel of \eqref{Eqn:CKOp} is called a \textit{conformal Killing tensor}. 

If $g = \sigma^{-2}\boldsymbol{g}$, a solution $k{}^{g} $ of \eqref{Eqn:CKTE} determines the weighted tensor $k = \sigma^{2\ell} k^g$ in the kernel of \eqref{Eqn:CKOp}. Under a change of metric $\hat{g} = \Omega^2g$, the scale becomes $\hat{\sigma} = \Omega^{-1}\sigma$, so the transformation law for $k^{\hat{g}}$ gives
\[
	\hat{\sigma}^{2\ell}k^{\hat{g}} = (\Omega^{-1}{\sigma})^{2\ell}\Omega^{2\ell}k^{{g}} = \sigma^{2\ell} k^g.
\]
Thus the weighted tensor $k$ is independent of the choice of metric. We now ask the converse question. Given a conformal Killing tensor $k$, for which metrics $g\in \boldsymbol{c}$ is $k{}^{g}$ the trace-free part of a Killing tensor for $g$? We call such a metric a \textit{Killing scale} for $k$. A main goal of this paper is to characterise Killing scales for a given conformal Killing tensor.

We also consider the form analogue of Killing tensors, introduced by Yano \cite{yano_remarks_1952}. An alternating $\ell$-form $\omega_{a_{1}\cdots a_\ell}$ ($1\le \ell\le n-1$) on a pseudo Riemannian manifold is a \textit{Killing--Yano form} if it satisfies
\[
\nabla{}_{(a_{0}}\omega{}_{a_{1})a_{2}\cdots a_\ell}=0.  
\] 
The Killing--Yano equation is not conformally invariant, but (when given conformal weight $\ell+1$) a Killing-Yano form solves the invariant \textit{conformal Killing--Yano equation} (see \cite{semmelmann_conformal_2003} for a modern account)
\begin{equation}\label{Eqn:CKYEq}
\operatorname{trace-free}\left( \nabla{}_{(a_{0}} \omega{}_{a_{1})a_{2}\cdots a_\ell}   \right) =0,
\end{equation} 
where the $\operatorname{trace-free}$ of a tensor with any symmetries is defined in Section \ref{Sec:MLA}. Since a solution of \eqref{Eqn:CKYEq}, $\omega\in \Gamma(\mathcal{E}{}_{[a_{1}\cdots a_\ell]}[\ell+1] )$ is skew, it is already trace-free. So, we call $g \in \boldsymbol{c}$ a Killing scale for $\omega$ if $\omega$ is a Killing--Yano form for $g$.

Killing tensors for a pseudo-Riemannian manifold $(M,g)$ give rise to first integrals of the geodesic flow. Conformal Killing tensors do determine first integrals for null geodesics~\cite{gover_distinguished_2022} since the trace terms are invisible. For a Killing scale, conformal Killing tensors provide first integrals of the full geodesic flow by restoring the appropriate trace terms. Rank two Killing tensors also govern separation of variables for the Hamilton-Jacobi equations \cite{kalnins_killing_1980} and superintegrability, where conformal Killing tensors with multiple Killing scales play a special role \cite{miller_classical_2013,goodwin_superintegrable_2026}. 

All of the above Killing and conformal Killing equations are finite type overdetermined partial differential equations \cite{branson_prolongations_2006}, and are naturally studied by prolongation. Prolongation replaces an overdetermined PDE with a connection on a vector bundle. Parallel sections of the \textit{prolongation bundle} then correspond to solutions of the PDE. A classic example, and the prototype for all the prolongations studied in this paper, is the Killing equation on one forms, $\nabla{}_{(a} k{}_{b)} = 0$. The prolongation bundle is then $T^*M\oplus \Lambda^2(T^*M)$, together with the connection
\[
	\nabla_a{}^{\text{prol}}\begin{pmatrix} k{}_{b} \\ \mu{}_{bc}  \end{pmatrix} = \begin{pmatrix} \nabla_ak_b - \mu{}_{ab} \\ \nabla_a \mu{}_{bc} -R{}_{bc}{}^{d}{}_{a} k_d \\  \end{pmatrix}.
\] 
The prolongation map then sends a Killing one-form $k{}_{b}$ to the pair $(k_b, \nabla{}_{[b} k{}_{c]}) $.

 The Riemannian Killing operators are first Bernstein--Gelfand--Gelfand (BGG) operators for projective differential geometry and the conformal operators are first BGG operators for conformal geometry. The BGG machinery \cite{cap_bernstein-gelfand-gelfand_2001, calderbank_differential_2001} takes as input a tractor bundle $\mathcal{V}$ (see Section \ref{Sec:ProjConf}) and constructs an invariant differential operator $\mathscr{D}_0\colon \Gamma(E_1)\to \Gamma(E_2) $, between sections of weighted tensor bundles. 

The BGG machinery also provides a canonical splitting operator, $L \colon \Gamma({E}_1)\to \Gamma(\mathcal{V}) $. A section $s$ solves $\mathscr{D}_0s = 0$ if and only if the \textit{prolongation} $L(s)$ is parallel for the \textit{prolongation connection} as derived in~\cite{hammerl_new_2012}. In the flat case this connection agrees with the standard tractor connection~\cite{bailey_thomass_1994}. 

Each metric $g = \sigma^{-2}\boldsymbol{g}\in \boldsymbol{c}$ is uniquely related to $\boldsymbol{g}$ by the positive density $\sigma\in \Gamma(\mathcal{E}_+[1])$. So, using the tractor $D$ (see Section \ref{Sec:ConTract}) each metric in the conformal class corresponds to its \textit{scale tractor} $I_A=D_A\sigma$, which is parallel for the tractor connection exactly when $g$ is Einstein. We also allow scale densities with zeros on a nowhere dense set, and consider the associated metric on the open set where $\sigma\neq 0$.

Our approach is to characterise the Killing scales by comparing the prolongation $K$ of the conformal Killing tensor with the scale tractor $I$. When both the prolongation and the scale tractor are parallel, the criterion is algebraic. In rank two we determine the differential terms required for arbitrary $K$ and $I$.

An Einstein scale $\sigma$, corresponding to a parallel scale tractor $I{}_{A} $ determines a connection preserving isomorphism between the projective tractor bundle of the projective manifold $(M\setminus\mathcal{Z}(\sigma),[\nabla^g])$ and the subbundle $I^{\perp}$ of the conformal tractor bundle \cite{gover_detecting_2014} (Theorem \ref{Thm:Iperp}). Through this isomorphism, we can reformulate Killing scale conditions in terms of the trace-free map acting on tractors in $I^{\perp}$. 

The question of Killing scales was studied in \cite{gover_conformal_nodate}, for rank one and two Killing tensors on conformally flat manifolds. Their main theorem for rank two Killing tensors is the following.
\begin{previousresult}[\!\!{\cite[Theorem 4.25]{gover_conformal_nodate}}]
	On a conformally flat manifold, let $k{}_{ab}$ be a conformal Killing tensor with full prolongation $K_{ABCD}$. Let $\sigma$ be an Einstein scale with corresponding scale tractor satisfying $|I|^2 = I_AI^A\neq 0$. Then $\sigma$ is a Killing scale for $k{}_{ab}$ if and only if
	\[
		I^AK_{AB[CD}I{}_{E]} =0.
	\] 
\end{previousresult} 

We extend and generalise this result. We show that the same algebraic criterion holds on curved conformal manifolds for normal conformal Killing tensors. We also admit null scale tractors, corresponding to Ricci-flat scales. A conformal Killing tensor is \textit{normal} if its full prolongation is parallel for the standard tractor connection. When $(M,\boldsymbol{c})$ is conformally flat, all conformal Killing tensors are normal.

\begin{theoremA}[Theorem \ref{Thm:IWI}]\label{Thm:Int:A}
	Let $(M,\boldsymbol{c})$ be a conformal manifold, let $k{}_{bd}  $ be a normal solution of the conformal Killing equation with prolongation $K{}_{ABCD} $ and let $I{}_{A} = {D}{}_{A} \sigma$ be a parallel tractor. Then $\sigma$ is a local Killing scale for $k{}_{bd} $ if and only if
	\begin{equation}\label{Eqn:Int:A}
		I{}^{A} K{}_{AB[CD} I{}_{E]} =0.
	\end{equation} 
	If there exists a parallel tractor $\hat{I}{}_{A} $ such that $\hat{I}^AI_A=1$, then $\sigma$ is a Killing scale for $k{}_{bd} $ if and only if \eqref{Eqn:Int:A} holds.
\end{theoremA}
Here a \textit{local Killing scale} is a globally defined scale where the trace adjustment may exist only locally. When $I$ is not null, we can always take $\hat{I} = |I|^{-2}I$.
We also derive the full required adjustment to \eqref{Eqn:Int:A} when either $K$ or $I$ is not parallel. For the tractor projectors $X$ and $Z$ and the tractor $D$ operator see 
Section \ref{Sec:ConTract}.
\begin{theoremA}[Theorem \ref{Thm:LKS} and Corollary \ref{Cor:Rk2}]
	Let $(M,\boldsymbol{c})$ be a conformal manifold, let $k{}_{bd}  $ be a conformal Killing tensor with prolongation $K{}_{ABCD} $. Let $\sigma\in \Gamma(\mathcal{E}[1])$ be a scale with scale tractor $I{}_{A}$. Then $\sigma$ is a local Killing scale for $k{}_{bd} $ if and only if
	\begin{equation}
		X{}^{E} Z{}^{A}{}_{a} Z{}^{B}{}_{b}  X{}^{D} \left(3I{}^{C}I{}_{[E} K_{AB]CD}{}-\sigma(D{}_{[E} I{}^{C}) K{}_{AB]CD} - \sigma I{}^{C} D{}_{[E} K{}_{AB]CD} \right)=0.
	\end{equation}
	When $K{}_{ABCD}$ is parallel, this reduces to
	\begin{equation}
		X{}^{E} Z{}^{A}{}_{a} Z{}^{B}{}_{b}  X{}^{D} \left(3I{}^{C}I{}_{[E} K{}_{AB]CD} - \sigma (D{}_{[E} I{}^{C}) K{}_{AB]CD} \right)=0.
	\end{equation}
\end{theoremA}
Together with Jonathan Kress, the first two authors have applied this criterion to superintegrable systems in the case where $K$ is parallel, with no parallelism assumption on $I$ \cite{goodwin_superintegrable_2026}.

When considering Killing scales for conformal Killing tensors of rank $\ell\ge 3$, the space of possible trace adjustments is symmetric tensors of rank $\ell-2$. For a conformal Killing tensor of rank \(\ell\geq 3\), we seek a Killing tensor of the form
\[
	r_{a_{1}\cdots a_\ell} = k^g_{a_{1}\cdots a_\ell}+{g}_{(a_{1}a_{2}}q_{a_{3}\cdots a_\ell)},\quad \quad q\in \Gamma(S^{\ell-2}T^*M).
\] 
Since \(\nabla^g g=0\), the Killing equation for \(r\) becomes
\[
	\nabla^g{}_{(a_1}q_{a_2\cdots a_{\ell-1})}=-\frac{\ell}{n+2\ell-2}\,\nabla^{g}{}_b(k^g){}_{a_1\cdots a_{\ell-1}}{}^b.
\]
Thus \(g\) is a Killing scale precisely when the divergence of \(k^g\) lies in the image of the rank \(\ell-2\) projective Killing operator. In the projectively flat case the BGG sequence is locally exact, so, this image is the kernel of the generalised Calabi operator (also known as the generalised Saint-Venant operator in the study of elasticity).

\begin{theoremA}[Proposition \ref{Prp:HR}]
	Let $(M,\boldsymbol{c}) $ be a conformal manifold. Let $\ell\ge 3$. Let $k{}_{a_{1}\cdots a_\ell} \in \Gamma(\mathcal{E}_{(a_{1}\cdots a_\ell)_0}[2\ell])$ be a conformal Killing tensor and let $\sigma\in \Gamma(\mathcal{E}[1])$. Then $\sigma$ is a Killing scale for $k$ if and only if, after using $\sigma$ to trivialise densities, $\nabla^g_ak{}^a{}_{a_{2}\cdots a_\ell} $ is in the image of the rank $\ell-2$ projective Killing operator.

	When $(M,\boldsymbol{p}=[\nabla^g])$ is projectively flat, $\sigma$ is a local Killing scale if and only if $\mathscr{D}_1(\nabla^g_ak{}^a{}_{a_{2}\cdots a_\ell})=0$, where $\mathscr{D}{}_{1}$ is the next operator in this projective BGG sequence including the rank $\ell-2$ projective Killing operator.
\end{theoremA}
Recent work \cite{costanza_range_2025,costanza_calabi_2026} has made progress on characterising the image of the rank one Killing operator, which, by the above result finds an application to Killing scales of rank 3 conformal Killing tensors.

Despite the increased complexity of higher rank, we prove a higher rank algebraic analogue of Theorem \ref{Thm:Int:A} for Einstein scales.

\begin{theoremA}[Theorem \ref{Thm:HO1}]\label{Thm:Int:B}
	Let $(M,\boldsymbol{c})$ be conformally flat. Let $k{}_{a_{1}\cdots a_\ell}\in \Gamma(\mathcal{E}{}_{(a_{1}\cdots a_\ell)_0}[2\ell] ) $ be a conformal Killing tensor and let $K{}_{A_{1} B_{1}\cdots A_\ell B_\ell} $ be its full prolongation. Let $I_A = D_A\sigma$ be a non-null parallel scale tractor. 

	Then $\sigma$ is a Killing scale for $k$ if and only if the tractor
	\begin{equation}
		I{}^{A_{1}}K{}_{A_{1}B_{1}A_{2}B_{2}\cdots A_\ell B_\ell}
	\end{equation} 
	vanishes after skewing each pair $A_iB_i$ with $I_{C_i}$ for each $i\ge 2$; that is, with an abuse of notation,
	\begin{equation}\label{Eqn:Int:B}
		 I{}^{A_{1}}K{}_{A_1B_{1}[A_{2}B_{2}}I_{C_{2}]}\cdots {}_{[A_\ell B_\ell}I_{C_{\ell}]}=0,
	 \end{equation} 
	where all of the indices $A_i B_i$ belong to $K$.
\end{theoremA}
We show that the sufficiency direction extends to normal solutions on arbitrary conformal manifolds by giving an explicit construction of the prolongation of the Killing tensor $r$, from the tractor $K$.
\begin{theoremA}[Theorem \ref{Thm:HO2}]\label{Thm:Int:C}
	Let $(M,\boldsymbol{c}) $ be a conformal manifold and let $k_{a_{1}\cdots a_\ell}\in \Gamma(\mathcal{E}{}_{(a_{1}\cdots a_\ell)_0}[2\ell] )$ be a normal solution of the conformal Killing equation with (parallel) full prolongation $K_{A_{1}B_{1}\cdots A_\ell B_\ell}$. Let $I_A = D_A\sigma$ be a non-null parallel scale tractor. If \eqref{Eqn:Int:B} holds, then $\sigma$ is a Killing scale for $k$.
\end{theoremA}

The Killing--Yano equations , with $1\le \ell\le n-1$, also fit into the BGG framework. For these equations, prolongations lie in the completely skew tractor bundles. Theorem \ref{Thm:KYKS} gives a criterion for a metric to be a Killing scale for a given conformal Killing–Yano form. Since Killing--Yano forms are skew, no trace adjustment is possible. Killing--Yano forms of any valency can be squared to obtain a rank two Killing tensor,
\[
	r_{ab} = \omega_{a}{}^{a_{2}\cdots a_\ell}\omega_{ba_{2}\cdots a_\ell},
\] 
and similarly, the trace-free part of the square of a conformal Killing--Yano form is a rank two conformal Killing tensor. We show that for a normal solution of the conformal Killing--Yano equation, the resulting conformal Killing tensor is also normal, and we construct its full prolongation in Theorem \ref{Thm:KYProd}.

This paper is organised as follows. In Section \ref{Sec:ProjConf}, we introduce projective and conformal differential geometry, their tractor bundles, and the BGG machinery. In Section \ref{Sec:SEP}, we discuss special Einstein products, such as the conformal product manifolds discussed in \cite{armstrong_definite_2007, gover_sub-product_2009, gover_class_2010}. We give an explicit construction of the tractor bundle, connection and projectors for these products, in a way which includes Ricci-flat products. In Section \ref{Sec:KT}, we give our results for rank two conformal Killing tensors. Section \ref{Sec:CKY} contains our results for conformal Killing--Yano forms and Section \ref{Sec:HO} contains results for higher order Killing tensors of all ranks $\ell\ge 2$. In Section \ref{Sec:Ex}, we use special Einstein products to construct examples that are not conformally flat.

\section{Conformal and Projective Geometry}\label{Sec:ProjConf}
Throughout this paper we will frequently need to consider tensors with symmetries given by a Young diagram, including trace and trace-free parts. In the following section, we set out the notation and conventions we will use.
\subsection{Multilinear Algebra}\label{Sec:MLA}
Let $V$ be a finite dimensional real vector space. Let $v\in V$ and $\omega\in V^*$. The abstract index notation \cite{penrose_spinors_1984} attaches upper indices $v{}^{a}$ to vectors and lower indices $\omega_a$ to covectors. Multiple distinct indices indicate higher rank tensors, such as $\eta_{abc}\in \otimes^3(V^*)$ or $\phi{}_{a}{}^{b} \in V\otimes V^*$.

Abstract indices allow for compact notation of many natural operations. Such as trace and tensor product,
\[
	\omega(v) = \omega_av^a,\quad \phi{}_{a}{}^{a} = tr(\phi),\quad  v\otimes u = v^au^b,
\] 
The identity endomorphism of $V$ is written $\delta{}_{a}{}^{b} $, satisfying $\delta{}_{a}{}^{b} v^a = v^b$ for all $v\in V$.

Symmetrisation and skew-symmetrisation across abstract indices $a\cdots c$ is written $(a \cdots c)$ and $[a\cdots c]$, respectively, with the normalisation $\frac{1}{\ell!}$, where there are $\ell$ indices in $a\cdots c$. For example,
\[
	\omega{}_{[ab]}  = \tfrac{1}{2}\left(\omega{}_{ab} - \omega{}_{ba} \right), \quad u{}^{(a} v{}^{b)}  = \tfrac{1}{2}\left(u^av^b+u^bv^a\right),
\] 
\[
\eta{}_{(abc)}  = \tfrac{1}{6}\left( \eta{}_{abc} +\eta{}_{acb} +\eta{}_{bac} +\eta{}_{bca} +\eta{}_{cab} +\eta{}_{cba}  \right).
\] 
The group $GL(V)$ and its subgroups act on $V$, and its tensor products.
Trace and symmetrisation are both $GL(V)$-equivariant maps, so they define $GL(V)$-submodules of tensor products, such as
\[
	S^2(V) = \{u{}^{ab} \in V\otimes V: u{}^{ab} = u{}^{(ab)} \},\quad \Lambda^2(V) = \{u{}^{ab} \in V\otimes V: u{}^{ab}  = u{}^{[ab]} \}
\]
and $\operatorname{trace-free}(V\otimes V^*) :=  \{\eta{}^{a}{}_{b} \in V\otimes V^*: \eta{}^{a}{}_{a} =0\} $. The first pair give the full $GL(V)$-equivariant decomposition of $V\otimes V$ into $GL(V)$-irreducible submodules,
\[
	V\otimes V = S^2(V)\oplus \Lambda^2(V)
\]
Adding the complimentary submodule to the trace-free endomorphisms gives the full decomposition of $V\otimes V^*$ into $GL$-irreducible submodules
\[
	V\otimes V^* = \operatorname{trace-free}(V\otimes V^*)\oplus \mathbb{R}\cdot \delta{}_{a}{}^{b}
\]

For higher valence tensor spaces, we can use Young diagrams to determine all irreducible $GL(V)$-equivariant submodules. A Young diagram \cite{fulton_young_1996} is a drawing of $\ell$ boxes, arranged in a connected left-justified grid such that the row lengths are non-increasing. For example,
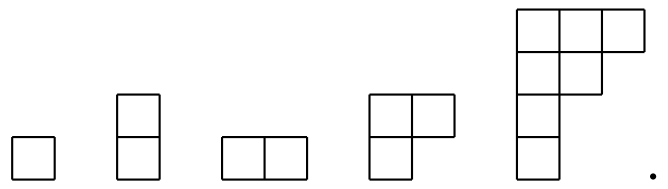
\[
\begin{picture}(8, 8)
	\put(0, 8){\line(1,0){8}}
	\put(0, 0){\line(1,0){8}}
	\put(0, 8){\line(0,-1){8}}
	\put(8, 8){\line(0,-1){8}}
\end{picture}\quad \begin{picture}(8, 16)
	\put(0, 16){\line(1,0){8}}
	\put(0, 8){\line(1,0){8}}
	\put(0, 0){\line(1,0){8}}
	\put(0, 16){\line(0,-1){16}}
	\put(8, 16){\line(0,-1){16}}
\end{picture}\quad \begin{picture}(16, 8)
	\put(0, 8){\line(1,0){16}}
	\put(0, 0){\line(1,0){16}}
	\put(0, 8){\line(0,-1){8}}
	\put(8, 8){\line(0,-1){8}}
	\put(16, 8){\line(0,-1){8}}
\end{picture}\quad \begin{picture}(16, 16)
	\put(0, 16){\line(1,0){16}}
	\put(0, 8){\line(1,0){16}}
	\put(0, 0){\line(1,0){8}}
	\put(0, 16){\line(0,-1){16}}
	\put(8, 16){\line(0,-1){16}}
	\put(16, 16){\line(0,-1){8}}
\end{picture}\quad \begin{picture}(24, 32)
	\put(0, 32){\line(1,0){24}}
	\put(0, 24){\line(1,0){24}}
	\put(0, 16){\line(1,0){16}}
	\put(0, 8){\line(1,0){8}}
	\put(0, 0){\line(1,0){8}}
	\put(0, 32){\line(0,-1){32}}
	\put(8, 32){\line(0,-1){32}}
	\put(16, 32){\line(0,-1){16}}
	\put(24, 32){\line(0,-1){8}}
\end{picture}.
\] 

Each Young diagram represents an isomorphism class of irreducible representations of $GL(V)$. The Young diagram with one box \begin{picture}(8, 7.0)(0, 1.0)
	\put(-0.2, 8){\line(1,0){8.4}}
	\put(-0.2, 0){\line(1,0){8.4}}
	\put(0, 8){\line(0,-1){8}}
	\put(8, 8){\line(0,-1){8}}
\end{picture} represents $V$. A single row of $\ell$ boxes represents $S^\ell(V)$ and a single column of $\ell$ boxes represents $\Lambda^\ell(V)$. 

Another example which we will routinely use is
\[
\begin{picture}(10, 10)(0,2)
	\put(0, 10){\line(1,0){10}}
	\put(0, 5){\line(1,0){10}}
	\put(0, 0){\line(1,0){10}}
	\put(0, 10){\line(0,-1){10}}
	\put(5, 10){\line(0,-1){10}}
	\put(10, 10){\line(0,-1){10}}
\end{picture}(V^*) = \{R{}_{abcd} :R{}_{abcd} = R{}_{[ab][cd]}\quad \text{and} \quad R{}_{[abc]d} =0\}, 
\] 
which implies the further symmetry $R{}_{abcd} = R{}_{cdab} $, that is, $\begin{picture}(10, 8.0)(0, 2.0)
	\put(-0.2, 10){\line(1,0){10.4}}
	\put(-0.2, 5){\line(1,0){10.4}}
	\put(-0.2, 0){\line(1,0){10.4}}
	\put(0, 10){\line(0,-1){10}}
	\put(5, 10){\line(0,-1){10}}
	\put(10, 10){\line(0,-1){10}}
\end{picture}(V^*)\subset S^2(\Lambda^2(V^*))$. These are the symmetries of the Riemann curvature tensor. In Section \ref{Sec:HO}, we will also need \[\begin{picture}(26, 8.0)(0, 2.0)
\put(-0.2, 10){\line(1,0){25.4}}
	\put(-0.2, 5){\line(1,0){25.4}}
	\put(-0.2, 0){\line(1,0){25.4}}
	\put(0, 10){\line(0,-1){10}}
	\put(5, 10){\line(0,-1){10}}
	\put(9, 2.3){${}_{\cdots}$}
	\put(9, 7.3){${}_{\cdots}$}
	\put(20, 10){\line(0,-1){10}}
	\put(25, 10){\line(0,-1){10}}
\end{picture}(V^*) = \{R_{a_{1}b_{1}a_{2}b_{2}\cdots a_\ell b_\ell}\in S^\ell(\Lambda^2(V^*)): R_{[a_{1}b_{1}a_{2}]b_{2}\cdots a_\ell b_\ell}=0\}.
\]
Here $R$ is skew in each pair $a_ib_i$ and symmetric under interchange of any two pairs. The machinery of abstract indices passes immediately from the vector space $V$ to a vector bundle $E\to M$ with fibre $V$ by considering all operations, tensor products, trace and symmetrisation fibrewise.

A symmetric non-degenerate bilinear form $g{}_{ab} \in S^2(V^*)$, with inverse $g^{ab}$ determines an isomorphism $V\cong V^*$, which we write $v^a\mapsto v_a := g{}_{ab} v^b$ and the inverse by $\omega_a\mapsto \omega^a := g{}^{ab} \omega_b$. Using $g$, we can define a trace on pairs of upper indices $u^{ab}\mapsto g_{ab}u^{ab}$ and similarly for pairs of lower indices. 

Raising, lowering and tracing indices are all $O(V,g)$-equivariant operations. These, combined with $GL(V)$-equivariant maps, give the $O(V,g)$-irreducible decompositions of tensor powers of $V$. Let $\lambda$ be a Young diagram and let $V_\lambda$ be the tensor space determined by $\lambda$ such as the above examples. Then we define

\[
	V^0_\lambda := \{u^{a_{1}\cdots a_\ell}\in V^\lambda : u^{a_{1}\cdots a_\ell} g_{a_i a_j}=0\textrm{ for all }1\le i<j\le \ell \} 
\] 
Then $V^0_\lambda$ is a $O(V,g)$-irreducible submodule of $\otimes^\ell(V)$, as long as $\lambda$ is in the \textit{stable range} of Young diagrams defined by the first two columns lengths summing to at most $n$. Later, we will use a Young diagram with a $0$ subscript instead of $V_\lambda$. For example,
\[
\begin{picture}(16, 10)(0,2)
	\put(0, 10){\line(1,0){10}}
	\put(0, 5){\line(1,0){10}}
	\put(0, 0){\line(1,0){10}}
	\put(0, 10){\line(0,-1){10}}
	\put(5, 10){\line(0,-1){10}}
	\put(10, 10){\line(0,-1){10}}
	\put(10, 0){${}_{\,0}$}
\end{picture}(V^*) = \{W{}_{abcd} :W{}_{abcd}  = W{}_{[ab][cd]},\quad W{}_{[abc]d} =0\quad \text{and} \quad W{}_{ab}{}^{a}{}_{d} =0\}.
\] 
These are the symmetries of the conformal Weyl curvature, see Section \ref{Sec:ConGeo}. The bilinear form $g$ also defines a non-degenerate bilinear form on $V{}_{\lambda} $ for any $\lambda$ given by tracing all indices,
\[
	\langle \eta^{a\cdots c}, \xi^{a\cdots c} \rangle \mapsto \eta^{a_{1}\cdots c_{1}}\xi^{a_{2}\cdots c_{2}}g_{a_{1}a_{2}}\cdots g_{c_{1}c_{2}}
\]

The space of trace tensors is then the subspace of $V^\lambda$ orthogonal to $V^\lambda_0$.
\[
	V^\lambda_g := (V^\lambda_0)^{\perp} = \{u^{a_{1}\cdots a_\ell}\in V^\lambda:u_{a_{1}\cdots a_\ell}w^{a_{1}\cdots a_\ell}=0\text{ for all } w\in V^\lambda_0  \},\quad V^\lambda = V^\lambda_0\oplus V^\lambda_g.
\] 
The trace subspace $V^\lambda_g$ is generally not irreducible with respect to the action of $O(V,g)$ (the full branching rules are due to Littlewood, see e.g. \cite[\S 25.3]{fulton_representation_2004}). This direct sum gives a canonical $O(V,g)$-equivariant projection
\[
	\operatorname{trace-free}\colon V{}^{\lambda} \to V^\lambda_0
\]
For the special case of totally symmetric $s\in S^\ell(V^*)$, we use the standard notation
\[
	s{}_{(a_{1}\cdots a_\ell)_0} = \operatorname{trace-free}(s{}_{(a_{1}\cdots a_\ell)})  
\] 
For example, trace-free map for $s{}_{ab} \in S^2(V^*)$ is
\[
\operatorname{trace-free}(s)_{ab} = s{}_{(ab)_0} =s{}_{ab} - \tfrac{1}{n}s{}_{c}{}^{c} g{}_{ab}.
\] 
Later we will use a Young diagram with a $g$ subscript to indicate the trace part of the space determined by the Young diagram, for example, due to the Kulkarni--Nomizu product, we have that
\[
\begin{picture}(17, 10)
	\put(0, 10){\line(1,0){10}}
	\put(0, 5){\line(1,0){10}}
	\put(0, 0){\line(1,0){10}}
	\put(0, 10){\line(0,-1){10}}
	\put(5, 10){\line(0,-1){10}}
	\put(10, 10){\line(0,-1){10}}
	\put(10, 0){${}_{\,g}$}
\end{picture}(V^*) = \{g\owedge s = s_{ac}g_{bd}+s{}_{bd} g{}_{ac} -s{}_{ad} g{}_{bc} -s{}_{bc} g{}_{ad}  : s{}_{ac} \in S^2(V^*)\}.
\] 

Then the map $\operatorname{trace-free}\colon \begin{picture}(10, 8.0)(0, 2.0)
	\put(-0.2, 10){\line(1,0){10.4}}
	\put(-0.2, 5){\line(1,0){10.4}}
	\put(-0.2, 0){\line(1,0){10.4}}
	\put(0, 10){\line(0,-1){10}}
	\put(5, 10){\line(0,-1){10}}
	\put(10, 10){\line(0,-1){10}}
\end{picture}(V^*)\to \begin{picture}(16, 8.0)(0, 2.0)
	\put(-0.2, 10){\line(1,0){10.4}}
	\put(-0.2, 5){\line(1,0){10.4}}
	\put(-0.2, 0){\line(1,0){10.4}}
	\put(0, 10){\line(0,-1){10}}
	\put(5, 10){\line(0,-1){10}}
	\put(10, 10){\line(0,-1){10}}
	\put(10, 0){${}_{\,0}$}
\end{picture}(V^*) $ is given by
\[
	\operatorname{trace-free}(R)_{abcd} = R{}_{abcd} - P{}_{ac} g{}_{bd} -P{}_{bd} g{}_{ac} +P{}_{ad} g{}_{bc} +P{}_{bc} g{}_{ad},
\] 
where
\[
	P{}_{ac}  = \frac{1}{n-2}\left(R{}_{abc}{}^{b} - \frac{R{}_{ab}{}^{ab} }{2(n-1)}g{}_{ac} \right).
\]
In calculations, we usually want to understand the kernel of the $\operatorname{trace-free}$. Any tensor in $V{}^{\lambda} $ which has two free metric indices, such as $\eta{}_{(ab} g{}_{cd)} $, is in the kernel of $\operatorname{trace-free}$.

When considering multiple vector spaces or vector bundles at the same time, we will use indices from different alphabets, such as capital Latin letters, $A,B,C$ or lower case Greek letters, $\alpha,\beta,\gamma$. So if we use indices $a,b,c$ for $V$ and $A,B,C$ for $W$, then $\phi\in \operatorname{Hom}(V,W)$ would have indices $\phi{}_{a}{}^{A} $.

In particular, in later sections, lower case Latin indices will be reserved for sections of the tangent bundle, $v^a\in \Gamma(TM)$ and $\omega_a\in \Gamma(T^*M)$, capital Latin indices will be used for conformal tractors and lower case Greek indices for projective tractors.

We also use indices to denote the vector bundles themselves, the trivial bundle is denoted $\mathcal{E} = M\times\mathbb{R}$, the tangent and cotangent bundles are $\mathcal{E}^a = TM$ and $\mathcal{E}_a = T^*M$. Higher valence bundles have more indices, e.g., $\mathcal{E}{}_{a}{}^{b}  = \mathcal{E}_a\otimes\mathcal{E}^b = T^*M\otimes TM = \operatorname{End}(TM)$. Symmetries can also be encoded in this notation, so $s{}^{ab} \in \Gamma(\mathcal{E}{}^{(ab)} )$ is a symmetric contravariant tensor field, and $\omega{}_{abc} \in \Gamma(\mathcal{E}{}_{[abc]} )$ is a $3$-form.

\subsection{Projective Geometry}\label{Sec:ProGeo}
Let $M$ be a smooth manifold of dimension $n\ge 2$. Two torsion free affine connections on a manifold are \textit{projectively equivalent} if they have the same geodesics as unparametrised curves. An equivalence class of projectively equivalent connections is called a \textit{projective class}, usually denoted $\boldsymbol{p} = [\nabla]$. The pair $(M,\boldsymbol{p})$ is a \textit{projective manifold}. 

Two torsion free affine connections $\nabla,\hat{\nabla}$ are projectively equivalent if and only if they related by a one form $\Upsilon{}_{a}\in \Gamma(\mathcal{E}_a)$:
\begin{equation}\label{Eqn:ProEq}
\hat{\nabla}_a\omega_b = \nabla_a\omega_b-\Upsilon_a\omega_b-\omega_a\Upsilon_b
\end{equation} 
for all $\omega\in \Gamma(\mathcal{E}_a)$. Conversely, two torsion free affine connections related by a one form as above are projectively equivalent. Other transformation laws then follow, such as
\[
\hat{\nabla}_au^b = \nabla_au^b+\Upsilon_au^b+\Upsilon_cu^c \delta{}_{a}{}^{b} 
\] 
The transformation law for a top form $\epsilon\in \Gamma(\Lambda^n(T^*M))$ is
\[
	\hat{\nabla}_a\epsilon = \nabla_a\epsilon-(n+1)\Upsilon_a\epsilon
\] 
The bundle $(\Lambda^n(T^*M))^{\otimes2}$ carries a canonical orientation, where squares are positive.  We then define $\mathcal{E}(-2(n+1)) = \Lambda^n(T^*M)^{\otimes 2}$, \textit{the projective density bundle} of weight $-2(n+1)$, with the weight corresponding to the transformation law. We then let $\mathcal{E}(w) = \mathcal{E}(-2(n+1))^{-\frac{w}{2(n+1)}}$. These real powers are well defined and oriented, since the structure group $\mathcal{F}(TM)$ acts by positive scalars. The density bundles can alternatively be constructed as an associated bundle of the frame bundle using the representation $\rho_w\colon GL(n) \to \mathbb{R}_+$ defined by $\rho_w(A) = |\det A|^{w / (n+1)}$. An affine connection determines a connection on all bundles associated to the frame bundle, including the density bundles. Then for any $\tau\in \Gamma(\mathcal{E}(w))$ and connections $\nabla,\hat{\nabla}\in \boldsymbol{p}$ related by \ref{Eqn:ProEq},
\[
\hat{\nabla}_a\tau = \nabla_a\tau+w\Upsilon_a\tau.
\] 

Projective weights are very useful for defining invariant differential operators. From the transformation law for one forms, we see that $\hat{\nabla}_a\omega_b-\nabla_a\omega_b$ is symmetric, so the exterior derivative $\omega\mapsto \nabla{}_{[a} \omega{}_{b]} $ is projectively invariant. However, if we instead consider $\omega\in \Gamma(\mathcal{E}_a(2))$, the transformation law becomes
\[
	\hat{\nabla}_a\omega_b = \nabla_a\omega_b+\Upsilon_a\omega_b-\omega_a\Upsilon_b,
\] 
and the difference is skew symmetric. It follows that the \textit{rank-one Killing} operator  acting on sections $\omega_a\in \Gamma(\mathcal{E}_a(2))$:
\[
\omega_a\mapsto \nabla{}_{(a} \omega{}_{b)} 
\] 
is projectively invariant.

Another transformation law we need is for the Riemann curvature, defined by
\[
R{}_{ab}{}^{c}{}_{d} u^d = \nabla_a\nabla_bu^c-\nabla_b\nabla_au^c
\] 
since all connections are torsion free. If $\hat{\nabla}, \nabla\in \boldsymbol{p}$ are related by $\Upsilon$, then their curvatures are related by
\[
	\hat{R}{}_{ab}{}^{c}{}_{d}  = R{}_{ab}{}^{c}{}_{d} +2\delta{}_{[a}{}^{c} \Upsilon{}_{b]} \Upsilon{}_{d} +2\nabla{}_{[a} \Upsilon{}_{b]} \delta{}^{c}{}_{d} +2\delta{}_{[b}{}^{c}\nabla{}_{a]}  \Upsilon_d.
\] 
What is immediately discernible from the above is that the trace-free part of the curvature is projectively invariant. We can decompose the curvature tensor uniquely under the action of $SL_n(\mathbb{R})$, into one trace-free part, and two trace terms, one symmetric and one skew symmetric. 

\begin{proposition}\label{Prp:ProCur}
Let $\nabla$ be a torsion free affine connection with curvature tensor $R{}_{ab}{}^{c}{}_{d} $. Let
\[
\tilde{\mathsf{P}}{}_{bd}  = \frac{1}{n-1}\operatorname{Ric}{}_{(bd)} +\frac{1}{n+1}\operatorname{Ric}{}_{[bd]},\quad  \beta{}_{bd} = -\frac{2}{n+1}\operatorname{Ric}{}_{[bd]} = -2\tilde{\mathsf{P}}{}_{[bd]},  \] 
where $\operatorname{Ric}{}_{bd}  = R{}_{ab}{}^{a}{}_{d} $. Then
\[
R{}_{ab}{}^{c}{}_{d}  = \tilde{W}{}_{ab}{}^{c}{}_{d} + 2 \delta{}_{[a}{}^{c} \tilde{\mathsf{P}}{}_{b]d} +\beta{}_{ab} \delta{}^{c}{}_{d},
\]
where $\tilde{W}$ is completely trace-free. The tensors $\operatorname{Ric},\tilde{\mathsf{P}}$ and $\beta$ are related by
\[
	(n-1)\tilde{\mathsf{P}}{}_{ab}  = \operatorname{Ric}{}_{ab} +\beta{}_{ab}.
\] 
The tensor $\beta$ acts as the curvature on densities; that is, if $\sigma\in \Gamma(\mathcal{E}(w))$, then
\[
	2\nabla{}_{[a} \nabla{}_{b]} \sigma = w\beta{}_{ab} \sigma.
\] 
Let $\hat{\nabla},\nabla\in \boldsymbol{p}$ be related by $\Upsilon$. Then
	\[
		\hat{\tilde{W}}{}_{ab}{}^{c}{}_{d} = \tilde{W}{}_{ab}{}^{c}{}_{d},\quad \hat{\tilde{\mathsf{P}}}{}_{bd}  = \tilde{\mathsf{P}}{}_{bd} -\nabla_b\Upsilon_d+\Upsilon_b\Upsilon_d\quad \text{and} \quad \hat{\beta}{}_{bd}  = \beta{}_{bd} +2\nabla{}_{[b} \Upsilon{}_{d]} 
	\] 
\end{proposition}
The tensor $\tilde{W}$ is the \textit{projective Weyl tensor} and $\tilde{\mathsf{P}}$ is the \textit{projective Schouten tensor}. Since we will also be considering conformal manifolds in this article, which carry their own Weyl and Schouten tensors, we will use tildes to denote the projective tensors. There is no $\beta$ curvature tensor in conformal geometry, so $\beta$ does not require a tilde.

\subsection{Projective Tractors}
Projective manifolds lack a distinguished connection on the tangent bundle by design. The tractor construction associates to a projective manifold a natural vector bundle equipped with a connection.

The standard projective cotractor bundle is $\tilde{\mathcal{T}}^* = J^1(\mathcal{E}(1))$. We will use lower case Greek indices for projective tractors, so $\tilde{\mathcal{T}}^* = \mathcal{E}_\alpha$ and its dual $\tilde{\mathcal{T}} = \mathcal{E}^\alpha$. The tractor bundle fits into the canonical first jet exact sequence for $\mathcal{E}(1)$,
\begin{equation}\label{Eqn:PTS}
	\begin{tikzcd}[column sep=1.8em]
	0 \arrow[r]
	& \mathcal{E}{}_{a}(1) \arrow[r, "Z{}_{\alpha}{}^{a} "]
	& \mathcal{E}_\alpha \arrow[r, "X^\alpha"] 	
	& \mathcal{E}(1) \arrow[r] 	& 0.
	\end{tikzcd}
\end{equation}
The canonical projectors are sections $X^\alpha\in \Gamma(\mathcal{E}^\alpha(1))$ and $Z{}_{\alpha}{}^{a} \in \Gamma(\mathcal{E}{}_{\alpha}{}^{a} (-1))$. An affine connection $\nabla\in \boldsymbol{p}$ induces a connection on $\mathcal{E}(1)$, which is equivalent to a splitting of \eqref{Eqn:PTS} \cite[4.6.3]{saunders_geometry_1989}. So each $\nabla\in \boldsymbol{p}$ determines sections $Y_\alpha\in \Gamma(\mathcal{E}_\alpha(-1))$ and $W{}^{\alpha}{}_{a} \in \Gamma(\mathcal{E}{}^{\alpha}{}_{a} (1))$, which satisfy
\[
X^\alpha Y_\alpha=1\quad \text{and} \quad W{}^{\alpha}{}_{a} Z{}_{\alpha}{}^{b} =\delta{}_{a}{}^{b}.
\]
If $\nabla$ and $\hat{\nabla}$ are related by $\Upsilon{}_{a}$, then 
\[
	\hat{W}{}^{\alpha}{}_{a} = W{}^{\alpha}{}_{a} +X^\alpha \Upsilon_a\quad \text{and} \quad \hat{Y}_\alpha = Y_\alpha - Z{}_{\alpha}{}^{a} \Upsilon_a.
\]

Fixing a connection then determines isomorphisms $\mathcal{E}(1)\oplus \mathcal{E}_a(1)\cong \mathcal{E}_\alpha$, given by $(\sigma,\mu_a)\mapsto \sigma Y_\alpha+\mu_a Z{}_{\alpha}{}^{a} = U_\alpha$, which we write as $U_\alpha \overset{\nabla}{=} (\sigma,\mu_a)$, or as a vertical vector.

The projectively invariant tractor connection can be given through this isomorphism,
\[
\nabla_a\begin{pmatrix} \sigma\\ \mu{}_{a}  \end{pmatrix} = \begin{pmatrix} \nabla_a\sigma-\mu_a\\ \nabla_a\mu_b+\sigma\tilde{\mathsf{P}}{}_{ab}  \end{pmatrix},
\] 
although we will more commonly use the connection expressed in terms of the projectors,

	\[
		\nabla_a X^\alpha = W{}^{\alpha}{}_{a} ,\quad \nabla_a W{}^{\alpha}{}_{b}  =-\tilde{\mathsf{P}}{}_{ab} X^\alpha, \quad \nabla_a Z{}_{\alpha}{}^{b}  = -\delta{}_{a}{}^{b} Y_\alpha,\quad \nabla_a Y_\alpha =  \tilde{\mathsf{P}}_{ab}Z{}_{\alpha}{}^{b},
	\] 
since all of these projectors are weighted, the tractor connection applied to them is not invariant.

There is a first order projectively invariant differential operator $D\colon \Gamma(\mathcal{E}(1))\to \Gamma(\mathcal{E}{}_{\alpha} ) $, given by 
\begin{equation*}
	D_\alpha\sigma|_p = j^1_p(\sigma) = (\sigma,\nabla_a\sigma),
\end{equation*}
the projective Thomas $D$.

The projective cotractor bundle and connection form a prolongation bundle for the projectively invariant partial differential operator
\[
	\tilde{A}(\sigma){}_{ab} = \nabla{}_{(a} \nabla{}_{b)} \sigma + \sigma\tilde{\mathsf{P}}{}_{(ab)}, \]
so parallel sections of $\mathcal{E}{}_{\alpha} $ correspond to solutions of $\tilde{A}$. The correspondence is given by $\sigma \mapsto D_\alpha \sigma$.

\subsection{Conformal Geometry}\label{Sec:ConGeo}
From this point, we require $n\ge 3$. Two pseudo-Riemannian metrics $g{}_{ab} ,\hat{g}{}_{ab} \in \Gamma(\mathcal{E}{}_{(ab)} )$ are conformally equivalent if there is a positive smooth function $\Omega>0$ such that $\hat{g}{}_{ab}  = \Omega^2 g{}_{ab} $. An equivalence class of conformally equivalent metrics is called \textit{conformal class} and is usually denoted $\boldsymbol{c} = [g]$. The pair $(M,\boldsymbol{c})$ is a \textit{conformal manifold}.

Let $(M,\boldsymbol{c})$ be a conformal manifold and let $g,\hat{g}\in \boldsymbol{c}$. If $\hat{g} = \Omega^2g\in \boldsymbol{c}$, then the Levi-Civita connections of $g$ and $\hat{g}$ are related by
\[
	\hat{\nabla}_a\omega_b = \nabla_a\omega_b -\Upsilon_a\omega_b-\omega_a\Upsilon_b+\Upsilon_d\omega_c g{}^{cd} g{}_{ab},
\] 
where $\Upsilon_a = \Omega^{-1}\nabla_a\Omega$. The transformation law for vectors is
\[
	\hat{\nabla}_av^b = \nabla_av^b + \Upsilon_av^b-v_a\Upsilon^b+\Upsilon_cv^c\delta_a^b,
\] 
where indices are currently being raised and lowered using any metric in the conformal class.

Reminiscent of the projective case, we have for top forms, $\epsilon\in \Gamma(\Lambda^n(T^*M))$,
\[
\hat{\nabla}_a\epsilon = \nabla_a\epsilon -n\Upsilon_a\epsilon.
\] 
So we define the oriented density bundle $\mathcal{E}[-2n] = (\Lambda^n(T^*M))^2 = \mathcal{E}(-2(n+1))$, with roots $\mathcal{E}[w] = \mathcal{E}\left(\frac{n+1}{n}w\right)$ for any $w\in \mathbb{R}$. Then if $\tau\in \Gamma(\mathcal{E}[w])$, we have
\[
\hat{\nabla}_a\tau = \nabla_a\tau+w \Upsilon_a\tau.
\]
Any metric $g\in \boldsymbol{c}$ determines a unique density $\tau^g\in\Gamma((\Lambda^n(T^*M))^2) =  \Gamma(\mathcal{E}[-2n])$, the metric density, which when $M$ is oriented will simply be the square of the Riemannian volume form. The density $\tau^g$ is parallel for the Levi-Civita connection of $g$ and satisfies
\[
	\frac{1}{n!}\tau_{a_{1}\cdots a_nb_{1}\cdots b_n}g^{a_{1}b_{1}}\cdots g^{a_nb_n}=(-1)^q,
\] 
where $g$ has signature $(p,q)$.
So if $\hat{g} = \Omega^2g $, then $\tau^{\hat{g}} = \Omega^{2n}\tau^g$. We call $\sigma_g = (\tau^g)^{-1 / 2n}$ the scale density of $g\in \boldsymbol{c}$, which is a nowhere vanishing section $\sigma_g\in \Gamma(\mathcal{E}_+[1])$. If $\hat{g} = \Omega^2 g$, then $\sigma_{\hat{g}} = \Omega^{-1}\sigma_g$.

The \textit{conformal metric} is defined by $\boldsymbol{g}{}_{ab}  = \sigma_g^2 g{}_{ab} $ for any $g\in \boldsymbol{c}$. If $\hat{g} = \Omega^2 g$, then
\[
	\sigma_{\hat{g}}^2 \hat{g}{}_{ab}  = \Omega^{-2}\sigma_g^2 \Omega^2 g{}_{ab} = \sigma{}_{g}^2  g{}_{ab}
\]
So the conformal metric is an invariant weighted non-degenerate section; $\boldsymbol{g}{}_{ab} \in \Gamma(\mathcal{E}{}_{(ab)} [2])$ with inverse $\boldsymbol{g}{}^{ab} \in \Gamma(\mathcal{E}{}^{(ab)} [-2])$.

The conformal metric is also parallel for any metric in the conformal class, so we will use $\boldsymbol{g}$ to raise and lower indices from now on, meaning that raising an index lowers the weight of a tensor by two and lowering an index increases the weight by two.

A positive density $\sigma\in \Gamma(\mathcal{E}_+[1])$ determines a metric via $g{}_{ab} = \sigma^{-2}\boldsymbol{g}{}_{ab} $ with $\sigma_g = \sigma$. So metrics in $\boldsymbol{c}$ are in one-to-one correspondence with positive scale densities.

More generally if $\sigma\in \Gamma(\mathcal{E}[1])$ and the complement of the zero set of $\sigma$, $M\setminus\mathcal{Z}(\sigma)$ is dense in $M$, then we call $\sigma$ a \textit{scale}. The metric $g = \sigma^{-2}\boldsymbol{g}$ is well defined on the open dense set $M\setminus\mathcal{Z}(\sigma)$.

\begin{proposition}\label{Prp:ConCur}
	Let $g$ be a metric with curvature $R{}_{abcd} \in \Gamma(\mathcal{E}{}_{[ab][cd]} [2])$. Let
\[
	\mathsf{P}{}_{bd} = \frac{1}{n-2}\left( \operatorname{Ric}_{bd}-\frac{\operatorname{Sc}}{2(n-1)}\boldsymbol{g}{}_{bd}  \right),
\] 
where $\operatorname{Ric}{}_{bd}  = R{}_{ab}{}^{a}{}_{d} \in \Gamma(\mathcal{E}{}_{(bd)} )$ and $\operatorname{Sc} = R{}_{ab}{}^{ab} \in \Gamma(\mathcal{E}[-2])$. Then
\[
	R{}_{abcd} = W{}_{abcd} +2\mathsf{P}{}_{a[c} \boldsymbol{g}{}_{d]b} +2\mathsf{P}{}_{b[d} \boldsymbol{g}{}_{c]a},
\] 
where $W$ is completely trace-free. If $\Omega\in C^\infty_+(M) $ and $\hat{g} = \Omega^2g$, then
\[
\hat{W}{}_{abcd}  = W{}_{abcd}, \quad \hat{\mathsf{P}}{}_{bd} =  \mathsf{P}{}_{bd} - \nabla_b\Upsilon_d +\Upsilon_b\Upsilon_d-\tfrac{1}{2}\boldsymbol{g}{}_{bd} \Upsilon_c\Upsilon^c.
\] 
Let $J = \mathsf{P}{}_{b}{}^{b} \in \Gamma(\mathcal{E}[-2])$ and $\overset{\circ }{\mathsf{P}}{}_{bd}  = \mathsf{P}{}_{bd} -\frac{1}{n} J\boldsymbol{g}{}_{bd} $. Then
\[
	\hat{J} = J-\nabla_a\Upsilon^a+\frac{2-n}{2}\Upsilon_a\Upsilon^a\quad \text{and} \quad{\hat{\mathsf{P}}}_{(ab)_0} = {\mathsf{P}}_{(ab)_0} -\nabla{}_{(a} \Upsilon{}_{b)_0}+\Upsilon{}_{(a} \Upsilon{}_{b)_0}.
\] 
Let $(M,[\nabla^g])$ be the projective manifold containing the Levi-Civita connection of $g$. Then
\[
\tilde{\mathsf{P}}{}_{ab}  = 2\mathsf{P}{}_{ab} - \frac{n}{n-1}\overset{\circ }{\mathsf{P}}{}_{ab}  = \mathsf{P}{}_{ab} +\frac{J}{n}\boldsymbol{g}{}_{ab} -\frac{1}{n-1}\overset{\circ }{\mathsf{P}}{}_{ab} 
\] 
\end{proposition}

The conformally invariant tensor $W{}_{abcd} $ is called the \textit{conformal Weyl curvature}, it is completely trace-free by construction. The \textit{conformal Schouten curvature} is $\mathsf{P}{}_{bd} $.

A pseudo-Riemannian metric is \textit{Einstein} if its Ricci curvature is pure trace, $\operatorname{Ric}{}_{(ab)_0} =0$. Which is clearly equivalent to $\mathsf{P}{}_{(ab)_0}=0 $. A Lemma of Schur states that an Einstein metric has constant scalar curvature on connected components, that is, $\nabla_a \operatorname{Sc}=0$ or $\nabla_aJ=0$. 
\subsection{Conformal Tractors}\label{Sec:ConTract}
Our references for this section are \cite{bailey_thomass_1994,curry_introduction_2018}.
Einstein metrics play an central role in conformal geometry. The scale density of an Einstein metrics in the conformal class satisfies the conformally invariant differential equation
\[
 A{}_{ab} \sigma = \nabla{}_{(a} \nabla{}_{b)_0}\sigma+\mathsf{P}{}_{(ab)_0}\sigma =0.
\] 
Conversely, a scale $\sigma$ satisfying $A{}_{ab} \sigma=0$ determines an Einstein metric $g{}_{ab}  = \sigma^{-2}\boldsymbol{g}{}_{ab} $ on $M\setminus\mathcal{Z}(\sigma)$. Calculating using the Levi-Civita connection $\nabla{}^{\sigma}$ away from $\mathcal{Z}(\sigma)$ gives
\[
	A{}_{ab}  \sigma = \nabla{}_{(a} \nabla{}_{b)_0}\sigma+\mathsf{P}^g{}_{(ab)_0} \sigma = \mathsf{P}^g_{(ab)_0}\sigma=0
\] 
The partial differential operator $A$ is called the \textit{almost Einstein operator}, and it gives a conformally invariant partial splitting of the second jet exact sequence

\begin{equation}\label{Eqn:CTS}
\begin{tikzcd}[column sep = 1.8em]
0\arrow[r]
	&\mathcal{E}{}_{(ab)}[1]\arrow[r,"\iota"]
	&J^2(\mathcal{E}[1])\arrow[r]\arrow[l,bend left=43,"A", pos=0.53]
	&J^1(\mathcal{E}[1])\arrow[r]
	&0,
\end{tikzcd}
\end{equation}
where $A\circ \iota$ is the identity on $\mathcal{E}{}_{(ab)_0}[1] $. The \textit{conformal cotractor bundle} is $\mathcal{E}_A = \ker(A)\subset J^2(\mathcal{E}[1])$. We have $\ker(A) \cong J^2(\mathcal{E}[1]) / \iota(A(J^2(\mathcal{E}[1]))) = J^2(\mathcal{E}[1])/\iota(\mathcal{E}{}_{(ab)_0}[1]) $. The sequence \eqref{Eqn:CTS} reduces to the \textit{tractor exact sequence}
\[
\begin{tikzcd}[column sep = 1.8em]
0\arrow[r]
	&\mathcal{E}{}[-1]\arrow[r,"X_A"]
	&\mathcal{E}_A \arrow[r]
	&J^1(\mathcal{E}[1])\arrow[r]
	&0,
\end{tikzcd}
\]
using the identification $\mathcal{E}[-1] \cong \mathcal{E}[-1]\otimes\boldsymbol{g}{}_{ab} \subset \mathcal{E}{}_{(ab)}[1] $. 

\begin{figure}[tbh]
\begin{tikzcd}
            &  &  & 0 \arrow[d]   &   \\
            & & & {\mathcal{E}_a[1]} \arrow[d] \arrow[ld, "{Z}{}^{a}{}_{A}"', dashed, bend right, shift right] &   \\
0 \arrow[r] & {\mathcal{E}[-1]} \arrow[r, "X_A"] & \mathcal{E}_A \arrow[r] \arrow[ru, "{Z}{}^{A}{}_{a}"', dashed] \arrow[rd, "X^A"] \arrow[l, "Y^A", dashed, bend left, shift left] & {\mathcal{J}^1(\mathcal{E}[1])} \arrow[r] \arrow[d]  & 0 \\
	    &  &  & {\mathcal{E}[1]} \arrow[d] \arrow[lu, "Y_A", dashed, bend left, shift left]   &   \\
            &  &  & 0   &  
\end{tikzcd}
\caption{\label{Fig:T} The tractor projectors}
\end{figure}
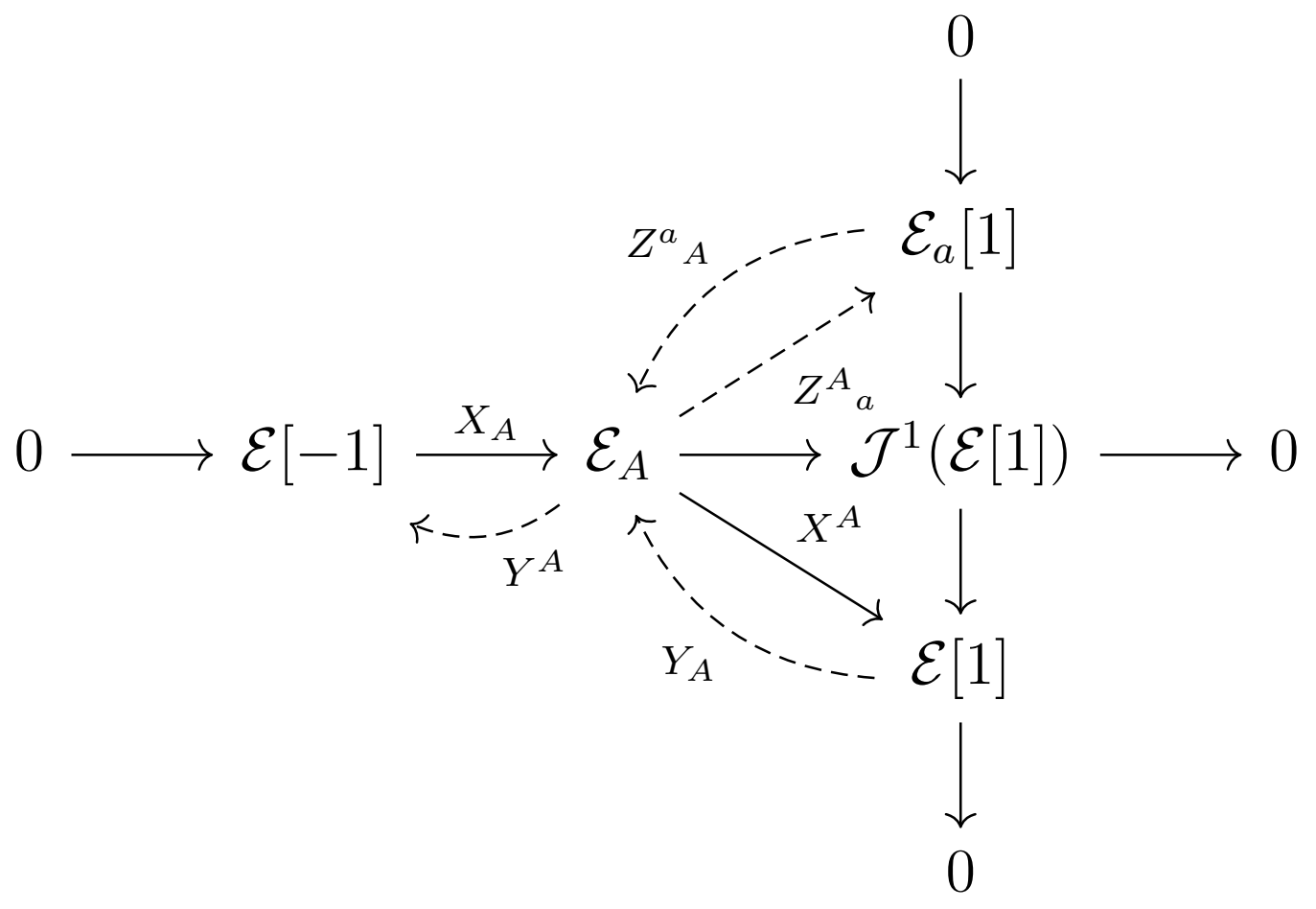

Let $g\in \boldsymbol{c}$ and let $j^2_p(\sigma)\in \ker(A)|_p = \mathcal{E}_A|_p$.
The canonical projector in conformal tractor calculus is $X^A\in \Gamma(\mathcal{E}^A[1])$, given by $X(j^2_p(\sigma)) = \sigma(p)$. The first non-invariant projector is $Z{}^{A}{}_{a} \colon \mathcal{E}_A\to \mathcal{E}_a[1] $ given by $Z(j^2_p(\sigma)) = \nabla_a\sigma$. The last projector is $Y^A \colon \mathcal{E}_A\to \mathcal{E}[-1] $, given by $Y(j^2_p(\sigma)) = -\frac{1}{n}(\Delta \sigma+J\sigma)$. Using the transformation laws for the Levi-Civita connection, we have
\[
\hat{Z}{}^{A}{}_{a} = Z{}^{A}{}_{a} +X^A \Upsilon_a,\quad  \hat{Y}^A = Y^A - Z{}^{A}{}_{a} \Upsilon^a - \tfrac{1}{2} \Upsilon_a\Upsilon^a X^A,
\] 
The conformal tractor bundle is equipped with the invariant \textit{tractor metric},
\[
	h{}^{AB} = X^AY^B+Y^AX^B+ Z{}^{A}{}_{a} Z{}^{B}{}_b\boldsymbol{g}{}^{ab},
\] 
which has signature $(p+1,q+1)$. We use the metric to raise and lower indices, giving rise to the projectors $X_A$, $Z{}_{A}{}^{a} $ and $Y_A$, as illustrated in Figure \ref{Fig:T}, all solid lines are canonical, and the dashed lines depend on a choice of metric $g\in \boldsymbol{c}$.
From the formula for $h{}_{AB} $ follows the tractor contractions, which we will use constantly,
\[
X^AY_A=1,\quad Z{}^{A}{}_{a} Z{}_{A}{}^{b} = \delta{}_{a}{}^{b},
\] 
and all other contractions between projectors are zero.

A metric $g\in \boldsymbol{c}$ determines an isomorphism $\mathcal{E}_A\cong \mathcal{E}[1]\oplus \mathcal{E}_a[1]\oplus \mathcal{E}[-1]$, given by
\[
	U_A = \sigma Y_A + \mu_a Z{}_{A}{}^{a} +\rho X_A \overset{g}{=} (\sigma,\mu_a,\rho),
\]
or as a vertical vector. The invariant \textit{tractor connection} is

\[
\nabla_a\begin{pmatrix} \sigma\\ \mu{}_{a} \\ \rho \end{pmatrix} = \begin{pmatrix} \nabla_a\sigma-\mu_a \\ \nabla_a\mu_b+\sigma\mathsf{P}_{ab}+\rho\boldsymbol{g}{}_{ab}   \\ \nabla_a\rho-\mathsf{P}_a{}^{b} \mu_b\end{pmatrix},
\] 
or in terms of projectors,
\[
\nabla_aX^A = Z^A_a,\quad  \nabla_b Z{}^{A}{}_{a}  = -\mathsf{P}{}_{ab} X^A -\boldsymbol{g}{}_{ab} Y^A\quad \text{and} \quad \nabla_aY^A = \mathsf{P}{}_{a}{}^{b} Z{}^{A}{}_{b}.
\]
The tractor metric is parallel for the tractor connection.

There is a second order conformally invariant differential operator $D\colon \Gamma(\mathcal{E}[1])\to \Gamma(\mathcal{E}_A) $ given by
\begin{equation}\label{Eqn:CTD}
	D_A\sigma|_p =  (\sigma,\nabla_a\sigma,-\tfrac{1}{n} (\Delta\sigma+J\sigma)),
\end{equation}

The operator $D$ gives a one-to-one correspondence between parallel tractors and solutions of $A{}_{ab}  \sigma=0$. Indeed, calculating gives
\[
	\nabla_aD_A \sigma = \begin{pmatrix} \nabla_a \sigma-\nabla_a\sigma\\ \nabla_a\nabla_b\sigma + \mathsf{P}{}_{ab}  \sigma -\frac{1}{n}\left( \Delta \sigma +J \sigma \right)\boldsymbol{g}{}_{ab}   \\  -\frac{1}{n}\left( \nabla_a \Delta \sigma + \nabla_a \left( J\sigma \right)\right) - \mathsf{P}{}_{a}{}^{b} \nabla_b \sigma  \end{pmatrix} = \begin{pmatrix} 0\\ A{}_{ab}  \sigma \\ * \end{pmatrix}. 
\] 
So it is clear that when $\nabla_a D_A\sigma=0$ that $A{}_{ab} \sigma=0$. To show the converse, we suppose that $\sigma$ is a scale, and calculate using the metric $g = \sigma^{-2}\boldsymbol{g}$ on $M\setminus\mathcal{Z}(\sigma)$, giving
\[
 \nabla_a D_A\sigma = \begin{pmatrix} 0\\ {\mathsf{P}^g}{}_{(ab)_0} \sigma \\ -\frac{1}{n}\nabla_a J \sigma \end{pmatrix}.
 \]
 The differential equation $A{}_{ab}  \sigma = 0$ reduces to $\mathsf{P}{}_{(ab)_0} = 0 $ on $M\setminus\mathcal{Z}(\sigma)$, or otherwise stated, $(M\setminus\mathcal{Z}(\sigma), g)$ being Einstein. Then a Lemma of Schur gives $\nabla_a J=0$, so $\nabla_a D{}_{A}\sigma = 0$ on $M\setminus\mathcal{Z}(\sigma)$, and by continuity and density, $D_A \sigma$ is parallel on $M$. 

Let $\sigma\in \Gamma(\mathcal{E}[1])$. The tractor $I_A = D_A \sigma$ is called the \textit{scale tractor} of $\sigma$. The metric $g = \sigma^{-2}\boldsymbol{g}$ does not extend to $\mathcal{Z}(\sigma)$. However, the scale tractor corresponding to $\sigma$ and $g$ is well defined on all of $M$. Using the Levi-Civita connection of $g$ on $M\setminus\mathcal{Z}(\sigma)$, we have $\nabla_a\sigma=0$, so 
\[
D_A\sigma = I_A = \begin{pmatrix} \sigma\\ 0 \\ -\frac{J^g}{n}\sigma \end{pmatrix}\quad \text{and} \quad I_AI_B h{}^{AB}  = |I|^2 = -\tfrac{2}{n}J^g.
\] 

The invariant operator $D$ has an invariant extension to general weighted tractors, the (normalised) \textit{Tractor $D$} operator is $D\colon \Gamma(\mathcal{V}[w])\to \Gamma(\mathcal{V}\otimes\mathcal{E}_A[w-1])$, where $\mathcal{V}$ is any conformal tractor bundle and $w\in \mathbb{R}$ is such that $n+2w-2\neq 0$ is given by
\[
D_A V = w Y_A V + \nabla_a V Z{}_{A}{}^{a} -\frac{1}{n+2w-2}(\Delta V +w JV)X_A,
\] 
where $\nabla{}_{a}$ is the coupled tractor and Levi-Civita connection of some metric from $\boldsymbol{c}$ and $\Delta = \nabla^a\nabla_a$. 

\subsection{Projective to Conformal Tractors}
In this article we will rely on the following result which allows us to identify projective and conformal tractors on an almost Einstein manifold. The following result first appeared in \cite{gover_detecting_2014}, we include the proof in order to be self contained, and to show that the obstruction to such an isomorphism is precisely $\mathsf{P}{}_{(ab)_0} $.

\begin{theorem}\label{Thm:Iperp}
	Let $(M,g)$ be an Einstein manifold and consider the conformal manifold $(M,\boldsymbol{c} = [g])$ and the projective manifold $(M,\boldsymbol{p} = [\nabla^g])$. Let $\sigma$ be the conformal scale density corresponding to $g$ and use $\sigma$ to identify density bundles. The linear bundle map $\phi\colon \mathcal{E}_\alpha\to \mathcal{E}_A $ given by
	\[
		\phi\left(\begin{pmatrix} \tau \\ \mu{}_{a}  \end{pmatrix}\right) = \begin{pmatrix} \tau \\ \mu{}_{a}  \\ \frac{J}{n} \tau  \end{pmatrix} 
	\] 
	is connection preserving. Let $I{}_{A} = {D}_A \sigma$ be the parallel tractor corresponding to the Einstein scale $\sigma$. Then $\phi$ is a connection preserving isomorphism between $\mathcal{E}{}_{\alpha} $ and $I^{\perp}\subset \mathcal{E}_A$.
\end{theorem} 
\begin{proof}
	The bundles $\mathcal{E}{}_{\alpha}$ and $I^{\perp}$ are both rank $n+1$, and $\phi$ is clearly injective, hence a linear bijection. 

	We need to show that if $U{}_{\alpha} \in \Gamma(\mathcal{E}_\alpha)$, then 
	\[
		\nabla_a (\phi{}_{A}{}^{\alpha} U_\alpha) = \phi{}_{A}{}^{\alpha} \nabla_a U_\alpha, 
	\] 
	which is equivalent to showing that $\nabla_a \phi{}_{A}{}^{\alpha} =0$. In the scale $\sigma$ we have
	\[
	\phi{}_{A}{}^{\alpha}  = Y_AX^\alpha+ Z{}_{A}{}^{b} W{}^{\alpha}{}_{b} +\frac{J}{n}X_AX^\alpha.
	\] 	
	\begin{align*}
		\nabla_a\phi{}_{A}{}^{\alpha} 
		&= \mathsf{P}{}_{ab} Z{}_{A}{}^{b} X{}^{\alpha} + Y{}_{A} W{}^{\alpha}{}_{a} -\mathsf{P}{}_{a}{}^{b} X{}_{A} W{}^{\alpha}{}_{b} -\boldsymbol{g}{}_{a}{}^{b} Y{}_{A} W{}^{\alpha}{}_{b} -\tilde{\mathsf{P}}{}_{ab} Z{}_{A}{}^{b} X{}^{\alpha} \\ 
		&\mathrel{\phantom{=}}+\frac{\nabla_a J}{n}X_A X^\alpha +\frac{J}{n}Z{}_{Aa} X^\alpha + \frac{J}{n}X_A W{}^{\alpha}{}_{a} \\ 
		&= -\overset{\circ }{\mathsf{P}}{}_{a}{}^{b} X{}_{A} W{}^{\alpha}{}_{b} +\frac{\nabla_a J}{n}X_AX^\alpha + \frac{1}{n-1}\overset{\circ }{\mathsf{P}}{}_{ab} Z{}_{A}{}^{b} X{}^{\alpha}.
	\end{align*}
	For an Einstein metric $\overset{\circ }{\mathsf{P}}=0$ by definition and $\nabla_aJ=0$ by Schur's Lemma. 
	Since $I_A = {D}_A\sigma = (\sigma,\nabla_a\sigma,-\frac{J}{n}\sigma)$, clearly $\phi$ maps into $I^{\perp}$, which is a parallel subbundle.
\end{proof}

The above calculation of $\nabla_a \phi{}_{A}{}^{\alpha} $ holds for any pseudo-Riemannian manifold and shows that the obstruction to $\phi$ being connection preserving is measured precisely by $\overset{\circ }{\mathsf{P}} $ and $\nabla_a J$, the obstructions to $g$ being Einstein.

The above theorem has an interpretation in terms of prolongation. The conformal tractor bundle is the prolongation bundle for the almost Einstein operator $A$, whereas the projective tractor bundle is the prolongation bundle for the projective equivalent $\tilde{A}$. 

It turns out that one can construct the projective tractor bundle analogously to the construction of conformal tractors, take $\tilde{\mathcal{T}}^* = \ker(\tilde{A})\subset J^2(\mathcal{E}(1))$. In this case, $\tilde{}A\colon J^2(\mathcal{E}(1))\to \mathcal{E}{}_{(ab)} (1) $ is a genuine splitting of the second jet exact sequence for $\mathcal{E}(1)$, so we find that $\ker(\tilde{A})\cong J^2(\mathcal{E}(1)) / \mathcal{E}{}_{(ab)} (1)\cong J^1(\mathcal{E}(1))$. 

In the Einstein case, the almost Einstein equation reduces to 
\[
A{}_{ab} \tau = \nabla{}_{(a} \nabla{}_{b)_0} \tau=0.
\] 
The projective equation also simplifies since $\tilde{\mathsf{P}}{}_{ab}  = \frac{1}{n-1}\operatorname{Ric}{}_{ab} = \frac{1}{n(n-1)}\operatorname{Sc}\boldsymbol{g}{}_{ab} $,
\[
\tilde{A}{}_{ab}  \tau = \nabla{}_{(a}\nabla{}_{b)} \tau+\frac{1}{n(n-1)}\operatorname{Sc}\tau \boldsymbol{g}{}_{ab} =0.
\] 
Then clearly a solution of $\tilde{A}{}_{ab} \tau=0$ will also satisfy $A{}_{ab} \tau=0$.

In fact, it is easy to show that the inclusion $\phi$ above is simply the inclusion of the subbundles $\ker(\tilde{A})\subset \ker(A)\subset J^2(\mathcal{E}[1])$, trivialising density bundles using $\sigma$.

For the PDE $\tilde{A}\tau=0$, the above Theorem characterises the solutions.

 \subsection{BGG Sequences}\label{Sec:BGG}

Bernstein-Gelfand-Gelfand (BGG) sequences in parabolic geometry \cite{cap_bernstein-gelfand-gelfand_2001, calderbank_differential_2001} are invariant sequences of differential operators associated to any irreducible tractor bundle.
The tractor bundles we consider in conformal geometry, are $\mathfrak{g} = \mathfrak{so}(p+1,q+1)$ invariant subbundles of a tensor power of the standard tractor bundle, equipped with the signature $(p+1,q+1)$ bilinear form $h{}_{AB}$. We will exclude spin representations throughout. In particular the irreducible tractor bundles we consider are the trace-free tractors with symmetries given by a Young diagram.

The tractor bundles we consider in projective geometry are $\mathfrak{sl}(n+1)$ invariant subbundles of a tensor product of the standard tractor and cotractor bundles. The irreducible projective tractor bundles we consider have Young diagram symmetries, one diagram for upper and one for lower indices, and are totally trace free, for contractions between upper and lower indices.

Let $\mathcal{V}$ be an irreducible (conformal or projective) tractor bundle.
Associated to $\mathcal{V}$, equipped with its tractor connection $\nabla$ is the \textit{twisted de Rham} sequence,
\[
	0\longrightarrow \mathcal{V} \overset{d^\nabla}{\longrightarrow} T^*M\otimes \mathcal{V} \overset{d^\nabla}{\longrightarrow} \Lambda^2T^*M\otimes\mathcal{V} \overset{d^\nabla}{\longrightarrow} \cdots,
\] 
where $d^{\nabla}(V)_{a_{0}a_{1}\cdots a_k} = (k+1)\nabla{}_{[a_{0}} V_{a_{1}\cdots a_k]}$, we have suppressed the tractor indices of $\mathcal{V}$ and $\nabla$ is the tractor connection paired with any torsion free affine connection on $M$. When the tractor connection is flat, the twisted de Rham sequence is a complex. 

The Kostant codifferential is an algebraic operator acting in the opposite direction,
\[
0\overset{\partial^*}{\longleftarrow} \mathcal{V} \overset{\partial^*}{\longleftarrow} T^*M\otimes \mathcal{V} \overset{\partial^*}{\longleftarrow} \Lambda^2T^*M\otimes\mathcal{V} \overset{\partial^*}{\longleftarrow} \cdots,
\] 
which satisfies $\partial^*\circ \partial^*=0$.
For conformal tractor bundles $\mathcal{V}\subset \mathcal{E}_{A_{1}\cdots A_j}$ and $V\in \Gamma(\Lambda^kT^*M\otimes\mathcal{V})$, we let
\[
	\partial^*(V)_{a_{1}\cdots a{}_{k-1} \; A_{1}\cdots A_j} = \sum_{i=1}^{j} 2\mathbb{X}_{A_i}{}^{B a_k}V_{a_{1}\cdots a_k A_{1}\cdots A_{i-1} B A_{i+1}\cdots A_j},
\]
where $\mathbb{X}{}_{AB}{}^{b}  = X{}_{[A} Z{}_{B]}{}^{b} $. The tractor $\mathbb{X}{}_{AB}{}^{b} $ is unweighted and invariant since $\hat{\mathbb{X}}{}_{AB}{}^{b} = \mathbb{X}{}_{AB}{}^{b} +X{}_{[A} X{}_{B]} \Upsilon^b = \mathbb{X}_{AB}{}^b$ and provides an invariant bundle inclusion $\mathbb{X}\colon \mathcal{E}_b\hookrightarrow \mathcal{E}_{[AB]}$. For projective tractor bundles $\mathcal{V}\subset \mathcal{E}_{\alpha_{1}\cdots \alpha_j}{}^{\gamma_{1}\cdots \gamma_\ell}$ and $V\in \Gamma(\Lambda^k T^*M\otimes \mathcal{V})$, we let
\[
	\begin{aligned}\partial^*(V)_{a_{1}\cdots a{}_{k-1} \; \alpha_{1}\cdots \alpha_j}^{\gamma_{1}\cdots \gamma_\ell} 
&= \sum_{i=1}^{j} 2\mathbb{X}^{\beta}{}_{\alpha_i}{}^{a_k}V_{a_{1}\cdots a_k \alpha_{1}\cdots \alpha_{i-1} \beta \alpha_{i+1}\cdots \alpha_j}{}^{\gamma_{1}\cdots \gamma_\ell} \\
&\mathrel{\phantom{=}}- \sum_{i=1}^{\ell} 2\mathbb{X}^{\gamma_i}{}_{\beta}{}^{a_k}V_{a_{1}\cdots a_k\alpha_{1}\cdots \alpha_j}{}^{\gamma_{1}\cdots \gamma_{i-1} \beta \gamma_{i+1}\cdots \gamma_\ell} \end{aligned}
\]
where $\mathbb{X}^{\beta}{}_{\alpha}{}^{b}  = X^{\beta} Z{}_{\alpha}{}^{b} $ is unweighted invariant since both $X{}^{\beta} $ and $Z{}_{\alpha}{}^{b} $ are canonical. Similarly to the conformal case, $\mathbb{X}{}^{\beta}{}_{\alpha}{}^{b} $ provides an invariant bundle inclusion, $\mathbb{X}\colon \mathcal{E}_b\hookrightarrow \operatorname{trace-free}(\mathcal{E}{}^{\beta}{}_{\alpha}) $.

The homology bundles are 
\[
	H_k(\mathcal{V}) = \operatorname{ker}(\partial^*\colon \Lambda^k T^*M \otimes \mathcal{V}\to \Lambda^{k-1}T^*M\otimes\mathcal{V} )/\operatorname{im}(\partial^*\colon \Lambda^{k+1}T^*M\otimes \mathcal{V}\to \Lambda^k T^*M\otimes\mathcal{V} ).
\]
We denote the projection onto homology by $\Pi\colon  \operatorname{ker}(\partial^*)\to H_k(\mathcal{V})$ and note that for $k=0$, this gives a map $\Pi\colon \mathcal{V}\to H_0(\mathcal{V}) $. For irreducible parabolic geometries, including conformal and projective, the homology bundles $H_k(\mathcal{V})$ are irreducible $\mathfrak{g}_0$ bundles, except in the even dimensional conformal case for $k=\frac{n}{2}$. For conformal geometry, $\mathfrak{g}_0 = \mathfrak{co}(p,q) = \mathfrak{so}(p,q)\oplus \mathbb{R}$, so the homology bundles are trace-free irreducible weighted tensor bundles. In projective geometry, $\mathfrak{g}_{0} = \mathfrak{gl}(n)$, so the homology bundles are weighted irreducible tensor bundles.

The Bernstein-Gelfand-Gelfand machinery gives two important families of operators. The first are the BGG splitting operators. 
\begin{theorem}[Theorem 2.5 and Lemma 2.7 in \cite{cap_bernstein-gelfand-gelfand_2001}]\label{Thm:BGGS}
	There exist invariant splitting operators, $L\colon \Gamma(H_k(\mathcal{V}))\to \Gamma(\operatorname{ker}(\partial^*))$ such that $\Pi\circ L = \operatorname{id}_{H_k(\mathcal{V})}$ and $\partial^*\circ d^{\nabla}\circ L=0$. The subspace $\operatorname{im}(L)\subset \Gamma(\Lambda^k(T^*M)\otimes\mathcal{V})$ consists of sections $s$ such that $\partial^*(s)=0$ and $\partial^*\circ d^{\nabla}(s)=0$.
\end{theorem}
We call $L{}_{k}(s) $ the \textit{full prolongation} of $s$.
The above conditions are sufficient to calculate the splitting operators. Take the $\mathcal{V} = \mathcal{E}_A$, the standard conformal cotractor bundle. We can compute the action of $\partial^*$ on $V{}_{aA} = \omega_aY_A+\nu_{ab}Z{}_{A}{}^{b} +\tau_aX_A \in \Gamma(\mathcal{E}{}_{aA} )$ in a splitting,
\begin{equation}\label{Eqn:dsT}
	\partial^*(V) = 2\mathbb{X}{}_{A}{}^{Ba} (\omega_aY_B+\nu_{ab}Z{}_{B}{}^{b} +\tau_aX_B) = -\omega{}_{b} Z{}_{A}{}^{b}  +\nu{}_{a}{}^{a} X_A.
\end{equation} 
So suppose $\sigma\in \Gamma(H_0(\mathcal{E}_A)) = \Gamma(\mathcal{E}[1])$ and that $L(\sigma)_A = (\sigma,\mu_a,\rho)$. Then $d^{\nabla}\circ L(\sigma) = \nabla_a(\sigma,\mu_b,\rho) = (\nabla_a\sigma-\mu_a,\nabla_a\mu_b+ \sigma\mathsf{P}{}_{ab} +\rho\boldsymbol{g}{}_{ab} , \nabla_a\rho- \mathsf{P}{}_{a}{}^{b} \mu_b) $. So the condition $\partial^*\circ d^{\nabla}\circ L=0$ becomes $\nabla_a\sigma-\mu_a = 0$ and $\nabla_a\mu^a+J\sigma+n\rho=0$. Hence, $\mu_a = \nabla_a\sigma$ and $\rho=-\frac{\Delta+J}{n}\sigma$. 
So $L\colon \Gamma(\mathcal{E}[1])\to \Gamma(\mathcal{E}_A) $ is given by $L(\sigma) = D_A\sigma$, as in \eqref{Eqn:CTD}. Similarly, the splitting operator $L\colon \Gamma(\mathcal{E}(1))\to \Gamma(\mathcal{E}{}_{\alpha} ) $ is the projective tractor $D$.

The splitting operators give rise to the second family of operators, the BGG differential operators $\mathscr{D}_k:= \Pi\circ d^{\nabla}\circ L \colon \Gamma(H_k(\mathcal{V})) \to \Gamma(H_{k+1}(\mathcal{V}))$ and the BGG sequence,
\[
	0\longrightarrow \Gamma(H_0(\mathcal{V}))\overset{\mathscr{D}_0}{\longrightarrow} \Gamma(H_1(\mathcal{V}))\overset{\mathscr{D}_1}{\longrightarrow} \Gamma(H_2(\mathcal{V}))\overset{\mathscr{D}_2}{\longrightarrow}\cdots 
\]
which in the flat case is a complex. If $U$ is parallel, then $\partial^*(U)=0$ and $\partial^*\circ d^{\nabla}(U)=0$, so $U = L_0\circ \Pi(U)$ and $\mathscr{D}_0\circ \Pi(U) = \Pi\circ d^{\nabla}(U)=0$. A solution of $\mathscr{D}_0s = 0$ is \textit{normal} if $L{}_{0}(s) $ is parallel. 

For flat geometries, $d^{\nabla}\circ L_k = L_{k+1}\circ \mathscr{D}_k$, so the squares in the following diagram commute
\begin{equation}\label{Eqn:BGG}
\begin{array}{ccccccccccc}
	0
&\longrightarrow
&\Gamma(\mathcal{V})
& \xrightarrow{\; d^\nabla \;}
& \Gamma\!\big(\Lambda^1T^*M\otimes \mathcal{V}\big)
& \xrightarrow{\; d^\nabla \;}
& \Gamma\!\big(\Lambda^2T^*M\otimes \mathcal{V}\big)
& \xrightarrow{\; d^\nabla \;}
& \cdots
\\[0.2em]
&&\big\uparrow\, L_0
&
&
\big\uparrow\, L_1
&
&
\big\uparrow\, L_2
&
\\[0.2em]
0
&\longrightarrow
&\Gamma\!\big({H}_0(\mathcal{V})\big)
& \xrightarrow{\; \mathscr{D}_0 \;}
& \Gamma\!\big({H}_1(\mathcal{V})\big)
& \xrightarrow{\; \mathscr{D}_1 \;}
& \Gamma\!\big({H}_2(\mathcal{V})\big)
& \xrightarrow{\; \mathscr{D}_2 \;}
& \cdots
\end{array}
\end{equation}
However, for general curved geometries, the above diagram does not commute.
\begin{theorem}\label{Thm:BGGL}
	Let $\mathcal{V}$ be an irreducible tractor bundle on a conformally flat or projectively flat manifold. Then the diagram \eqref{Eqn:BGG} commutes. 

	The map $L_k$ gives a one-to-one correspondence between tensors in the kernel of the $k$-th BGG operator $\mathscr{D}_k \colon \Gamma(H_k(\mathcal{V}))\to \Gamma(H_{k+1}(\mathcal{V}))$ and sections of $\Lambda^k(T^*M)\otimes\mathcal{V}$ which are in the kernel of $\partial^*$ and $d^{\nabla}$. 

	In particular, for $k=0$ there is a one-to-one correspondence between solutions of $\mathscr{D}_0 s=0$ and parallel tractors in $\mathcal{V}$. So all solutions of $\mathscr{D}_0s=0$ are normal.
\end{theorem}
\begin{proof}
	In the flat case we have $d^{\nabla}\circ d^{\nabla}=0$, then $d^{\nabla}\circ L_k(s)\in \ker(\partial^*)\cap \ker(\partial^*\circ d^{\nabla})$, so $d^{\nabla}\circ L_k(s)\in \operatorname{im}(L_{k+1})$ by Theorem \ref{Thm:BGGS}. Then observe that $\Pi\circ d^{\nabla}\circ L_k(s) = \mathscr{D}_k(s) = \Pi\circ L_{k+1}\circ \mathscr{D}_k(s)$, and $d^{\nabla}\circ L_k(s) = L_{k+1}\circ \mathscr{D}_k(s)$ now follows since $\Pi$ is injective on $\operatorname{im}(L_{k+1})$.

	Let $s\in\Gamma(H_k(\mathcal{V}))$. Then $L_k(s)\in \ker(\partial^*)$, and if $\mathscr{D}_k(s)=0$, then $L_{k+1}\circ \mathscr{D}_k(s) = d^{\nabla}\circ L{}_{k} (s) = 0$, so $L_k(s)$ is in the kernel of both $\partial^*$ and $d^{\nabla}$. 

	The converse direction doesn't require the geometry to be flat. Let $V\in \Gamma(\Lambda^kT^*M\otimes\mathcal{V})$ be in the kernel of both $\partial^*$ and $d^{\nabla}$, then $V\in \operatorname{im}(L_k)$ by Theorem \ref{Thm:BGGS}. So $V = L_k(s)$ for some $s\in\Gamma(H_k(\mathcal{V}))$ and $\mathscr{D}_k(s) = \Pi\circ d^{\nabla}(V) = 0$.
\end{proof}

For $\mathcal{V} = \mathcal{E}{}_{A} $, the image of $\partial^*$ can be determined from \eqref{Eqn:dsT} to be the subspace spanned by $Z$ and $X$. So the projection onto homology is given by $\sigma Y_A+\mu_a Z{}_{A}{}^{a} +\rho X_A\mapsto \sigma$, equivalently $U_A\mapsto X^AU_A$.

 In conformal geometry for an irreducible tractor bundle with symmetries given by a Young diagram, the projection $\Pi\colon \mathcal{V}\to H_0(\mathcal{V}) $ is given by first applying as many $X^A$ projectors as are allowed by the symmetries and then applying $Z{}^{A}{}_{a} $ projectors to the remaining indices. For projective tractor bundles, the projection on lower indices is given by first applying as many $X{}^{\alpha}$ projectors as are allowed by the symmetries and then applying $W{}^{\alpha}{}_{a} $ projectors to the remaining indices. For upper indices, the invariant $Z{}_{\alpha}{}^{a} $ projector is applied to all indices.

 For a tractor $V\in \Gamma(\mathcal{V})$, the tensor $\Pi(V)$ is called the \textit{top slot} of $V$. The projection is always invariant. The following Lemma is frequently useful to move from tensor to tractor equations.

\begin{lemma}\label{Lem:BGG}
	Let $\mathcal{V}$ be an irreducible tractor bundle and let $V\in\Gamma(\mathcal{V})$ be parallel. If $\Pi(V)=0$, then $V=0$.
\end{lemma}
\begin{proof}
	By Theorem \ref{Thm:BGGS}, $V = L\circ \Pi(V) = 0$, since $d^{\nabla}(V) = 0$ and $ \partial^*(V) = 0$.
\end{proof}
For any particular tractor bundle, one can prove the above lemma by brute force by differentiating and showing that all slots of $V$ vanish. For example, if $V{}_{A}\in \Gamma(\mathcal{E}_A)$ is parallel and $X^AV_A=0$, then, in any scale,
\[
	\nabla_a(X^AV_A) = Z{}^{A}{}_{a}  V_A = 0 \implies \nabla_b( Z{}^{A}{}_{a} V_A) = -\mathsf{P}{}_{ab}  X^AV_A-\boldsymbol{g}_{ab} Y^AV_A = -\boldsymbol{g}_{ab} Y^AV_A=0,
\] 
so $X^AV_A = 0$, $Z{}^{A}{}_{a} V_A=0$ and $Y^AV_A=0$, so $V=0$. The calculations get more tedious as the number of tractor indices increases, making the BGG machinery a convenient route to take. 

The tractor bundles and BGG operators we will use in this article are detailed below.

The projective and conformal Killing vector operators are the BGG operators for $\mathcal{E}{}_{[\alpha\beta]} $ and $\mathcal{E}{}_{[AB]} $, respectively. For the bundle $\mathcal{E}{}_{[AB]} $ the BGG sequence is, when $n\ge 4$,
\[
	\begin{tikzcd}[column sep = 1.8em]
		\Gamma({\mathcal{E}_a[2]}) & \Gamma({\mathcal{E}_{(ab)_0}[2]}) & {\Gamma\left(\begin{picture}(22, 11.0)(0, 5.0)
	\put(-0.2, 16){\line(1,0){16.4}}
	\put(-0.2, 8){\line(1,0){16.4}}
	\put(-0.2, 0){\line(1,0){16.4}}
	\put(0, 16){\line(0,-1){16}}
	\put(8, 16){\line(0,-1){16}}
	\put(16, 16){\line(0,-1){16}}
	\put(16, 0){${}_{\,0}$} 
\end{picture}(\mathcal{E}_a)[2]\right)} & {\cdots},
	\arrow["{\mathscr{D}_0}", from=1-1, to=1-2]
	\arrow["\mathscr{D}_1", from=1-2, to=1-3]
	\arrow["\mathscr{D}_2", from=1-3, to=1-4]
\end{tikzcd}
\] 
where the operator $\mathscr{D}_0$ is the conformally invariant \textit{conformal Killing operator}
\[
k_a\mapsto \nabla{}_{(a} k{}_{b)_0}.
\]
A tensor $k_a\in \Gamma(\mathcal{E}_a[2])$ in the kernel of the above conformal Killing operator is a \textit{conformal Killing vector}.
In the projective case, the BGG sequence for $\mathcal{E}{}_{[\alpha\beta]} $ is
\[
\begin{tikzcd}[column sep = 1.8em]
	\Gamma(\mathcal{E}_a(2))\arrow[r, "\mathscr{D}_0"] &\Gamma(\mathcal{E}{}_{(ab)}(2)) \arrow[r, "\mathscr{D}_1"] & \Gamma\left({\begin{picture}(18, 11.0)(0, 5.0)
	\put(-0.2, 16){\line(1,0){16.4}}
	\put(-0.2, 8){\line(1,0){16.4}}
	\put(-0.2, 0){\line(1,0){16.4}}
	\put(0, 16){\line(0,-1){16}}
	\put(8, 16){\line(0,-1){16}}
	\put(16, 16){\line(0,-1){16}}
\end{picture}(\mathcal{E}_a)(2)}\right)\arrow[r, "\mathscr{D}_2"]&\cdots,
\end{tikzcd}
\] 
where the operator $\mathscr{D}_0$ is the projectively invariant \textit{projective Killing operator},
\[
r_a\mapsto \nabla{}_{(a} r{}_{b)}.
\] 
A tensor $r_a\in \Gamma(\mathcal{E}_a(2))$ in the kernel of the above operator is a \textit{projective Killing vector}. 

When the projective structure is metric projective, that is, when $\boldsymbol{p} = [\nabla^g]$ for a pseudo-Riemannian metric $g$, then a projective Killing vector for $(M,\boldsymbol{p})$ is the same as a Killing vector for $(M,g)$, once we use the volume form of $g$ to trivialise the density bundles. After removing the weights, the above sequence becomes the Riemannian deformation sequence, or Kr\"oner sequence \cite{kroner_dislocations_1975}, the operator $\mathscr{D}_1$ is often referred to as the incompatibility operator, and is important to the study of elasticity \cite{arnold_mixed_2007}.

There are symmetric and skew generalisations of Killing vectors. Let $1\le \ell\le n-1$. Consider the completely skew tractor bundles $\mathcal{E}_{[A_{0}A_{1}\cdots A_\ell]}$ and $\mathcal{E}_{[\alpha_{0}\alpha_{1}\cdots \alpha_\ell]}$. The homology bundles are $H_0(\mathcal{E}_{[A_{0}A_{1}\cdots A_\ell]}) = \mathcal{E}_{[a_{1}\cdots a_\ell]}[\ell+1]$ and $H_0(\mathcal{E}_{[\alpha_{0}\alpha_{1}\cdots \alpha_\ell]}) = \mathcal{E}_{[a_{1}\cdots a_\ell]}(\ell+1)$. These homology bundles can be determined by first applying one $X$ projector and observing that any further $X$ projectors will result in zero due to skew indices. Therefore the remaining projectors must be $Z{}^{A}{}_{a} $ in the conformal case and $W{}^{\alpha}{}_{a} $ in the projective case.

If $n$ is even and $\frac{n}{2}+q$ is even, then taking $\ell = \frac{n}{2}$ results in the conformal tractor bundle $\mathcal{E}_{[A_{0}A_{1}\cdots A_\ell]}$ and the tensor bundle $\mathcal{E}_{[a_{1}\cdots a_\ell]}[\ell+1]$ being of middle degree, so they are not irreducible under the action of $\mathfrak{so}(p+1,q+1)$ and $\mathfrak{so}(p,q)$. While this is of independent interest, we will not address this here. The BGG operator for middle conformal Killing--Yano forms in this case is the sum of the BGG operator on the self-dual anti-self-dual forms. For the projective operators, $\ell = \frac{n}{2}$ is not a special case.

The conformal BGG sequence is
\[
\begin{tikzcd}[column sep = 1.8em]
	\Gamma(\mathcal{E}_{[a_{1}\cdots a_\ell]}[\ell+1])
\arrow[r, "\mathscr{D}_0"] &\Gamma\left(\begin{picture}(16, 13.0)(0, 7.0)
	\put(-0.2, 20){\line(1,0){10.4}}
	\put(-0.2, 15){\line(1,0){10.4}}
	\put(1.5, 12.5){${}_{\cdot}$}
	\put(1.5, 10){${}_{\cdot}$}
	\put(1.5, 7.5){${}_{\cdot}$}
	\put(-0.2, 5){\line(1,0){5.4}}
	\put(-0.2, 0){\line(1,0){5.4}}
	\put(0, 20){\line(0,-1){20}}
	\put(5, 20){\line(0,-1){20}}
	\put(10, 20){\line(0,-1){5}}
	\put(10, 0){${}_{\,0}$}
\end{picture}(\mathcal{E}_a)[\ell+1]\right) \arrow[r, "\mathscr{D}_1"] & \cdots,
\end{tikzcd}
\] 
where the operator $\mathscr{D}_0$ is the conformally invariant \textit{conformal Killing--Yano operator},
\[
	\omega_{a_{1}\cdots a_\ell}\mapsto \operatorname{trace-free}(\nabla_{(a_{0}}\omega_{a_{1})a_{2}\cdots a_\ell}).
\] 
The kernel of the above operator consists of the \textit{conformal Killing--Yano forms}.
The projective BGG sequence is
\[
\begin{tikzcd}[column sep = 1.8em]
	\Gamma(\mathcal{E}_{[a_{1}\cdots a_\ell]}(\ell+1))
\arrow[r, "\mathscr{D}_0"] &\Gamma\left(\begin{picture}(16, 13.0)(0, 7.0)
	\put(-0.2, 20){\line(1,0){10.4}}
	\put(-0.2, 15){\line(1,0){10.4}}
	\put(1.5, 12.5){${}_{\cdot}$}
	\put(1.5, 10){${}_{\cdot}$}
	\put(1.5, 7.5){${}_{\cdot}$}
	\put(-0.2, 5){\line(1,0){5.4}}
	\put(-0.2, 0){\line(1,0){5.4}}
	\put(0, 20){\line(0,-1){20}}
	\put(5, 20){\line(0,-1){20}}
	\put(10, 20){\line(0,-1){5}}
\end{picture}(\mathcal{E}_a)(\ell+1)\right) \arrow[r, "\mathscr{D}_1"] & \cdots,
\end{tikzcd}
\] 
where the operator $\mathscr{D}_0$ is the projectively invariant \textit{projective Killing--Yano operator},
\[
	\omega_{a_{1}\cdots a_\ell}\mapsto \nabla_{(a_{0}}\omega_{a_{1})a_{2}\cdots a_\ell}.
\] 
The kernel of the above operator consists of the \textit{projective Killing--Yano forms}. In the metric projective case, we just refer to \textit{Killing--Yano forms} for $(M,g)$.

The symmetric generalisation of Killing vectors are the Killing tensors. Let $\ell\ge 1$ and consider the tractor bundles \begin{picture}(31, 8.0)(0, 2.0)
	\put(-0.2, 10){\line(1,0){25.4}}
	\put(-0.2, 5){\line(1,0){25.4}}
	\put(-0.2, 0){\line(1,0){25.4}}
	\put(0, 10){\line(0,-1){10}}
	\put(5, 10){\line(0,-1){10}}
	\put(9, 2.3){${}_{\cdots}$}
	\put(9, 7.3){${}_{\cdots}$}
	\put(20, 10){\line(0,-1){10}}
	\put(25, 10){\line(0,-1){10}}
	\put(25, 0){${}_{\,0}$}
\end{picture}$(\mathcal{E}_A)$ and \begin{picture}(31, 8.0)(0, 2.0)
	\put(-0.2, 10){\line(1,0){25.4}}
	\put(-0.2, 5){\line(1,0){25.4}}
	\put(-0.2, 0){\line(1,0){25.4}}
	\put(0, 10){\line(0,-1){10}}
	\put(5, 10){\line(0,-1){10}}
	\put(9, 2.3){${}_{\cdots}$}
	\put(9, 7.3){${}_{\cdots}$}
	\put(20, 10){\line(0,-1){10}}
	\put(25, 10){\line(0,-1){10}}
	\end{picture}$(\mathcal{E}_\alpha)$. We use the skew representation of these spaces. So tractors $K_{A_{1}B_{1}\cdots A_\ell B_\ell}\in \begin{picture}(31, 8.0)(0, 2.0)
	\put(-0.2, 10){\line(1,0){25.4}}
	\put(-0.2, 5){\line(1,0){25.4}}
	\put(-0.2, 0){\line(1,0){25.4}}
	\put(0, 10){\line(0,-1){10}}
	\put(5, 10){\line(0,-1){10}}
	\put(9, 2.3){${}_{\cdots}$}
	\put(9, 7.3){${}_{\cdots}$}
	\put(20, 10){\line(0,-1){10}}
	\put(25, 10){\line(0,-1){10}}
	\put(25, 0){${}_{\,0}$}
\end{picture}(\mathcal{E}_A)$ satisfy
\[
	K{}_{A_{1}B_{1}\cdots A_\ell B_\ell} = K{}_{[A_{1}B_{1}]\cdots [A_\ell B_\ell]}\quad \text{and} \quad K{}_{\cdots [A_{i}B_{i}A_{i+1}]B_{i+1}\cdots}=0 
\] 
and are completely trace-free. The projective tractor bundle consists of tractors with the same symmetries without the trace-free condition. To determine the homology bundles, observe that we can apply $X^{A_i}$ (or $X^{\alpha_i}$) for all $1\le i\le \ell$ without the resulting weighted tractor vanishing. However, the remaining projectors must all be $Z^{B_i}{}_{a_i}$ (or $W^{\beta_i}{}_{a_i}$) due to the skew index pairs.

The resulting homology bundles are $\mathcal{E}{}_{(a_{1}\cdots a_\ell)_0}[2\ell] $ and $\mathcal{E}{}_{(a_{1}\cdots a_\ell)}(2\ell) $. The assertion that the conformal case produces a trace-free tensor bundle needs justification. Let $K_{A_{1}B_{1}\cdots A_\ell B_\ell}\in \Gamma\left(\begin{picture}(31, 8.0)(0, 2.0)
	\put(-0.2, 10){\line(1,0){25.4}}
	\put(-0.2, 5){\line(1,0){25.4}}
	\put(-0.2, 0){\line(1,0){25.4}}
	\put(0, 10){\line(0,-1){10}}
	\put(5, 10){\line(0,-1){10}}
	\put(9, 2.3){${}_{\cdots}$}
	\put(9, 7.3){${}_{\cdots}$}
	\put(20, 10){\line(0,-1){10}}
	\put(25, 10){\line(0,-1){10}}
	\put(25, 0){${}_{\,0}$}
\end{picture}(\mathcal{E}_A)\right)$. Then
\[
	\begin{aligned}
		\boldsymbol{g}^{a_{1}a_2}(X^{A_1}Z^{B_1}{}_{a_1}X^{A_2}Z^{B_2}{}_{a_2}\cdots) K_{A_{1}B_{1}\cdots A_\ell B_\ell} = 
		&-(X^{A_1}X^{B_1}X^{A_2}Y^{B_2}\cdots) K_{A_{1}B_{1}\cdots A_\ell B_\ell} \\
		&- (X^{A_1}Y^{B_1}X^{A_2}X^{B_2}\cdots) K_{A_{1}B_{1}\cdots A_\ell B_\ell} = 0,
	\end{aligned}
\] 
since $K$ is trace-free and skew on all index pairs $A_iB_i$.

The BGG sequence for $ \begin{picture}(31, 8.0)(0, 2.0)
	\put(-0.2, 10){\line(1,0){25.4}}
	\put(-0.2, 5){\line(1,0){25.4}}
	\put(-0.2, 0){\line(1,0){25.4}}
	\put(0, 10){\line(0,-1){10}}
	\put(5, 10){\line(0,-1){10}}
	\put(9, 2.3){${}_{\cdots}$}
	\put(9, 7.3){${}_{\cdots}$}
	\put(20, 10){\line(0,-1){10}}
	\put(25, 10){\line(0,-1){10}}
	\put(25, 0){${}_{\,0}$}
\end{picture}(\mathcal{E}_A)$ is
\[
\begin{tikzcd}[column sep = 1.8em]
	\Gamma(\mathcal{E}_{(a_{1}\cdots a_\ell)_0}[2\ell])
	\arrow[r, "\mathscr{D}_0"] &\Gamma\left(\mathcal{E}_{(a_{0}a_{1}\cdots a_\ell)_0}[2\ell]\right) \arrow[r, "\mathscr{D}_1"] & \cdots,
\end{tikzcd}
\]
where the operator $\mathscr{D}_0$ is the conformally invariant \textit{conformal Killing tensor operator}
\[
k{}_{a_{1}\cdots a_\ell}\mapsto \nabla{}_{(a_{0}}k{}_{a_{1}\cdots a_\ell)_0}.
\] 
Tensors in the kernel of the above operator are \textit{rank $\ell$ conformal Killing tensors}.
The BGG sequence for $ \begin{picture}(31, 8.0)(0, 2.0)
	\put(-0.2, 10){\line(1,0){25.4}}
	\put(-0.2, 5){\line(1,0){25.4}}
	\put(-0.2, 0){\line(1,0){25.4}}
	\put(0, 10){\line(0,-1){10}}
	\put(5, 10){\line(0,-1){10}}
	\put(9, 2.3){${}_{\cdots}$}
	\put(9, 7.3){${}_{\cdots}$}
	\put(20, 10){\line(0,-1){10}}
	\put(25, 10){\line(0,-1){10}}
	\put(25, 0){${}_{\,0}$}
\end{picture}(\mathcal{E}_\alpha)$ is
\[
\begin{tikzcd}[column sep = 1.8em]
	\Gamma(\mathcal{E}_{(a_{1}\cdots a_\ell)}(2\ell))
	\arrow[r, "\mathscr{D}_0"] &\Gamma\left(\mathcal{E}_{(a_{0}a_{1}\cdots a_\ell)}(2\ell)\right) \arrow[r, "\mathscr{D}_1"] & \cdots,
\end{tikzcd}
\]
where the operator $\mathscr{D}_0$ is the projectively invariant \textit{projective Killing tensor operator}
\[
r{}_{a_{1}\cdots a_\ell}\mapsto \nabla{}_{(a_{0}}r{}_{a_{1}\cdots a_\ell)}.
\] 
Tensors in the kernel of the above operator are \textit{rank $\ell$ projective Killing tensors}. In the metric projective case, we instead refer to \textit{rank $\ell$ Killing tensors} for $(M,g)$.

The above examples share one important quality. That the conformal BGG operator is a \textit{weakening} of the projective BGG operator. Precisely, let $g$ be a metric on $M$, let $(M,\boldsymbol{c}=[g])$ and $(M,\boldsymbol{p}=[\nabla^g])$ use the scale density $\sigma\in \Gamma(\mathcal{E}[1])$ corresponding to $g$ to trivialise density bundles. 

Then a projective Killing--Yano form is a conformal Killing Yano form and the trace-free part of a Killing tensor is a conformal Killing tensor. Given a fixed conformal Killing tensor, form or vector, there is then a natural condition on scales. 

Let $(M,\boldsymbol{c})$ be a conformal manifold, let $\sigma\in \Gamma(\mathcal{E}[1])$ be a scale and consider the metric $g = \sigma^{-2}\boldsymbol{g}$, defined on $M\setminus\mathcal{Z}(\sigma)$. 
Let $\ell\ge 1$ and let $\omega_{a_{1}\cdots a_\ell}\in \Gamma(\mathcal{E}{}_{[a_{1}\cdots a_\ell}[\ell+1] )$ be a conformal Killing--Yano form. If, after using $\sigma$ to trivialise weight bundles, $\omega|_{M\setminus\mathcal{Z}(\sigma)}$ is a Killing--Yano form for $(M\setminus\mathcal{Z}(\sigma), g)$, then we call $\sigma$ (and $g$) a \textit{Killing scale for} $\omega$. This includes the case of Killing vectors. 

Let $\ell\ge 2$ and let $k_{a_{1}\cdots a_\ell}\in \Gamma(\mathcal{E}{}_{(a_{1}\cdots a_\ell)_0}[2\ell] )$ be a conformal Killing tensor. If, after using $\sigma$ to trivialise weight bundles, $k{}_{a_{1}\cdots a_\ell}$ is a Killing tensor for $(M\setminus\mathcal{Z}(\sigma),g)$, then we call $\sigma$ (and $g$) a \textit{strong Killing scale for} $k$. The strong condition is very restrictive, since most Killing tensors on a Riemannian manifold are not trace-free. If there is a Killing tensor $r{}_{a_{1}\cdots a_\ell}$ for $(M\setminus\mathcal{Z}(\sigma), g)$ with trace free part $k$, that is, $k_{a_{1}\cdots a_\ell} = r_{(a_{1}\cdots a_\ell)_0}$, then we call $\sigma$ (and $g$) a \textit{Killing scale for} $k$. 

The notions of Killing scales we have defined here extend the notions defined in \cite{gover_conformal_nodate}. 

\section{Special Einstein Products}\label{Sec:SEP}
 
 Let $(M_{1},g_{1})$ and $(M_{2},g_{2})$ be pseudo-Riemannian manifolds of dimensions $n_{1}\ge 3$ and $n_{2}\ge 3$ and signatures $(p_{1},q_{1})$ and $(p_{2},q_{2})$. Then the product manifold $M_{1}\times M_{2}$ has dimension $n = n_{1}+n_{2}$ and its tangent bundle decomposes as $TM = \pi_{1}^*TM_{1}\oplus \pi_{2}^*TM_{2}$, where $\pi_i\colon M\to M{}_{i} $ is the projection. We usually exclude the pullbacks from the notation. Then $M$ is equipped with the product metric $g = g_{1}\times g_{2}$, considering $g_i$ as a section of $TM_i\otimes TM_i)\subset TM\otimes TM$.
 
 Most Riemannian invariants pass nicely to the product. The curvatures are $R{}_{abcd} = R{}^{1}{}_{abcd} +R{}^{2}{}_{abcd} $, $\operatorname{Ric}{}_{ab}  = \operatorname{Ric}{}^{1}{}_{ab} +\operatorname{Ric}{}^{2}{}_{ab} $ and $\operatorname{Sc} = \operatorname{Sc}^1+\operatorname{Sc}^2$. The Levi-Civita connection of $g$ is also the direct sum $\nabla = \nabla^1+\nabla^2$.
 
 Let $(M_{1},\boldsymbol{c}_1)$ and $(M_{2},\boldsymbol{c}_2 )$ be conformal manifolds. Then the product manifold $M = M_{1}\times M_{2}$ is not canonically equipped with a conformal structure derived from $\boldsymbol{c}_1$ and $\boldsymbol{c}_2$. So we need to fix metrics $g_{1}\in \boldsymbol{c}_1$ and $g_{2}\in \boldsymbol{c}_2$ and form the product metric. 
 
 The conformal Schouten does not decompose nicely in general. Except for a special case,

 \begin{proposition}
 	Let $(M_{1},[g_{1}])$, $(M_{2},[g_{2}])$ and $(M=M_{1}\times M_{2},[g] = [g_{1}\times g_{2}])$ be conformal manifolds. Then 
 	$\mathsf{P}{}_{ab}  = \mathsf{P}{}^{1}{}_{ab} +\mathsf{P}{}^{2}{}_{ab} $ if and only if $g_{1}$ and $g_{2}$ are Einstein and $\frac{J^{1}}{n_{1}}+\frac{J^{2}}{n_{2}}=0$.
 \end{proposition}
 \begin{proof}
 We first look at the conformal scalar curvature,
 \[
 	2(n-1)J = \operatorname{Sc}=\operatorname{Sc}_1+\operatorname{Sc}_2 = 2(n_{1}-1)J_{1}+2(n_{2}-1)J_{2}.
 \]
 Then the conformal Schouten tensor of $(M,g)$ is 
 \begin{equation}\label{Eqn:SEP}
 \begin{aligned}
 	(n-2)\mathsf{P}{}_{ab} 
 	&= \operatorname{Ric}{}_{ab} -J\boldsymbol{g}{}_{ab} \\
 	&= \operatorname{Ric}{}^{1}{}_{ab}+\operatorname{Ric}{}^{2}{}_{ab} - \tfrac{n_{1}-1}{n-1}J_1 {g}{}^{1}{}_{ab} -\tfrac{n_{2}-1}{n-1}J_2 {g}{}^{2}{}_{ab} \\
 	&\mathrel{\phantom{=}}-\tfrac{n_{1}-1}{n-1}J_1 {g}{}^{2}{}_{ab} -\tfrac{n_{2}-1}{n-1} J_2 {g}{}^{1}{}_{ab}\\
 	&= (n_{1}-2)\mathsf{P}{}^{1}{}_{ab}+\left(\tfrac{n_{2}}{n-1}J_{1} -\tfrac{n_{2}-1}{n-1}J_{2}\right)g{}^{1}{}_{ab}\\ 
 	&\mathrel{\phantom{=}}+(n_{2}-2)\mathsf{P}{}^{2}{}_{ab} +\left(\tfrac{n_{1}}{n-1}J_{2} -\tfrac{n_{1}-1}{n-1}J_{1}\right)g{}^{2}{}_{ab}\\
 	&= (n-2)\left(\mathsf{P}{}^{1}{}_{ab} + \mathsf{P}{}^{2}{}_{ab} \right)-n_{2}\mathsf{P}{}^{1}{}_{ab}+\left(\tfrac{n_{2}}{n-1}J_{1} -\tfrac{n_{2}-1}{n-1}J_{2}\right)g{}^{1}{}_{ab}\\
 	&\mathrel{\phantom{=}}-n_{1}\mathsf{P}{}^{2}{}_{ab} +\left(\tfrac{n_{1}}{n-1}J_{2} -\tfrac{n_{1}-1}{n-1}J_{1}\right)g{}^{2}{}_{ab}
 \end{aligned}
 \end{equation}
 So $\mathsf{P}{}_{ab}  = \mathsf{P}{}^{1}{}_{ab} +\mathsf{P}{}^{2}{}_{ab} $ if and only if 
 	\[
 		-n_{2}\mathsf{P}{}^{1}{}_{ab} +\left(\tfrac{n_{2}}{n-1}J_{1} -\tfrac{n_{2}-1}{n-1}J_{2}\right) \boldsymbol{g}{}^{1}{}_{ab} =0\quad \text{and} \quad -n_{1}\mathsf{P}{}^{2}{}_{ab} +\left(\tfrac{n_{1}}{n-1}J_{2} -\tfrac{n_{1}-1}{n-1}J_{1}\right)\boldsymbol{g}{}^{2}{}_{ab} =0
 	\] 
 	Which is equivalent to both metrics being Einstein and satisfying
 	\[
 		-\tfrac{n_{2}}{n_{1}}J{}_{1} +\left(\tfrac{n_{2}}{n-1}J_{1} -\tfrac{n_{2}-1}{n-1}J_{2}\right)=0 \quad \text{and} \quad  -\tfrac{n_{1}}{n_{2}}J_2+\left(\tfrac{n_{1}}{n-1}J_{2} -\tfrac{n_{1}-1}{n-1}J_{1}\right)  =0,
 	\] 
 	which both reduce to $\frac{J_{1}}{n_{1}}+\frac{J_{2}}{n_{2}}=0$.
 \end{proof}
 If $(M_{1},g_{1})$ and $(M_{2},g_{2})$ are Einstein manifolds satisfying $\frac{J_1}{n_{1}}+\frac{J_{2}}{n_{2}}=0$, we call the product $M = (M_{1}\times M_2,[g_{1}\times g_{2}])$ a \textit{Special Einstein Product}. 
 
 We will use the scale densities of $g = g_{1}\times g_{2}$, $g_{1}$ and $g_{2}$ to trivialise all density bundles. Consider the first and second jet bundles of $\mathcal{E}[1]$. We have
 \[
 	J^1(\mathcal{E}[1]) = \mathcal{E}[1]\loplus T^*M = \mathcal{E}[1]\loplus T^*M_{1}\oplus T^*M_{2}
 \] 
 and 
 \begin{align*}
 	J^2(\mathcal{E}[1]) 
 	&= \mathcal{E}[1]\loplus T^*M\loplus S^2(T^*M) \\
 	&= \mathcal{E}[1]\loplus T^*M_{1}\oplus T^*M_2\loplus S^2(T^*M_{1})\oplus S^2(T^*M_2)\oplus T^*M_1\odot T^*M_{2}
 \end{align*} 
 Consider the almost Einstein operator of $M$. Let $p = (p_{1},p_{2})\in M$. Let $j^2_p(\tau) = (\tau,\nabla_a\tau, \nabla{}_{(a} \nabla{}_{b)} \tau)\in J^2_p(\mathcal{E}[1])$, then
 \[
 A{}_{ab} (j^2(\tau)) = \nabla{}_{(a} \nabla{}_{b)_0} \tau+\mathsf{P}{}_{(ab)_0} \tau.
 \] 
 Since $\mathsf{P}\in \Gamma(S^2(T^*M_{1})\oplus S^2(T^*M_{2}))$, if $A{}_{ab} (j^2_p(\tau))=0$, then $\nabla{}_{(a} \nabla{}_{b)} \tau\in S^2(T^*_pM_{1})\oplus S^2(T^*_pM_2)$, which is precisely the condition that there exist two jets $j^2_{p_{1}}(\tau_{1})\in J^2_{p_{1}}(\mathcal{E}_{M_{1}}) $ and $j^2_{p_{2}}(\tau_{2})\in J^2_{p_{2}}(\mathcal{E}_{M_{2}}) $ such that $j^2_p(\tau) = j^2_p(\tau_{1})+j^2_p(\tau_{2})$. Then for special Einstein products, we can further characterise the kernel of $A$.
 
 \begin{theorem}\label{Thm:SEP}
 	Let $(M_{1},[g_{1}])$, $(M_{2},[g_{2}])$ and $(M=M_{1}\times M_{2},[g] = [g_{1}\times g_{2}])$ be conformal manifolds and suppose that $(M,g)$ is a special Einstein product. Let $\hat{\mathcal{T}} = \pi_1^*\mathcal{T}_1\oplus \pi_2^*\mathcal{T}_2$ equipped with the metric $\hat{h} = h^1+h^2$. Let $I^1\in \Gamma(\mathcal{T}_1)$ and $I^2\in \Gamma(\mathcal{T}_2)$ be the parallel tractors corresponding to the metrics $g_{1}$ and $g_{2}$. Then, pulling back to $M$, the section $I^1-I^2\in \Gamma(\hat{\mathcal{T}})$ is null for $\hat{h}$ and the tractor bundle of $M$ is $\mathcal{T} = (I^1-I^2)^{\perp}/\langle I^1-I^2 \rangle $. The tractor connection on $\mathcal{T}$ is $\nabla = \nabla^1+\nabla^2$, the direct sum of the connections on $\mathcal{T}{}^{1} $ and $\mathcal{T}^2$.
 \end{theorem}
 \begin{proof}
 	First we calculate
	\[
	|I^1-I^2|^2 = |I^1|^2+|I^2|^2 = -2\left( \frac{J_{1}}{n_{1}}+\frac{J_{2}}{n_{2}} \right) = 0,
	\] 
	since $M$ is a special Einstein product. 

	 Let $j^2_p(\tau)\in \ker(A)_p$. By the above $j^2_p(\tau)$ decomposes as $j^2_p(\tau) = j^2_p(\tau_{1})+j^2_p(\tau_{2})$, where $\tau_{1}$ and $\tau_{2}$ are functions on $M_{1}$ and $M_{2}$, respectively. 
 	Since $M$ is a special Einstein product, we have $\mathsf{P} = \mathsf{P}{}^{1}{}_{ab} +\mathsf{P}{}^{2}{}_{ab} $. So we can write
 	\[
 		\nabla{}_{(a} \nabla{}_{b)} \tau+\mathsf{P}{}_{ab} \tau = \begin{pmatrix} \nabla{}_{(a} \nabla{}_{b)} \tau_{1}+\mathsf{P}{}^{1}{}_{ab} \tau_{1}+\mathsf{P}{}^{1}{}_{ab} \tau_{2} & 0 \\ 0&  \nabla{}_{(a} \nabla{}_{b)} \tau_{2}+\mathsf{P}{}^{2}{}_{ab} \tau_{2}+\mathsf{P}{}^{2}{}_{ab} \tau_{1}\end{pmatrix}.
 	\] 
 	Since $n_{1}P{}^{1}{}_{ab}  = J^1 g{}^{1}{}_{ab} $ and $n_{2}P{}^{2}{}_{ab}  = J^2 g{}^{2}{}_{ab} $, it follows that $j^2_{p_{1}}(\tau_{1})\in \ker(A^1)_p = \mathcal{T}^1_p$ and $j^2_{p_{2}}(\tau_{2})\in \ker(A^2)_p = \mathcal{T}^2_p$, so $j^2_p(\tau)\in \hat{\mathcal{T}}$. Whence we write $j^2_p(\tau) = U = U^1+U^2$, where
 	\[
 		j^2_{p_{1}}(\tau_{1}) = U^1 = \begin{pmatrix} \tau_{1}\\ \mu{}^{1}{}_{a}  \\ \rho_{1} \end{pmatrix}\quad \text{and} \quad
 		j^2_{p_{2}}(\tau_{2}) = U^2 = \begin{pmatrix} \tau_{2}\\ \mu{}^{2}{}_{a}  \\ \rho_{2} \end{pmatrix},
 	\] 
 	with $\rho^1 = -\frac{1}{n_{1}}(\Delta\tau_{1}+J^1\tau_{1})$ and $\rho^2 = -\frac{1}{n_{2}}(\Delta\tau_{2}+J^2\tau_{2})$. Therefore,
 	\begin{equation}\label{Eqn:Atau}
 		\nabla{}_{(a} \nabla{}_{b)} \tau+\mathsf{P}{}_{ab} \tau = \begin{pmatrix} \left(-\rho_{1}+\frac{J^1}{n_{1}}\tau_{2}\right)g{}^{1}{}_{ab}  & 0 \\ 0&  \left(-\rho_{2}+\frac{J^2}{n_{2}}\tau_{1}\right)g{}^{2}{}_{ab} \end{pmatrix}.
 	\end{equation} 
 	So
 	\begin{align*}
 		-\rho_{1}+\frac{J^1}{n_{1}}\tau_{2}&=-\rho_{2}+\frac{J^2}{n_{2}}\tau_{1} \\
 		\implies -\rho_{1}-\frac{J^2}{n_{2}}\tau_{2}&=-\rho_{2}-\frac{J^1}{n_{1}}\tau_{1}\\
 		\implies \rho_{1}-\frac{J^1}{n_{1}}\tau_{1} &= \rho_{2}-\frac{J^2}{n_{2}}\tau_{2}\\
 		\implies \langle I{}^{1} , U{}^{1}  \rangle &= \langle I{}^{2} , U{}^{2}  \rangle \implies \langle U^1+U^2, I^1-I^2 \rangle =0
 	\end{align*}
 	So $j^2_p(\tau)\in (I^1-I^2)^{\perp}$. The decomposition of $j^2_p(\tau)$ is unique up to adding a constant to $\tau_{1}$ and subtracting the same constant from $\tau_{2}$. The constant functions correspond to multiples of $j^2_p(\sigma_{1}) = I^1$ and $j^2_p(\sigma_{2}) = I^2$, so we have defined a canonical injective map $\mathcal{T} = \ker(A)\to (I^1-I^2)^{\perp} / \langle I^1-I^2  \rangle $. Since these two bundles both have rank $n+2$, this is an isomorphism.
 
 	The projectors $X$ and $Z$ in the scale $g$ are given by 
 	\[
 		X(j^2_p(\tau)) = \tau(p) = \tau_{1}(p)+\tau_{2}(p),\quad Z(j^2_p(\tau)) = \nabla_a \tau|_p = \nabla{}_{a} \tau_{1}|_p + \nabla{}_{a} \tau_{2}|_p,
 	\]
 	so $X = \hat{X}+\langle I^1-I^2 \rangle $ and $Z = \hat{Z}+\langle I^1-I^2  \rangle $, where $\hat{X} = X^1+X^2\in \Gamma((I^1-I^2)^{\perp}[1])$ and $\hat{Z} = Z^1+Z^2\in \Gamma((I^1-I^2)^{\perp}\otimes T^*M[1])$.
 	Let $\hat{Y} = Y^1+Y^2\in \Gamma(\hat{\mathcal{T}})$. In general, $\hat{Y}$ is not a section of $(I^1-I^2)^{\perp}$. By \eqref{Eqn:Atau}, the projector $Y$ satisfies
 		\begin{align*}
 			-nY(j^2_p(\tau)) 
 			&= \Delta\tau|_p+J\tau(p) \\
 			&= n\left( -\rho_{1}+\frac{J_{1}}{n_{1}}\tau_{2} \right) = n\left( -\rho_{2}+\frac{J_{2}}{n_{2}}\tau_{1} \right)\\
 			&=-n_{1}\rho_{1}+J_{1}\tau_{2}-n_{2}\rho_{2}+J_{2}\tau_{1}\\
 			&=-n_{1}Y^1(U)-n_{2}Y^2(U)+J_{1}X^2(U)+J_{2}X^1(U)\\
 			&=-n(Y^1(U)+Y^2(U))+n_{2}Y^1(U)+J_{2}X^1(U)+n_{1}Y^2(U)+J_{1}X^2(U)\\
 			&=-n \hat{Y}(U)+n_{2}\left( Y^1(U)-\frac{J_{1}}{n_{1}}X^1(U)  \right) +n_{1}\left( Y^2(U) -\frac{J_{2}}{n_{2}} X^2(U)  \right) \\
 			&=-n \hat{Y}(U)+n_{2}\langle I{}^{1} , U \rangle +n_{1}\langle I{}^{2} , U \rangle \\
 			&=-n \hat{Y}(U)+n\langle I{}^{1} , U \rangle =-n \hat{Y}(U)+n\langle I{}^{2} , U \rangle,
 		\end{align*}	
 		Hence, $Y = \hat{Y}-I^{1}=\hat{Y}-I^{2}\in \Gamma(\mathcal{T})$, since $\mathcal{T} = (I^1-I^2)^{\perp} / \langle I^1-I^2 \rangle $.

		The tractor connections of $\mathcal{T}_1$ and $\mathcal{T}_2$ induce pullback connections on the pullback bundles, and hence on the direct sum $\hat{\mathcal{T}}$. Then $I^1-I^2\in \Gamma(\hat{\mathcal{T}})$ is parallel, so there is a connection induced on $\mathcal{T} = (I^1-I^2)^{\perp} / \langle I^1-I^2 \rangle $. We now verify that this induced connection agrees with the intrinsic connection of $(M_{1}\times M_{2},[g_{1}\times g_{2}])$.

 	The tractor connection of $\mathcal{T}$ satisfies
 	\[
 	\nabla_a X^A = Z{}^{A}{}_{a} = Z{}^{1}{}^{A}{}_{a} +Z{}^{2}{}^{A}{}_{a} = \nabla{}^{1}{}_{a} X^A+\nabla{}^{2}{}_{a} X^A,
 	\]
 
 	\[
 		\begin{aligned}
 			\nabla_a Z{}^{A}{}_{b} 
 			&= -\boldsymbol{g}{}_{ab} Y^A- \mathsf{P}{}_{ab} X{}^{A} \\
 			&= -\boldsymbol{g}{}^{1}{}_{ab} Y{}^{1}{}^{A} -\boldsymbol{g}{}^{2}{}_{ab} Y{}^{2}{}^{A} -\boldsymbol{g}{}^{1}{}_{ab} Y{}^{2}{}^{A} -\boldsymbol{g}{}^{2}{}_{ab} Y{}^{1}{}^{A} +\boldsymbol{g}^2{}_{ab} I{}^{1}{}_{A} +\boldsymbol{g}^1{}_{ab} I{}^{2}{}_{A} \\
 			&\mathrel{\phantom{=}}-\mathsf{P}{}^{1}{}_{ab} X{}^{1}{}^{A} -\mathsf{P}{}^{2}{}_{ab} X{}^{2}{}^{A} -\mathsf{P}{}^{1}{}_{ab} X{}^{2}{}^{A} -\mathsf{P}{}^{2}{}_{ab} X{}^{1}{}^{A}\\
 			&= \nabla{}^{1}{}_{a} Z{}^{A}{}_{b} +\nabla{}^{2}{}_{a} Z{}^{A}{}_{b} -\left( Y{}^{2}{}^{A} -I{}^{2}{}^{A} +\tfrac{J_{1}}{n_{1}}X{}^{2}{}^{A}  \right)\boldsymbol{g}{}^{1}{}_{ab}-\left( Y{}^{1}{}^{A} -I{}^{1}{}^{A} +\tfrac{J_{2}}{n_{2}}X{}^{1}{}^{A}  \right)\boldsymbol{g}{}^{2}{}_{ab}\\
 			&= \nabla{}^{1}{}_{a} Z{}^{A}{}_{b} +\nabla{}^{2}{}_{a} Z{}^{A}{}_{b} -\left( Y{}^{2}{}^{A} -I{}^{2}{}^{A} -\tfrac{J_{2}}{n_{2}}X{}^{2}{}^{A}  \right)\boldsymbol{g}{}^{1}{}_{ab}-\left( Y{}^{1}{}^{A} -I{}^{1}{}^{A} -\tfrac{J_{1}}{n_{1}}X{}^{1}{}^{A}  \right)\boldsymbol{g}{}^{2}{}_{ab}\\
 			&= \nabla{}^{1}{}_{a} Z{}^{A}{}_{b} +\nabla{}^{2}{}_{a} Z{}^{A}{}_{b}
 		\end{aligned}
 	\]
 	and
 	\[
	\nabla_aY^A = \mathsf{P}{}_{ab}  Z{}^{Ab} =\mathsf{P}{}^{1}{}_{ab} Z{}^{1}{}^{Ab} +\mathsf{P}{}^{2}{}_{ab} Z{}^{2}{}^{Ab} = \nabla^1_a Y^A+\nabla^2_a Y^A,
 	\] 
 	where the final equality follows from $I^1$ and $I^2$ being parallel.
 \end{proof}
 
 Suppose that $I{}^{1} $ and $I^2$ are not null. Then $I_{1},I_{2}\not\in (I^1-I^2)^{\perp}$, so neither induce Einstein structures on $M_{1}\times M_{2}$. We can choose as a complementary null tractor, $I^1+I^2$ and the subspace spanned by $I^1-I^2$ and $I^1+I^2$ is $\langle I^1,I^2\rangle$, which is non-degenerate in $\hat{\mathcal{T}}$. Then $\mathcal{T} = (I^1-I^2)^{\perp}/\langle I^1-I^2 \rangle \cong (I^1)^{\perp}\oplus (I^{2})^{\perp}$.

 When $I{}^{1}$ and $I{}^{2} $ are null, the parallel tractors induce a single parallel tractor $I^1 = I^2 \in \mathcal{T}$, corresponding to the Ricci flat product metric.
 
 Special Einstein products will generate a large class of curved examples with a subspace of parallel tractors. If we take one component $(M_{1},g_{1})$ to be conformally flat and the second component $(M_{2},g_{2})$ conformally curved, the result is a curved manifold that still has parallel tractors coming from $M_{1}$. This allows us to apply results later in this article to these curved products.

\section{Rank Two Conformal Killing Tensors}\label{Sec:KT}
Rank two Killing tensors characterise orthogonal separability of the Hamilton-Jacobi equations of classical mechanics, and conformal Killing tensors play a central role in the theory of second-order superintegrable systems. 

In \cite{gover_conformal_nodate}, the full prolongation of conformal Killing tensors was invariantly produced using a half prolongation and the tractor $D$ on conformally flat manifolds. We first review that construction. 

Our full and half prolongations differ from the ones in \cite{gover_conformal_nodate} by a factor of $\frac{1}{4}$, with ours being four times larger. We extract the top slot of a tractor $K{}_{ABCD} \in \Gamma(\begin{picture}(16, 8.0)(0, 2.0)
	\put(-0.2, 10){\line(1,0){10.4}}
	\put(-0.2, 5){\line(1,0){10.4}}
	\put(-0.2, 0){\line(1,0){10.4}}
	\put(0, 10){\line(0,-1){10}}
	\put(5, 10){\line(0,-1){10}}
	\put(10, 10){\line(0,-1){10}}
	\put(10, 0){${}_{\,0}$}
\end{picture}(\mathcal{E}_A))$ by taking
\[
k{}_{bd}  = X^A Z{}^{B}{}_{b} X{}^{C} Z{}^{D}{}_{d} K{}_{ABCD},
\] 
whereas the full prolongation in \cite{gover_conformal_nodate} is of the form
\[
	k{}_{bd} Y{}_{[A} Z{}_{B]}{}^{b} Y{}_{[C} Z{}_{D]}{}^{d} + \textrm{lower order terms}
\] 

Let $k{}_{ab} \in \Gamma(\mathcal{E}{}_{(ab)_0}[4] )$ and let $K{}_{ABCD}  = L(k)_{ABCD}\in \Gamma(\begin{picture}(16, 8.0)(0, 2.0)
	\put(-0.2, 10){\line(1,0){10.4}}
	\put(-0.2, 5){\line(1,0){10.4}}
	\put(-0.2, 0){\line(1,0){10.4}}
	\put(0, 10){\line(0,-1){10}}
	\put(5, 10){\line(0,-1){10}}
	\put(10, 10){\line(0,-1){10}}
	\put(10, 0){${}_{\,0}$}
\end{picture}(\mathcal{E}_A))$ be its full prolongation. Then define the \textit{half prolongation} of $k{}_{ab} $ by $\mathbb{K}{}_{AC} = X^BX^D K{}_{ABCD} $. We first confirm that the half prolongation produced here via the BGG machinery matches the one produced in \cite{gover_conformal_nodate}.

\begin{proposition}\label{Prp:mbbK}
	Let $k{}_{ab} \in \Gamma(\mathcal{E}{}_{(ab)_0}[4] )$ have full prolongation $K{}_{ABCD} $ and half prolongation $\mathbb{K}{}_{BD}$. Then
	\[
	Z{}^{B}{}_{b} Z{}^{D}{}_{d} \mathbb{K}{}_{BD}  = k{}_{bd},\quad  Y^B Z{}^{D}{}_{d} \mathbb{K}{}_{BD}  = -\frac{1}{n+2}\nabla_a k{}^{a}{}_{d} 	\] 
	\[
	Y^BY^D \mathbb{K}{}_{BD} = \frac{1}{n+1}\left( \frac{1}{n+2}\nabla^d\nabla^b k{}_{bd}+\mathsf{P}{}^{bd} k{}_{bd}  \right). 
	\] 
	This gives a one-to-one correspondence between sections $k{}_{ab} \in \Gamma(\mathcal{E}{}_{(ab)_0} [4])$ and sections $\mathbb{K}{}_{AB} \in \Gamma(\mathcal{E}{}_{(AB)_0}[2])$ satisfying
	\[
		X{}^{A} \mathbb{K}{}_{AB} =0 \quad \text{and} \quad D^A \mathbb{K}{}_{AB} =0.
	\] 
\end{proposition}
\begin{proof}
	The composition $\partial^*\circ d^{\nabla}(K)$ for any $K\in \Gamma(\begin{picture}(16, 8.0)(0, 2.0)
	\put(-0.2, 10){\line(1,0){10.4}}
	\put(-0.2, 5){\line(1,0){10.4}}
	\put(-0.2, 0){\line(1,0){10.4}}
	\put(0, 10){\line(0,-1){10}}
	\put(5, 10){\line(0,-1){10}}
	\put(10, 10){\line(0,-1){10}}
	\put(10, 0){${}_{\,0}$}
\end{picture}(\mathcal{E}_A))$ is given by
	\[
		\begin{aligned}
			\partial^*\circ d^{\nabla}(K) 
		&= 2X{}_{[B} Z{}^{Ea} \nabla{}_{a} K{}_{A]ECD} - 2Z{}_{[B}{}^{a} X{}^{E} \nabla{}_{a} K{}_{A]ECD} \\
		&\mathrel{\phantom{=}} + 2X{}_{[D} Z{}^{Ea} \nabla{}_{a} K{}_{C]EAB} - 2Z{}_{[D}{}^{a} X{}^{E} \nabla{}_{a} K{}_{C]EAB},
		\end{aligned}
	\] 
	which vanishes for a full prolongation $K = L(k)$. Let $k{}_{bd} \in \Gamma(\mathcal{E}{}_{(bd)_0}[4] )$ and let $K{}_{ABCD} = L(k)_{ABCD}$. By Theorem \ref{Thm:BGGS}, $\partial^*\circ d^{\nabla}(K)_{ABCD} = 0$. Applying projectors to this expression gives the following identities for $K = L(k)$.
	\[
	X^AY^BX^C Z{}^{D}{}_{d} (\partial^*\circ d^{\nabla})(K)_{ABCD} = X^A Z{}^{Ba} X{}^{C} Z{}^{D}{}_{d} \nabla{}_{a} K{}_{ABCD} =0
	\] 	
	\[
		Z{}^{A}{}_{a} Z{}^{B}{}_{b} X{}^{C} Z{}^{D}{}_{d} (\partial^*\circ d^{\nabla})(K)_{ABCD} = 2Z{}^{A}{}_{[a} X{}^{B}X{}^{C} Z{}^{D}{}_{|d|}\nabla{}_{b]} K{}_{ABCD} =0 
	\] 
	\[
		X^AY^BX^CY^D (\partial^*\circ d^{\nabla})(K)_{ABCD} = 2X^AY^BX^C Z{}^{Da} \nabla{}_{a} K{}_{ABCD} =0.
	\] 
	Other projectors can be applied, but these are what we will need here. Then we have by definition
	\[
	Z{}^{B}{}_{b} Z{}^{D}{}_{d} \mathbb{K}{}_{BD} = X^A Z{}^{B}{}_{b} X^C Z{}^{D}{}_{d} K{}_{ABCD} =k{}_{bd} .
	\] 
	Then we take the divergence of $k{}_{bd} $,
	\begin{align*}
		\nabla{}_{a} k{}^{a}{}_{d} 
		&= -nX^A Y^BX^C Z{}^{D}{}_{d} K{}_{ABCD} +X^A Z{}^{B}{}_{a} Z{}^{Ca} Z{}^{D}{}_{d} K{}_{ABCD} \\ 
		&\mathrel{\phantom{=}}- X^A Z{}^{B}{}_{d} X^CY^D K{}_{ABCD}+X^A Z{}^{Ba} X{}^{C} Z{}^{D}{}_{d} \nabla{}_{a} K{}_{ABCD} \\
		&= -(n+2)X^AY^BX^C Z{}^{D}{}_{d} K{}_{ABCD} = -(n+2)Y^B Z{}^{D}{}_{d} \mathbb{K}{}_{BD},
	\end{align*}
	where we have used the symmetries of $K{}_{ABCD}$, being skew in the pair $AB$ and $CD$, as well as being trace-free. We also used that $Z{}_{A}{}^{a} Z{}_{Ba}  = h{}_{AB} -2X{}_{(A} Y{}_{B)} $. 
	Taking another divergence gives
	\begin{align*}
		\tfrac{-1}{n+2}\nabla_c\nabla_a k{}^{ac}  
		&= Z{}^{Ac} Y^BX^C Z{}^{D}{}_{c} K{}_{ABCD} +\mathsf{P}{}^{ac} X{}^{A} Z{}^{B}{}_{a} X{}^{C} Z{}^{D}{}_{c} K{}_{ABCD} \\ 
		&\mathrel{\phantom{=}}-n X^AY^BX^CY^D K{}_{ABCD} +X^AY^BX^C Z{}^{Dc} \nabla{}_{c} K{}_{ABCD} \\ 
		&= -(n+1) X^AY^BX^CY^D K{}_{ABCD} +\mathsf{P}{}^{ac} k{}_{ac}.
	\end{align*}
	Rearranging and relabelling indices gives the final equality. The claim that the map gives a one-to-one correspondence between $\Gamma(\mathcal{E}{}_{(ab)_0}[4])$ and the given subspace of $\Gamma(\mathcal{E}{}_{(AB)_0}[2] )$ is proved in \cite{gover_conformal_nodate}.
\end{proof}

Now consider when $k{}_{ab} $ is a conformal Killing tensor. Let $g\in \boldsymbol{c}$ have Levi-Civita connection $\nabla$ and consider the tensor $\nabla_a k{}^{a}{}_{b} $. Since $\nabla{}_{(a} k{}_{bc)_0} =0$, it follows that
\[
\nabla{}_{(a} k{}_{bc)} = \rho{}_{(a} \boldsymbol{g}{}_{bc)} 
\] 
for some $\rho\in \Gamma(\mathcal{E}_a[2])$. Then taking a trace gives $2\nabla_a k{}^{a}{}_{b} = (n+2)\rho_b$. 

Let $\sigma\in \Gamma(\mathcal{E}[1])$ be a scale density and let $g = \sigma^{-2}\boldsymbol{g}$ be the corresponding metric, which is defined on $M\setminus\mathcal{Z}(\sigma)$. We say that $\sigma$, or $g$, is a \textit{Killing scale} for a conformal Killing tensor $k{}_{bd} \in \Gamma(\mathcal{E}{}_{(bd)_0} [4])$ if there is a density $f\in \Gamma(\mathcal{E}[2])$ such that $r{}_{ab}  = k{}_{ab} +f\boldsymbol{g}{}_{ab} $ is a Killing tensor for $(M\setminus\mathcal{Z}(\sigma), g)$, that is,
\[
\nabla^g{}_{(a} r{}_{bd)} =0
\] 
on $M\setminus\mathcal{Z}(\sigma)$. When a scale has zeros, the Killing tensor, its trace adjustment, and its projective prolongation are all understood on \(M\setminus\mathcal Z(\sigma)\). We call $\sigma$, or $g$, a \textit{strong Killing scale} for $k{}_{bd} $ if $k{}_{bd} $ itself is a Killing tensor for $(M\setminus \mathcal{Z}(\sigma), g)$, that is,
\[
\nabla{}^{g}{}_{(a} k{}_{bd)} =0
\] 
on $M\setminus\mathcal{Z}(\sigma)$. 

\begin{proposition}\label{Prp:SKS}
	Let $K{}_{ABCD} $ be the full prolongation of a tensor $k{}_{bd} \in \Gamma(\mathcal{E}{}_{(bd)_0}[4])$ and let $\sigma\in \Gamma(\mathcal{E}_+[1])$ be a scale with scale tractor $I{}_{A} $. Then
	\[
		\frac{1}{n+2}\nabla_a k{}^{a}{}_{d} = \sigma^{-1}I^AX^BX^C Z{}^{D}{}_{d} K{}_{ABCD}.
	\] 
	Suppose $k{}_{bd} $ is a conformal Killing tensor. Then $\sigma$ is a strong Killing scale for $k{}_{bd} $ if and only if  
	\[
		I^AX^BX^C Z{}^{D}{}_{d} K{}_{ABCD}=0.
	\] 
	If $\sigma$ is an Einstein scale and $k{}_{bd} $ is a normal solution, that is when $I$ and $K$ are parallel, then $\sigma $ is a strong Killing scale if and only if
	\[
		I^A K{}_{ABCD} =0.
	\] 
\end{proposition}
\begin{proof}
	The identity follows from Proposition \ref{Prp:mbbK}. Calculating using the Levi-Civita connection of $g = \sigma^{-2}\boldsymbol{g}$ gives $I_A = \sigma Y_A-\frac{J}{n}\sigma X{}_{A} $. So
	\[
	I^AX^BX^C Z{}^{D}{}_{d} K{}_{ABCD} = -\sigma X^AY^BX^C Z{}^{D}{}_{d} K{}_{ABCD}  = -\sigma Y^B Z{}^{D}{}_{d} \mathbb{K}{}_{BD} = \frac{\sigma}{n+2}\nabla_a k{}^{a}{}_{d} .
	\] 

	When $I$ and $K$ are parallel, $I{}^{A} K{}_{ABCD} \in \Gamma(\begin{picture}(16, 8.0)(0, 2.0)
	\put(-0.2, 10){\line(1,0){10.4}}
	\put(-0.2, 5){\line(1,0){10.4}}
	\put(-0.2, 0){\line(1,0){5.4}}
	\put(0, 10){\line(0,-1){10}}
	\put(5, 10){\line(0,-1){10}}
	\put(10, 10){\line(0,-1){5}}
	\put(10, 0){${}_{\,0}$}
\end{picture}(\mathcal{E}_A))$ is a parallel section of an irreducible tractor bundle so if the top slot 
\[
I^AX^BX^C Z{}^{D}{}_{d} K{}_{ABCD} 
\] 
vanishes, then $I^A K{}_{ABCD} =0$. 

\end{proof}

\begin{lemma}\label{Lem:KS}
	Let $k{}_{bd}\in \Gamma(\mathcal{E}{}_{(bd)_0} [4]) $ be a conformal Killing tensor and let $\sigma\in \Gamma(\mathcal{E}[1])$ be a scale. Then $\sigma$ is a Killing scale for $k{}_{bd} $ if and only if 
	\[
	\sigma^{-2}\nabla_{a} k{}^{a}{}_{d}  = -(n+2)\sigma^{-3}I{}^{A} X{}^{B} Z{}^{C}{}_{d} X{}^{D} K{}_{ABCD} 
	\]
	is an exact one form on $M\setminus\mathcal{Z}(\sigma)$.
\end{lemma}
\begin{proof}
Let $f\in \Gamma(\mathcal{E}[2])$ and let $r{}_{bd}   =k{}_{bd} +f\boldsymbol{g}_{bd}$. Then, calculating in the scale $\nabla = \nabla{}^{g}$,
\[
	\nabla{}_{(a} r{}_{bd)}  = \nabla{}_{(a} k{}_{bd)} +\nabla{}_{(a} f \boldsymbol{g}{}_{bd)} \implies \nabla{}_{(a} r{}_{bd)_0} =0 \quad \text{and} \quad  3\boldsymbol{g}{}^{ab} \nabla{}_{(a} r_{bd)} = 2\nabla{}_{a} k{}^{a}{}_{d} + (n+2)\nabla_d f.
\]
So $r{}_{bd} $ is Killing if and only if $\nabla_d f = -\frac{2}{n+2}\nabla{}_{a} k{}^{a}{}_d$. Since $\sigma$ is parallel for $\nabla{}^{g} $, we use $\sigma$ to trivialise density bundles, and existence of such a $f\in \Gamma(\mathcal{E}[2])$ is equivalent to exactness of the unweighted one form $\sigma^{-2}\nabla{}_{a} k{}^{a}{}_{d}  $.
\end{proof}

We say that a scale $\sigma\in \Gamma(\mathcal{E}[1])$, or the metric $g = \sigma^{-2}\boldsymbol{g}$ is a \textit{local Killing scale} for a conformal Killing tensor $k{}_{bd} \in \Gamma(\mathcal{E}{}_{(bd)_0} [4])$ if for all $p\in M\setminus\mathcal{Z}(\sigma)$ there exists an open subset $U\subset M\setminus\mathcal{Z}(\sigma)$ containing $p$ on which $\sigma|_U$ is a Killing scale.

Therefore we immediately have that $\sigma$ is a local Killing scale if and only if $\sigma^{-2}\nabla_a k{}^{a}{}_d $ is closed. This has a compact form in terms of tractors.

\begin{theorem}\label{Thm:LKS}
	Let $K{}_{ABCD} $ be the full prolongation of a conformal Killing tensor $k{}_{bd} \in \Gamma(\mathcal{E}{}_{(bd)_0} [4])$and let $\sigma$ be a scale with scale tractor $I{}_{A}$. Then $\sigma$ is a local Killing scale for $k{}_{bd}$ if and only if
	\begin{equation}\label{Eqn:IIDIDK}
		X{}^{E} Z{}^{A}{}_{a} Z{}^{B}{}_{b}  X{}^{D} \left(3I{}^{C}I{}_{[E} K_{AB]CD}{}-\sigma(D{}_{[E} I{}^{C}) K{}_{AB]CD} -\sigma I{}^{C} D{}_{[E} K{}_{AB]CD} \right)=0.
	\end{equation} 
\end{theorem}

\begin{proof}
	The following calculation holds in all scales, despite neither side being invariant, 
	\begin{align*}
		\nabla{}_{[a} (I^AX^B Z{}^{C}{}_{c]} X{}^{D} K{}_{ABCD})  
		&= I^A Z{}^{B}{}_{[a} Z{}^{C}{}_{c]} X{}^{D} K{}_{ABCD} +I^AX^B Z{}^{C}{}_{c} Z{}^{D}{}_{a} K{}_{ABCD}\\ 
		&\mathrel{\phantom{=}}+ \nabla{}_{[a} I{}^{A} X{}^{B} Z{}^{C}{}_{c]} X{}^{D} K{}_{ABCD} + I^A X^B Z{}^{C}{}_{[c}X{}^{D} \nabla{}_{a]} K{}_{ABCD}\\ 
		&=\tfrac{3}{2}I^A X^B Z{}^{C}{}_{c} Z{}^{D}{}_{a} K{}_{ABCD} + (D_E I{}^{A})Z{}^{E}{}_{[a} X{}^{B} Z{}^{C}{}_{c]} X{}^{D} K{}_{ABCD} \\ 
		&\mathrel{\phantom{=}}+ I{}^{A} X{}^{B} Z{}^{C}{}_{[c} X{}^{D} Z{}^{E}{}_{a]} D{}_{E} K{}_{ABCD} 
	\end{align*}
	since $\nabla{}_{[a} Z{}^{B}{}_{b]} =0$ and by the symmetries of $K$. Now calculate in the scale $\sigma$, so $Z{}^{A}{}_{a} I_A=0$ and $\sigma$ is a local Killing scale if and only if
	\begin{align*}
		0&=2\nabla_{[a} (I^AX^B Z{}^{C}{}_{c]} X{}^{D} K{}_{ABCD})\\ 
		&= 9 \sigma^{-1}I^AX^B Z{}^{C}{}_{c} Z{}^{D}{}_{a} X^E K{}_{AB[CD} I{}_{E]}-3 X^B Z{}^{C}{}_{c} Z{}^{D}{}_{a} X{}^{E} K{}_{AB[CD} D{}_{E]} I^A \\ 
		&\mathrel{\phantom{=}}-3I^AX^B Z{}^{C}{}_{c} Z{}^{D}{}_{a} X{}^{E} D{}_{[E} K{}_{|AB|CD]} 
	\end{align*}
	The result now follows by applying $K{}_{ABCD}  = K{}_{CDAB} $ and renaming indices.
\end{proof}
Both of the tractors $K$ and $I$ are unweighted, so $D_A I_B = Z{}_{A}{}^{a} \nabla_a I_B+X_A(*)$ and $D_A K{}_{BCDE} = Z{}_{A}{}^{a} \nabla_a K{}_{BCDE} +X_A(*)$, and the projectors in \eqref{Eqn:IIDIDK} annihilate the $X$ term in both cases. So when $I$, $K$ or both are parallel, we obtain special cases of the above Theorem.

We state explicitly the following, which has already been useful in the theory of superintegrable systems \cite{goodwin_superintegrable_2026}, which often take place on conformally flat manifolds.
\begin{corollary}\label{Cor:Rk2}
	Let $K{}_{ABCD} $ be the full prolongation of a conformal Killing tensor and let $\sigma$ be a scale with scale tractor $I{}_{A}$. Suppose that $(M,\boldsymbol{c})$ is conformally flat, or more generally that $k{}_{bd} $ is a normal solution of the conformal Killing equation. Then $\sigma$ is a local Killing scale for $k{}_{bd}$ if and only if
	\[
		X{}^{E} Z{}^{A}{}_{a} Z{}^{B}{}_{b}  X{}^{D} \left(3I{}^{C}I{}_{[E} K{}_{AB]CD} - \sigma (D{}_{[E} I{}^{C}) K{}_{AB]CD} \right)=0.
	\] 
\end{corollary}
The strongest case is when both $K$ and $I$ are parallel. Then all projectors can be removed.
The second part of the below Theorem was proved for conformally flat non-Ricci-flat Einstein scales in \cite{gover_conformal_nodate}. In this case our proof can be reduced to their proof. But in the Ricci flat case, our methods are required.

For symmetric tractors $Q{}_{AB} , P{}_{AB} \in \Gamma(\mathcal{E}{}_{(AB)} )$, we define
\[
	(Q\owedge P)_{ABCD} = Q{}_{AC} P{}_{BD} +Q{}_{BD} P{}_{AC} -Q{}_{AD} P{}_{BC} -Q{}_{BC} P{}_{AD}, 
\] 
which is the standard Kulkarni-Nomizu product, adapted to tractor indices. The Kulkarni-Nomizu product has curvature symmetries, $(Q\owedge P)_{ABCD}\in \Gamma(\begin{picture}(10, 8.0)(0, 2.0)
	\put(-0.2, 10){\line(1,0){10.4}}
	\put(-0.2, 5){\line(1,0){10.4}}
	\put(-0.2, 0){\line(1,0){10.4}}
	\put(0, 10){\line(0,-1){10}}
	\put(5, 10){\line(0,-1){10}}
	\put(10, 10){\line(0,-1){10}}
\end{picture}(\mathcal{E}_A))$. In fact, the Kulkarni-Nomizu product is a projection of the following isomorphism
\[
S^2(S^2(V)) \cong  \begin{picture}(16, 11.0)(0, 5.0)
	\put(-0.2, 16){\line(1,0){16.4}}
	\put(-0.2, 8){\line(1,0){16.4}}
	\put(-0.2, 0){\line(1,0){16.4}}
	\put(0, 16){\line(0,-1){16}}
	\put(8, 16){\line(0,-1){16}}
	\put(16, 16){\line(0,-1){16}}
	\end{picture}\;(V)\oplus \begin{picture}(32, 7.0)(0, 1.0)
	\put(-0.2, 8){\line(1,0){32.4}}
	\put(-0.2, 0){\line(1,0){32.4}}
	\put(0, 8){\line(0,-1){8}}
	\put(8, 8){\line(0,-1){8}}
	\put(16, 8){\line(0,-1){8}}
	\put(24, 8){\line(0,-1){8}}
	\put(32, 8){\line(0,-1){8}}\end{picture}\;(V)
\] 
onto the first factor, for a vector space $V$. This is particularly useful when $V$ is equipped with a non-degenerate bilinear form and we consider the trace decomposition
\[
\begin{picture}(22, 11.0)(0, 5.0)
	\put(-0.2, 16){\line(1,0){16.4}}
	\put(-0.2, 8){\line(1,0){16.4}}
	\put(-0.2, 0){\line(1,0){16.4}}
	\put(0, 16){\line(0,-1){16}}
	\put(8, 16){\line(0,-1){16}}
	\put(16, 16){\line(0,-1){16}}
\end{picture} (V) = \begin{picture}(22, 11.0)(0, 5.0)
	\put(-0.2, 16){\line(1,0){16.4}}
	\put(-0.2, 8){\line(1,0){16.4}}
	\put(-0.2, 0){\line(1,0){16.4}}
	\put(0, 16){\line(0,-1){16}}
	\put(8, 16){\line(0,-1){16}}
	\put(16, 16){\line(0,-1){16}}
	\put(16, 0){${}_{\,0}$}
\end{picture}(V)\oplus \begin{picture}(22, 7.0)(0, 1.0)
	\put(-0.2, 8){\line(1,0){16.4}}
	\put(-0.2, 0){\line(1,0){16.4}}
	\put(0, 8){\line(0,-1){8}}
	\put(8, 8){\line(0,-1){8}}
	\put(16, 8){\line(0,-1){8}}
	\put(16, 0){${}_{\,0}$}
\end{picture}(V)\oplus \mathbb{R}.
\] 
Now specialising to tractor indices, we can write
\[
R{}_{ABCD} = K{}_{ABCD} +(P\owedge h)_{ABCD},
\] 
where 
\[
	P{}_{BD}  = \frac{1}{n}\left(R{}_{AB}{}^{A}{}_{D} - \frac{1}{2(n+1)}R{}_{AC}{}^{AC} h{}_{BD}\right).
\]
Then the resulting $K{}_{ABCD} $ is trace-free, and is called the \textit{trace-free part} of $R{}_{ABCD}$.

\begin{theorem}\label{Thm:IWI}
	Let $(M,\boldsymbol{c})$ be a conformal manifold, let $k{}_{bd}  $ be a normal solution of the conformal Killing equation with parallel prolongation $K{}_{ABCD} $ and let $I{}_{A} = {D}{}_{A} \sigma$ be a parallel tractor. Then $\sigma$ is a local Killing scale for $k{}_{bd} $ if and only if
	\begin{equation}\label{Eqn:IWI}
		I{}^{A} K{}_{AB[CD} I{}_{E]} =0.
	\end{equation} 
	If there exists a parallel tractor $\hat{I}{}_{A} $ such that $\hat{I}^AI_A=1$, then $\sigma$ is a Killing scale for $k{}_{bd} $ if and only if \eqref{Eqn:IWI} holds.
\end{theorem}
\begin{proof}
	The first part of the Theorem follows from Theorem \ref{Thm:LKS}. The identity reduces to
	\[
		X^E Z{}^{A}{}_{a} Z{}^{B}{}_{b} X{}^{D} I{}^{C} I{}_{[E} K{}_{AB]CD} =0.
	\] 
	Changing the indices, this is equivalent to the $X^B Z{}^{C}{}_{c} Z{}^{D}{}_{d} X^E$ slot of the parallel tractor $I^A K{}_{AB[CD} I{}_{E]} $. This tractor satisfies
	\[
		3I^A K{}_{AB[CD} I{}_{E]} h{}^{BE} = I^AI^B K{}_{ABCD} =0 \quad \text{and} \quad I^A K{}_{A[BCD} I{}_{E]} =0
	\] 
	which are the symmetries of the tractor bundle 
	\[
	\begin{picture}(22, 15.0)(0, 9.0)
	\put(-0.2, 24){\line(1,0){16.4}}
	\put(-0.2, 16){\line(1,0){16.4}}
	\put(-0.2, 8){\line(1,0){8.4}}
	\put(-0.2, 0){\line(1,0){8.4}}
	\put(0, 24){\line(0,-1){24}}
	\put(8, 24){\line(0,-1){24}}
	\put(16, 24){\line(0,-1){8}}
	\put(16, 0){${}_{\,0}$}
\end{picture}(\mathcal{E}_A) = \{S\in \mathcal{E}{}_{BCDE} : S{}_{BCDE}  = S{}_{B[CDE]},\quad S{}_{[BCDE]} =0\quad \text{and} \quad h{}^{BE} S{}_{BCDE} =0\},
	\] 
\\
which is an irreducible tractor bundle. The top slot (homology) of a section of this bundle is given by the projector $X^B Z{}^{C}{}_{c} Z{}^{D}{}_{d} X^E$. So by Lemma \ref{Lem:BGG},
\begin{equation*}
	I^AX^B Z{}^{C}{}_{c} Z{}^{D}{}_{d} X{}^{E} K{}_{AB[CD} I{}_{E]} =0\iff I^A K{}_{AB[CD} I{}_{E]} =0
\end{equation*}

The remaining statement to prove is that $\sigma$ is a Killing scale when \eqref{Eqn:IWI} holds. For this, we construct a parallel projective tractor with $\begin{picture}(10, 8.0)(0, 2.0)
	\put(-0.2, 10){\line(1,0){10.4}}
	\put(-0.2, 5){\line(1,0){10.4}}
	\put(-0.2, 0){\line(1,0){10.4}}
	\put(0, 10){\line(0,-1){10}}
	\put(5, 10){\line(0,-1){10}}
	\put(10, 10){\line(0,-1){10}}
\end{picture}$ curvature symmetries and top slot $r{}_{bd}$ such that $r{}_{(bd)_0}  = k{}_{bd}$. Using Theorem \ref{Thm:Iperp}, it suffices to construct a parallel conformal tractor $R{}_{ABCD} \in \Gamma(\begin{picture}(10, 8.0)(0, 2.0)
	\put(-0.2, 10){\line(1,0){10.4}}
	\put(-0.2, 5){\line(1,0){10.4}}
	\put(-0.2, 0){\line(1,0){10.4}}
	\put(0, 10){\line(0,-1){10}}
	\put(5, 10){\line(0,-1){10}}
	\put(10, 10){\line(0,-1){10}}
\end{picture}(\mathcal{E}_A))$ with top slot $r{}_{bd} $ such that $I^A R{}_{ABCD} =0$. 
	Suppose that $\hat{I}$ is parallel and that $I_A\hat{I}^A=1$. Then let
	\[
		U{}_{BC} = I^A\hat{I}^E K{}_{ABCE} - I^A\hat{I}^E\hat{I}^F K{}_{AEFB} I{}_{C} 
	\] 
	First we show that $U$ is symmetric. We first apply \eqref{Eqn:IWI},
	\begin{equation}\label{Eqn:RhoS}
		0 = 3I^A\hat{I}^F K{}_{AE[CB} I{}_{F]}  = I^A K{}_{AECB} +2I^A\hat{I}^F K{}_{AEF[C} I{}_{B]}.
	\end{equation} 
	Since $K{}_{A[BCE]} =0$, 
	\[
		2I^A \hat{I}^E K{}_{A[BC]E} = -I^A \hat{I}^E K{}_{AEBC} = I^A \hat{I}{}^{E} K_{AECB}.
	\] 
	So contracting $\hat{I}{}^{E}$ into \eqref{Eqn:RhoS} gives
	\[
		\begin{aligned}
			0 
			&= I^A \hat{I}^E K_{AECB} + 2I^A \hat{I}^E \hat{I}^F K_{AEF[C}I{}_{B]} \\
			&= 2 I^A \hat{I}^E K{}_{A[BC]E} -2 I^A \hat{I}^E \hat{I}^F K{}_{AEF[B}I{}_{C]}= 2U{}_{[BC]} 
		\end{aligned}
	\] 
Then, by \eqref{Eqn:RhoS}, we have
\[
	2U{}_{B[C} I{}_{D]} =-2 I{}^{A} \hat{I}^E K{}_{ABE[C} I{}_{D]} = I^A K{}_{ABCD} 
\] 
with $U$ parallel. Then $U{}_{B}{}^{B} = -I^A \hat{I}^E \hat{I}^F I^B K{}_{AEFB}$, so
\begin{equation}\label{Eqn:UWU}
	I^B U{}_{BC} = I^AI^B\hat{I}^E K{}_{ABCE} -I^A\hat{I}^E\hat{I}^F I^B K{}_{AEFB} I_C = U{}_{B}{}^{B} I_C
\end{equation} 
Then let ${P}{}_{BC} = U{}_{BC} -\frac{1}{2}U{}_{A}{}^{A} h{}_{BC} $ and $R{}_{ABCD}  = K{}_{ABCD} +(h\owedge {P})_{ABCD}$. Then
\begin{align*}
	I{}^{A} R{}_{ABCD} 
	&= I{}^{A} K{}_{ABCD} +2{P}{}_{B[D} I{}_{C]} +2I{}^{A} {P}{}_{A[C} h{}_{D]B} \\ 
	&= 2U{}_{B[C} I{}_{D]} + 2U{}_{B[D} I{}_{C]} -U{}_{A}{}^{A} h{}_{B[D} I{}_{C]} +U{}_{A}{}^{A} I{}_{[C} h{}_{D]B} = 0
\end{align*} 
So $R\in I^{\perp} $, $R$ has curvature symmetries and $R$ is parallel. The top slot of $R$ is
\[
X{}^{A} Z{}^{B}{}_{b} X{}^{C} Z{}^{D}{}_{d} R{}_{ABCD}  = X{}^{A} Z{}^{B}{}_{b} X{}^{C} Z{}^{D}{}_{d} K{}_{ABCD}+X{}^{A} X{}^{C} {P}{}_{AC} \boldsymbol{g}{}_{bd} = k{}_{bd} + f \boldsymbol{g}{}_{bd} ,
\] 
where $f = X^AX^C P{}_{AC} \in \Gamma(\mathcal{E}[2])$.
\end{proof}
The condition that such a $\hat{I}$ exists is automatic when $I$ is not null. In this case, let $\hat{I}=I / |I|^2$. Then the tractors in the above proof are
\[
U{}_{BC}  = \frac{1}{|I|^2}I^AI^D K{}_{ABCD} ={P}{}_{BC},
\] 
which matches the proof in \cite{gover_conformal_nodate}. The remaining freedom for ${P}{}_{BC} $ in this case is multiples of $h{}_{BD} -\frac{2}{|I|^2} I_BI_D$, which correspond to adding constant multiples of $g{}_{bd} $ to $r{}_{bd} $. The choice of ${P}{}_{BC} $ in our proof is the unique choice which is orthogonal to $I$. There are many cases where this is not the most natural choice. Indeed, when we start with $R{}_{ABCD}$ and form $K{}_{ABCD} = R{}_{ABCD} -(P\owedge {h})_{ABCD}$, we take
\[
{P}{}_{BD}  = \frac{1}{n}\left( R{}_{AB}{}^{A}{}_{D} -\frac{1}{2(n+1)} R{}_{AC}{}^{AC}  h{}_{BD}  \right),
\] 
which is orthogonal to $I$ precisely when $R{}_{AC}{}^{AC} =0$. When this is not the case, the process of going from $R$ to $K$ and back will not return the same $R$.

For null $I$, that is, scalar flat metrics, the remaining gauge freedom is given by adding multiples of $I_BI_D$. In the presence of $\hat{I}$, we can still make a distinguished choice of ${P}{}_{BD} $. However, without $\hat{I}$, it is much less clear how to remove this ambiguity. The quantities
\[
	I^A {P}{}_{AB}, \quad {P}{}_{A[B} I{}_{C]} \quad \text{and} \quad {P}{}_{A}{}^{A} 
\] 
are all unchanged by adding a multiple of $I_BI_D$. So we cannot form a uniform normalisation condition. We conclude this section by showing that on every conformally flat manifold there locally exists a rank-two conformal Killing tensor, which has no true strong Killing scale, that is, a Killing scale $\sigma$ which vanishes nowhere.
\begin{proposition}
	Let $(M,\boldsymbol{c})$ be conformally flat and let $p\in M$. Then locally there is a rank two conformal Killing tensor such that every strong Killing scale for $k{}_{ab} $ vanishes at $p$.
\end{proposition}
\begin{proof}
	Let $p\in M$ and let $K{}_{ABCD}|_p\in \begin{picture}(16, 8.0)(0, 2.0)
	\put(-0.2, 10){\line(1,0){10.4}}
	\put(-0.2, 5){\line(1,0){10.4}}
	\put(-0.2, 0){\line(1,0){10.4}}
	\put(0, 10){\line(0,-1){10}}
	\put(5, 10){\line(0,-1){10}}
	\put(10, 10){\line(0,-1){10}}
	\put(10, 0){${}_{\,0}$}
\end{picture}(\mathcal{E}_A|_p)$. Let $\mathbb{K}{}_{AC} = X^BX^D K{}_{ABCD} \in \mathcal{E}{}_{(AC)_0}|_p $. Since $X^A \mathbb{K}{}_{AC} =0$, we can define a linear map $\lambda\colon \mathcal{E}_A|_p\to X^{\perp}|_p / \langle X  \rangle$ by
\[
\lambda(V_A) = \mathbb{K}{}_{AB} V^A \mod X
\] 
The kernel of $\lambda$ contains $X$. Looking at the slots of $\mathbb{K}$ from Proposition \ref{Prp:mbbK}, we see that if $V_A = \tau Y_A+ \mu_a Z{}_{A}{}^{a} +\omega X_A $ , then
\[
\lambda(V_A) = \left(k{}^{a}{}_{b} \mu_a  -\tfrac{1}{n+2} \tau\nabla_a k{}^{a}{}_{b}\right) Z{}_{B}{}^{b} \mod X.
\] 
So, if $\nabla_a k{}^{a}{}_{b}|_p \not\in \operatorname{im}(k(p))$ (the image of the endomorphism induced by $k_{ab}|_p$), then $\tau$ is forced to vanish at $p$, so the kernel of $\lambda$ is contained in $X^{\perp}$. Since we can choose a $K{}_{ABCD} $ arbitrarily at $p$, we can take it such that $k(p)$ is a trace-free degenerate endomorphism and that $\nabla_a k{}^{a}{}_{b}|_p \not\in \operatorname{im}(k(p))$. Now if $\sigma$ is a strong Killing scale, then $I_A = D_A\sigma$ is in the kernel of $\lambda$. Hence, $I_A|_p\in X^{\perp}$ and $\sigma = X^AI_A=0$.
\end{proof}

\section{Conformal Killing--Yano Forms}\label{Sec:CKY}

In this section, we consider the form analogue of Killing tensors. For a pseudo-Riemannian manifold, an $\ell$-form $\omega_{a_{1}\cdots a_\ell}$ $(1\le \ell\le n-1)$ is called a \textit{Killing--Yano form} if it satisfies
\[
\nabla{}_{(a_{0}}\omega{}_{a_{1})a_{2}\cdots a_\ell}=0, \text{ equivalently, } \nabla{}_{a_{0}}\omega{}_{a_{1}\cdots a_\ell} = \frac{1}{\ell+1}d\omega{}_{a_{0}\cdots a_\ell}.
\] 
Killing--Yano forms are interesting for many reasons. Firstly, on pseudo-Riemannian manifolds, they generate Killing tensors. If $\omega\in \Lambda^\ell(T^*M)$ is a Killing--Yano form, then $r{}_{ab}  = \omega_{a}{}^{a_{2}\cdots a_\ell}\omega{}_{ba_{2}\cdots a_\ell} $ is a symmetric tensor satisfying
\[
\nabla{}_{(c} r{}_{ab)} = \nabla{}_{(c}\omega{}_{a}{}^{a_{2}\cdots a_\ell}\omega{}_{b)a_{2}\cdots a_\ell} + \omega{}_{(a}{}^{a_{2}\cdots a_\ell}\nabla{}_{c}\omega{}_{b)a_{2}\cdots a_\ell}=0.
\] 
So $r{}_{ab}$ is a Killing tensor with ``square root'' $\omega$. The Killing--Yano operator is projectively invariant and is the BGG operator for the completely skew projective tractor bundle $\mathcal{E}_{[\alpha_{0}\alpha_{1}\cdots \alpha_{\ell}]}$, with top slot $\mathcal{E}_{[a_{1}\cdots a_\ell]}(\ell+1)$. There is also a natural conformal analogue, which is the BGG operator for the completely skew conformal tractor bundle $\mathcal{E}_{[A_{0}A_{1}\cdots A_{\ell}]}$
\[
	\Gamma(\mathcal{E}_{[a_{1}\cdots a_\ell]}[\ell+1])\ni\omega_{a_{1}\cdots a_\ell}\mapsto \operatorname{trace-free}(\nabla{}_{(a_{0}}\omega{}_{a_{1})a_{2}\cdots a_\ell})\in \Gamma\left(\begin{picture}(16, 13.0)(0, 7.0)
	\put(-0.2, 20){\line(1,0){10.4}}
	\put(-0.2, 15){\line(1,0){10.4}}
	\put(1.5, 12.5){${}_{\cdot}$}
	\put(1.5, 10){${}_{\cdot}$}
	\put(1.5, 7.5){${}_{\cdot}$}
	\put(-0.2, 5){\line(1,0){5.4}}
	\put(-0.2, 0){\line(1,0){5.4}}
	\put(0, 20){\line(0,-1){20}}
	\put(5, 20){\line(0,-1){20}}
	\put(10, 20){\line(0,-1){5}}
	\put(10, 0){${}_{\,0}$}
\end{picture}(T^*M)[\ell+1]\right),
\] 
called the conformal Killing--Yano operator. A tensor in the kernel of this operator is a \textit{conformal Killing--Yano form}.

If $\nabla$ is a Levi-Civita connection of a metric in $\boldsymbol{c}$, then a conformal Killing--Yano form satisfies
\[
	\nabla{}_{a_{0}}\omega{}_{a_{1}a_{2}\cdots a_\ell}-\frac{1}{\ell+1}(d\omega)_{a_{0}\cdots a_\ell} =\ell \boldsymbol{g}_{a_{0}[a_{1}}\rho_{a_{2}\cdots a_\ell]},
\] 
where $\rho_{a_{2}\cdots a_\ell} = \frac{1}{n-\ell+1}\nabla{}_{a}{\omega}{}^{a}{}_{a_{2}\cdots a_\ell} $.

In terms of Killing scales, there is no strong case, since no trace adjustment of a form is possible. Hence, a scale $\sigma$ is a \textit{Killing scale} for $\omega$ if $\rho=0$. Additionally, there is also no need to consider local Killing scales separately.

\begin{theorem}\label{Thm:KYKS}
	Let $\Omega_{A_{0}A_{1}\cdots A{}_{\ell} }$  be the full prolongation of an $\ell$-form $\omega_{a_{1}\cdots a_\ell}\in \Gamma(\mathcal{E}_{[a_{1}\cdots a_\ell]}[\ell+1])$. Let $\sigma\in \Gamma(\mathcal{E}_+[1])$ be a scale with scale tractor $I{}_{A}$. Then
	\[
		\frac{\sigma}{n-\ell+1}\nabla^b \omega{}_{ba_{2}\cdots a_\ell} = I{}^{A_{0}} X{}^{A_{1}}Z{}^{A_{2}}{}_{a_{2}}\cdots Z{}^{A_\ell}{}_{a_\ell} \Omega_{A_{0}A_{1}\cdots A_\ell}.
	\] 
	Suppose that $\omega{}_{a_{1}\cdots a_\ell} $ is a conformal Killing--Yano form. Then $\sigma$ is a Killing scale for $\omega$ if and only if
	\[
	I{}^{A_{0}} X{}^{A_{1}}Z{}^{A_{2}}{}_{a_{2}}\cdots Z{}^{A_\ell}{}_{a_\ell} \Omega_{A_{0}A_{1}\cdots A_\ell}=0.
	\] 
	If $\sigma$ is an Einstein scale and $\Omega$ is parallel, then $\sigma$ is a Killing scale if and only if
	\[
		I{}^{A_{0}} \Omega_{A_{0}A_{1}\cdots A_\ell}=0.
	\] 
\end{theorem}
\begin{proof}
	The full prolongation $\Omega_{A_{0}A_{1}\cdots A_\ell}$ satisfies $\partial^*\circ d^{\nabla}(\Omega)=0$, or equivalently
	\[
		\mathbb{X}{}^{B}{}_{[A_{0}|}{}^a\nabla_{a}\Omega_{B|A_{1}\cdots A_\ell]}=0.
	\] 
	This immediately gives that
	\[
		X{}^{A_{0}}Z{}^{A_{1}}{}_{a_{1}}\cdots Z^{A_\ell}{}_{a_\ell}\nabla_{a_{0}}\Omega_{A_{0}\cdots A_\ell}
	\] 
	has vanishing completely skew part and is also trace-free. Therefore,
	\begin{align*}
		\nabla^b\omega_{ba_{2}\cdots a_\ell} 
		&= \nabla^b\left( X{}^{A_{0}} Z{}^{A_{1}}{}_{b}\cdots Z{}^{A_\ell}{}_{a_\ell} \Omega_{A_{0}A_{1}\cdots A_\ell} \right) \\
		&= -nX^{A_{0}}Y^{A_{1}}Z{}^{A_{2}}{}_{a_{2}} \cdots Z{}^{A_\ell}{}_{a_\ell} \Omega_{A_{0}A_{1}\cdots A_\ell}\\
		&\mathrel{\phantom{=}}-\sum_{i=2}^{\ell} X{}^{A_{0}}Z^{A_{1}}{}_{a_i}Z^{A_{2}}{}_{a_{2}}\cdots Y^{A_i}\cdots Z^{A_\ell}{}_{a_\ell} \Omega_{A_{0}\cdots A_\ell}\\
		&= (n-(\ell-1))\sigma^{-1}I^{A_{0}}X^{A_{1}}Z^{A_{2}}{}_{a_{2}}\cdots Z^{A_\ell}{}_{a_\ell}\Omega_{A_{0}\cdots A_\ell},
	\end{align*}
	For the summed terms, we apply an even permutation of projectors, cycling $X^{A_0}Z^{A_1}_{a_i}Y^{A_i}\longrightarrow Y^{A_0}X^{A_1}Z^{A_i}_{a_i}$ and using $Y\wedge X = \sigma^{-1}I\wedge X$.

	Hence, $\sigma$ is a Killing scale if and only if the above quantity vanishes. When $I$ and $\Omega$ are parallel, then so is the completely skew tractor $I^{A_{0}}\Omega_{A_{0}A_{1}\cdots A_\ell}$. So by Lemma \ref{Lem:BGG}, it vanishes if and only if its top slot vanishes.
\end{proof}
An interesting property of conformal Killing--Yano forms is how they generate conformal Killing tensors. Let $\omega_{a_{1}\cdots a_\ell}\in \Gamma(\mathcal{E}_{[a_{1}\cdots a_\ell]}[\ell+1])$ be a conformal Killing--Yano form, and let $s{}_{bc}  = \omega_{b}{}^{a_{2}\cdots a_\ell}\omega{}_{ca_{2} \cdots a_\ell}\in \Gamma(\mathcal{E}_{(bc)}[4]) $. Then
\begin{align*}
\nabla{}_{(a} s{}_{bc)} 
&= 2\nabla{}_{(a} (\omega_{b}{}^{a_{2}\cdots a_\ell})\omega{}_{c)a_{2}\cdots a_\ell}\\
&=2\left(\boldsymbol{g}{}_{(ab} \rho^{a_{2}\cdots a_\ell}\omega{}_{c)a_{2}\cdots a_\ell}-(\ell-1)\rho{}_{(b}{}^{a_{3}\cdots a_\ell}\omega{}_{ac)a_{3}\cdots a_\ell}\right) \\
&= 2\boldsymbol{g}{}_{(ab} \rho^{a_{2}\cdots a_\ell}\omega{}_{c)a_{2}\cdots a_\ell},
\end{align*} 
since $d\omega$ and $\omega$ are both completely skew.
Therefore, $\nabla{}_{(a} s{}_{bd)_0} =0$ and $k{}_{bd}  = s{}_{(bd)_0} $ is a conformal Killing tensor. This also works for $\ell=1$ with no contraction; the symmetric trace-free product of two conformal Killing vectors is a conformal Killing tensor of rank $2$. When considering parallel tractor prolongations, this is also easy to see.
\begin{theorem}\label{Thm:KYProd}
	Let $\omega_{a_{1}\cdots a_\ell}\in \Gamma(\mathcal{E}{}_{[a_{1}\cdots a_\ell]}[\ell+1] )$ be a conformal Killing--Yano form and let $k{}_{bd}\in \Gamma(\mathcal{E}{}_{(bd)_0}[4] )$ be the trace-free part of $s{}_{bd}  = \omega{}_{(b}{}^{a_{2}\cdots a_\ell}\omega{}_{d)a_{2}\cdots a_\ell}$. Then $k{}_{bd} $ is a conformal Killing tensor. If $\omega$ is a normal solution and has full parallel prolongation $\Omega_{A_{0}\cdots A_\ell}\in \Gamma(\mathcal{E}_{[A_{0}\cdots A_\ell]})$, then the full prolongation $K{}_{ABCD} \in \Gamma(\begin{picture}(16, 8.0)(0, 2.0)
	\put(-0.2, 10){\line(1,0){10.4}}
	\put(-0.2, 5){\line(1,0){10.4}}
	\put(-0.2, 0){\line(1,0){10.4}}
	\put(0, 10){\line(0,-1){10}}
	\put(5, 10){\line(0,-1){10}}
	\put(10, 10){\line(0,-1){10}}
	\put(10, 0){${}_{\,0}$}
\end{picture}(\mathcal{E}_A))$ of $k{}_{bd} $ is the trace-free part of
	\[
		\tilde{K}_{ABCD} = \Omega{}_{AB}{}^{A_{2}\cdots A_\ell}\Omega{}_{CDA_{2}\cdots A_\ell}  -\Omega{}_{[AB}{}^{A_{2}\cdots A_\ell}\Omega{}_{CD]A_{2}\cdots A_\ell}\in \Gamma(\begin{picture}(10, 8.0)(0, 2.0)
	\put(-0.2, 10){\line(1,0){10.4}}
	\put(-0.2, 5){\line(1,0){10.4}}
	\put(-0.2, 0){\line(1,0){10.4}}
	\put(0, 10){\line(0,-1){10}}
	\put(5, 10){\line(0,-1){10}}
	\put(10, 10){\line(0,-1){10}}
\end{picture}(\mathcal{E}_A)).
	\] 
	Then $K{}_{ABCD} $ is parallel and $k{}_{bd}$ is a normal solution of the conformal Killing tensor equation.
\end{theorem}
\begin{proof}
	Let $K{}_{ABCD} $ be the trace-free part of 
	\[
		\Omega{}_{AB}{}^{A_{2}\cdots A_\ell}\Omega{}_{CDA_{2}\cdots A_\ell}  -\Omega{}_{[AB}{}^{A_{2}\cdots A_\ell}\Omega{}_{CD]A_{2}\cdots A_\ell}.
	\] 
	The tractor $\Omega{}_{AB}{}^{A_{2}\cdots A_\ell}\Omega{}_{CDA_{2}\cdots A_\ell}$ is a section of $S^2(\Lambda^2(\mathcal{E}_A))$. We have
	\[
	S^2(\Lambda^2(\mathcal{E}_A)) = \begin{picture}(16, 11.0)(0, 5.0)
	\put(-0.2, 16){\line(1,0){16.4}}
	\put(-0.2, 8){\line(1,0){16.4}}
	\put(-0.2, 0){\line(1,0){16.4}}
	\put(0, 16){\line(0,-1){16}}
	\put(8, 16){\line(0,-1){16}}
	\put(16, 16){\line(0,-1){16}}
\end{picture} \ (\mathcal{E}_A)\oplus \begin{picture}(8, 19.0)(0, 13.0)
	\put(-0.2, 32){\line(1,0){8.4}}
	\put(-0.2, 24){\line(1,0){8.4}}
	\put(-0.2, 16){\line(1,0){8.4}}
	\put(-0.2, 8){\line(1,0){8.4}}
	\put(-0.2, 0){\line(1,0){8.4}}
	\put(0, 32){\line(0,-1){32}}
	\put(8, 32){\line(0,-1){32}}\end{picture}\ (\mathcal{E}_A),
	\] 
	\\
	so removing the completely skew part yields a tractor with curvature symmetries. Removing trace yields a parallel tractor $K{}_{ABCD} $ with trace-free curvature symmetries.
	Then the top slot of $\tilde{K}{}_{ABCD} $ is 
	\begin{align*}
		X{}^{A} Z{}^{B}{}_{b} X{}^{C} Z{}^{D}{}_{d} \tilde{K}{}_{ABCD}
		&=X^A Z{}^{B}{}_{b} X^C Z{}^{D}{}_{d} \Omega{}_{AB}{}^{A_{2}\cdots A_\ell}\Omega{}_{CDA_{2}\cdots A_\ell} \\
		&= X^A Z{}^{B}{}_{b} X^C Z{}^{D}{}_{d}h^{A_{2}B_{2}}\cdots h^{A_\ell B_\ell} \Omega{}_{ABA_{2}\cdots A_\ell}\Omega{}_{CDB_{2}\cdots B_\ell}\\
		&= X^A Z{}^{B}{}_{b} X^C Z{}^{D}{}_{d}Z^{A_{2}a_{2}}Z^{B_{2}}{}_{a_{2}}\cdots Z^{A_\ell a_\ell} Z^{B_\ell}{}_{a_\ell} \Omega{}_{AB}{}^{A_{2}\cdots A_\ell}\Omega{}_{CDA_{2}\cdots A_\ell}\\
		&=\omega{}_{b}{}^{a_{2}\cdots a_\ell}\omega{}_{da_{2}\cdots a_\ell},
	\end{align*}
	since $X^A X^{A_i} \Omega{}_{ABA_{2}\cdots A_\ell} =0 $ and $X^AX^C$ applied to the totally skew part of $\tilde{K}{}_{ABCD} $ vanishes. The trace tractor which is then removed from $\tilde{K}$ has totally trace top slot, and its final top slot is trace-free, so it must be the trace free part of the final line above, which is $k{}_{bd} $.
	\end{proof}
	\section{Higher Rank Conformal Killing Tensors}\label{Sec:HO}

We now consider higher order conformal Killing tensors, that is, with $\ell\ge 3$. There are several things which make this case more complicated. We will discuss these challenges, while proving extensions of the results from Section \ref{Sec:KT} where possible.

The prolongation bundle for the conformal Killing tensor equation is the tractors with symmetries \begin{picture}(31, 8.0)(0, 2.0)
\put(-0.2, 10){\line(1,0){25.4}}
	\put(-0.2, 5){\line(1,0){25.4}}
	\put(-0.2, 0){\line(1,0){25.4}}
	\put(0, 10){\line(0,-1){10}}
	\put(5, 10){\line(0,-1){10}}
	\put(9, 2.3){${}_{\cdots}$}
	\put(9, 7.3){${}_{\cdots}$}
	\put(20, 10){\line(0,-1){10}}
	\put(25, 10){\line(0,-1){10}}
	\put(25, 0){${}_{\,0}$}
\end{picture}$(\mathcal{E}_A)$, equipped with the prolongation connection, as derived in \cite{hammerl_natural_2009}, which in the flat case is the standard tractor connection. We will use the skew symmetric realisation of this bundle, consisting of tractors satisfying
\[
	K{}_{A_{1}B_{1}\cdots A_\ell B_\ell} = K{}_{[A_{1}B_{1}]\cdots [A_\ell B_\ell]}\quad \text{and} \quad K{}_{\cdots [A_{i}B_{i}A_{i+1}]B_{i+1}\cdots}=0 
\] 
for all $i\in \{1\dots \ell-1\} $. This implies symmetry of adjacent pairs, $K{}_{\cdots A_i B_i A_{i+1} B_{i+1}\cdots } = K{}_{\cdots A{}_{i+1}  B{}_{i+1}  A_{i} B_{i}\cdots } $, and hence symmetry under any permutation of pairs. This follows in the same way as showing that tractors with curvature symmetries have symmetry of pairs. So this is a subbundle of $S^\ell(\Lambda^2(\mathcal{E}_A))$. As outlined in Section \ref{Sec:BGG}, the full prolongation of a tensor $k{}_{a_{1}\cdots a_\ell}\in \Gamma(\mathcal{E}{}_{(a_{1}\cdots a_\ell)_0}[2\ell] ) $ is a totally trace-free tractor $K{}_{A_{1}B_{1}\cdots A_\ell B_\ell} = L_0(k)$, characterised by
\[
\Pi(K) = X{}^{A_{1}}Z{}^{B_{1}}{}_{a_{1}}  \cdots X{}^{A_\ell} Z{}^{B_\ell}{}_{a_\ell}  K{}_{A_{1}B_{1}\cdots A_\ell B_\ell} = k{}_{a_{1}\cdots a_\ell}\quad \text{and} \quad \partial^*(\nabla_aK)=0 
\] 

We defined Killing scales analogously to the case of rank two conformal Killing tensors. Let $k{}_{a_{1}\cdots a_\ell}\in \Gamma(\mathcal{E}{}_{(a_{1}\cdots a_\ell)_0}[2\ell] ) $ be a conformal Killing tensor, let $\sigma\in \Gamma(\mathcal{E}[1])$ be a scale and let $g = \sigma^{-2}\boldsymbol{g}$ be the corresponding metric on $M\setminus\mathcal{Z}(\sigma)$. Then $g$ is a \textit{strong Killing scale} for $k$ if $k$ is a Killing tensor for $g$.
\begin{theorem}
	Let $k{}_{a_{1}\cdots a_\ell}\in \Gamma(\mathcal{E}{}_{(a_{1}\cdots a_\ell)_0}[2\ell] ) $ be a conformal Killing tensor and let $K{}_{A_{1}B_{1}\cdots A_\ell B_\ell} $ be its full prolongation. Let $I_A = D_A\sigma$ for some $\sigma\in \Gamma(\mathcal{E}[1])$. Then 
	\[
		I{}^{A_{1}}X{}^{B_{1}}X{}^{A_{2}}Z{}^{B_{2}}{}_{a_{2}}\cdots X{}^{A_{\ell}}Z{}^{B_\ell}{}_{a_\ell}K{}_{A_{1}B_{1}\cdots A_\ell B_\ell} = \frac{\sigma}{(n+2\ell-2)}\nabla_a k{}^{a}{}_{a_{2}\cdots a_\ell} 
	\] 
	Therefore, $\sigma$ is a strong Killing scale for $k$ if and only if
	\begin{equation}\label{Eqn:IXZK}
	I{}^{A_{1}}X^{B_{1}}X{}^{A_{2}} Z{}^{B_{2}}{}_{a_{2}}\cdots X{}^{A_\ell}Z{}^{B_\ell}{}_{a_\ell} K{}_{A_{1}B_{1}\cdots A_\ell B_\ell} = 0.
	\end{equation} 
	If $I$ and $K$ are parallel, then $\sigma$ is a strong Killing scale if and only if
	\[
	I{}^{A_{1}}K{}_{A_{1}B_{1}\cdots A_\ell B_\ell}=0  
	\] 
\end{theorem}
\begin{proof}
	We differentiate $k{}_{a_{1}\cdots a_\ell} = X{}^{A_{1}}Z{}^{B_{1}}{}_{a_1} \cdots X{}^{A_\ell} Z{}^{B_\ell}{}_{a_\ell} K{}_{A_{1}B_{1}\cdots A_\ell B_\ell} $. In the following, all omitted projectors are $X{}^{A_j}Z^{B_j}{}_{a_j} $.
	\begin{align*}
		\nabla{}^{a_{1}}  k{}_{a_{1}\cdots a_\ell} 
		&= \nabla{}^{a_{1}} \left(X{}^{A_{1}}Z{}^{B_{1}}{}_{a_{1}}  \cdots X{}^{A_\ell} Z{}^{B_\ell}{}_{a_\ell}  K{}_{A_{1}B_{1}\cdots A_\ell B_\ell} \right)\\ 
		&= -\sum_{i=1}^{\ell} \boldsymbol{g}{}_{a_i}{}^{a_{1}} \cdots X{}^{A_i} Y{}^{B_i}\cdots K{}_{A_{1}B_{1}\cdots A_\ell B_\ell} + \sum_{i=1}^{\ell} \cdots Z{}^{A_ia_{1}}Z{}^{B_i}{}_{a_i} \cdots K{}_{A_{1}B_{1}\cdots A_\ell B_\ell} \\
		&\mathrel{\phantom{=}} +X{}^{A_{1}} Z{}^{B_{1}}{}_{a_{1}} \cdots X{}^{A_{\ell}} Z{}^{B_\ell}{}_{a_\ell} \nabla{}^{a_{1}} K{}_{A_{1}B_{1}\cdots A_\ell B_\ell} \\
		&= -(n+\ell-1)X{}^{A_{1}} Y{}^{B_{1}} X{}^{A_{2}} Z{}^{B_{2}}{}_{a_{2}} \cdots X{}^{A_\ell} Z{}^{B_\ell}{}_{a_\ell}K{}_{A_{1}B_{1}\cdots A_\ell B_\ell} \\
		&\mathrel{\phantom{=}} -\sum_{i=2}^{\ell} X{}^{A_{1}} Y{}^{B_{1}} \cdots X{}^{A_{i}} Z{}^{B_i}{}_{a_i} \cdots K_{A_{1}B_{1}\cdots A_\ell B_\ell}\\
		&\mathrel{\phantom{=}} +X{}^{A_{1}} Z{}^{B_{1}}{}_{a_{1}} \cdots X{}^{A_{\ell}} Z{}^{B_\ell}{}_{a_\ell} \nabla{}^{a_{1}} K{}_{A_{1}B_{1}\cdots A_\ell B_\ell}\\
		&=\sigma^{-1}(n+2\ell-2)I{}^{A_{1}} X{}^{B_{1}} X{}^{A_{2}} Z{}^{B_{2}}{}_{a_{2}} \cdots X{}^{A_\ell} Z{}^{B_\ell}{}_{a_\ell}K{}_{A_{1}B_{1}\cdots A_\ell B_\ell} \\
		&\mathrel{\phantom{=}} +X{}^{A_{1}} Z{}^{B_{1}}{}_{a_{1}} \cdots X{}^{A_{\ell}} Z{}^{B_\ell}{}_{a_\ell} \nabla{}^{a_{1}} K{}_{A_{1}B_{1}\cdots A_\ell B_\ell},
	\end{align*}
	where we used $h{}_{AB}  = Z{}_{A}{}^{a} Z{}_{Ba} +X_AY_B+Y_AX_B$ and that $K$ is trace-free.
	So now all that remains is to show that $X{}^{A_{1}} Z{}^{B_{1}}{}_{a_{1}} \cdots X{}^{A_{\ell}} Z{}^{B_\ell}{}_{a_\ell} \nabla{}^{a_{1}} K{}_{A_{1}B_{1}\cdots A_\ell B_\ell}=0$
	Similar to the $\ell=2$ case, we can succinctly write $\partial^*\circ d^{\nabla}(K)$ as
	\[
	\partial^*\circ d^{\nabla}(K) = 2\ell\operatorname{Sym}_{A_iB_i}\left( X{}_{[B_{1}}Z{}^{Ca} \nabla{}_{a} K_{A_{1}]CA_{2}B_{2}\cdots A_\ell B_\ell} -Z{}_{[B_{1}}{}^a X^C \nabla_a K{}_{A_{1}]CA_{2}B_{2}\cdots A_\ell B_\ell}   \right),
	\] 
	by which we mean symmetrise over all skew pairs $A_iB_i$. Since $\partial^*\circ d^{\nabla}(K)=0$, 
	\[
	\begin{aligned}
		0&=X{}^{A_{1}} Y{}^{B_{1}} X{}^{A_{2}} Z{}^{B_{2}}{}_{a_{2}} \cdots X{}^{A_\ell}Z{}^{B_\ell}{}_{a_\ell}\partial^*\circ d^{\nabla}(K)_{A_{1}B_{1}\cdots A_\ell B_\ell}\\
		 &= X{}^{A_{1}}Z{}^{Ca} X{}^{A_{2}}Z{}^{B_{2}}{}_{a_{2}}\cdots X{}^{A_{\ell}}Z{}^{B_{\ell}}{}_{a_{\ell}}  \nabla_aK{}_{A_{1}CA_{2}B_{2}\cdots A_\ell B_\ell},
	\end{aligned}
	\] 
	which, after relabelling indices, is precisely the desired equality. 

	When $I$ and $K$ are parallel, the tractor $I{}^{A_{1}} K{}_{A_{1}B_{1}\cdots A_\ell B_\ell}$ is parallel and has symmetries given by the Young diagram \begin{picture}(31, 8.0)(0, 2.0)
	\put(-0.2, 10){\line(1,0){25.4}}
	\put(-0.2, 5){\line(1,0){25.4}}
	\put(-0.2, 0){\line(1,0){20.4}}
	\put(0, 10){\line(0,-1){10}}
	\put(5, 10){\line(0,-1){10}}
	\put(9, 2.3){${}_{\cdots}$}
	\put(9, 7.3){${}_{\cdots}$}
	\put(20, 10){\line(0,-1){10}}
	\put(25, 10){\line(0,-1){5}}
	\put(25, 0){${}_{\,0}$}
\end{picture}$(\mathcal{E}_A)$. The top slot is the left side of \eqref{Eqn:IXZK}, so $\sigma$ is a strong Killing scale if and only if $I{}^{A_{1}} K{}_{A_{1}B_{1}\cdots A_\ell B_\ell}=0$ by Lemma \ref{Lem:BGG}.
\end{proof}
The following appeared in \cite{gover_conformal_nodate} for Killing vectors.
\begin{proposition}
	Let $k_{a_{1}\cdots a_\ell}\in \Gamma(\mathcal{E}_{(a_{1}\cdots a_\ell)_0}[2\ell])$ be a conformal Killing tensor of rank $\ell$ and let $\sigma$ be a strong Killing scale. Then
	\[
		\mathcal{Z}(k_{a_{1}\cdots a_\ell})\subset \mathcal{Z}({\nabla}{}_{b}k{}^{b}{}_{a_{2}\cdots a_\ell}) \cup \mathcal{Z}(\sigma),
	\] 
	where ${\nabla}$ is the Levi-Civita connection of any metric in the conformal class.
\end{proposition}
\begin{proof}
	Let $p\in M$ be such that $\sigma(p)\neq 0$. Let $g\in \boldsymbol{c}$. Then, in terms of the projectors of $g$, we have $I_A = \sigma Y_A+(\nabla_a\sigma) Z{}_{A}{}^{a} -\frac{1}{n}(\Delta \sigma+J\sigma)X_A$. Then since $\sigma$ is a strong Killing scale,
	\begin{align*}
	0
	&= I{}^{A_{1}}X^{B_{1}}X{}^{A_{2}} Z{}^{B_{2}}{}_{a_{2}}\cdots X{}^{A_\ell}Z{}^{B_\ell}{}_{a_\ell} K{}_{A_{1}B_{1}\cdots A_\ell B_\ell} \\ 
	&= \sigma Y{}^{A_{1}}X^{B_{1}}X{}^{A_{2}} Z{}^{B_{2}}{}_{a_{2}}\cdots X{}^{A_\ell}Z{}^{B_\ell}{}_{a_\ell} K{}_{A_{1}B_{1}\cdots A_\ell B_\ell}\\
	&\mathrel{\phantom{=}}+(\nabla_a\sigma) Z{}^{A_{1}a}X^{B_{1}}X{}^{A_{2}} Z{}^{B_{2}}{}_{a_{2}}\cdots X{}^{A_\ell}Z{}^{B_\ell}{}_{a_\ell} K{}_{A_{1}B_{1}\cdots A_\ell B_\ell}\\
	&= \frac{\sigma}{n+2\ell-2}\nabla_a k{}^{a}{}_{a_{2}\cdots a_\ell}-(\nabla_a\sigma) k^{a}{}_{a_{2}\cdots a_\ell}
	\end{align*} 
	Hence, if $\sigma(p)\neq 0$ and $k|_p=0$, then $\nabla_b k{}^{b}{}_{a_{2}\cdots a_\ell} =0$. 
\end{proof}

A metric $g\in \boldsymbol{c}$ is a \textit{Killing scale} for a conformal Killing tensor $k_{a_{1}\cdots a_\ell}$ if there is a Killing tensor $r{}_{a_{1}\cdots a_\ell}\in \Gamma(\mathcal{E}{}_{(a_{1}\cdots a_\ell)} [2\ell]) $ for $g$ such that
\[
\nabla{}_{(a_{0}}r{}_{a_{1}\cdots a_\ell)}=0
\] 
and $r{}_{(a_{1}\cdots a_\ell)_0}=k{}_{a_{1}\cdots a_\ell}$. The last condition is equivalent to $r{}_{a_{1}\cdots a_\ell} = k{}_{a_{1}\cdots a_\ell}+\boldsymbol{g}{}_{(a_{1}a_{2}}\eta{}_{a_{3}\cdots a_\ell)}$ for some section $\eta_{a_{3}\cdots a_\ell}\in \Gamma(\mathcal{E}_{(a_{3}\cdots a_\ell)}[2\ell-2])$. The Killing condition for $r$ then becomes
\[
\nabla{}_{(a_{0}}r{}_{a_{1}\cdots a_\ell)} = \nabla{}_{(a_{0}} k{}_{a_{1}\cdots a_\ell)}+\boldsymbol{g}_{(a_{1}a_{2}}\nabla{}_{a_{0}}\eta{}_{a_{3}\cdots a_\ell)} = 0,
\] 
which then becomes
\[
	\nabla{}_{(a_{2}}\eta{}_{a_{3}\cdots a_\ell)}  = \rho_{a_{2}\cdots a_\ell} = -\frac{\ell}{n+2\ell-2}\nabla_a k{}^{a}{}_{a_{2}\cdots a_\ell}. 
\]
In the case $\ell=2$, this reduces to whether the form $\nabla_a k{}^{a}{}_{a_{2}} $ is exact. Then when $M$ is simply connected, or on a simply connected open subset of $M$, existence of $\eta$ is equivalent to $\rho$ being closed. Consider the case $\ell=3$. In this case, we want to know if $\rho{}_{ab}\in \Gamma(\mathcal{E}_{(ab)_0}[4])$ is in the image of the operator $\eta{}_{b} \mapsto \nabla{}_{(a} \eta{}_{b)} $. We use $\sigma$ to trivialise all density bundles. Then $\sigma$ is a Killing scale for $k{}_{abc} $ if and only if $\rho{}_{bc}  = -\frac{3}{n+4}\nabla_a k{}^{a}{}_{bc}$ is in the image of the projective Killing operator on one forms. The BGG sequence for the projective Killing operator is
\[
\begin{tikzcd}[column sep = 1.8em]
	0 \arrow[r]&\Gamma(\mathcal{E}_a(2))\arrow[r]&\Gamma(\mathcal{E}{}_{(ab)}(2) )\arrow[r]&\Gamma(\begin{picture}(16, 8.0)(-1, 2.0)
	\put(-0.2, 10){\line(1,0){10.4}}
	\put(-0.2, 5){\line(1,0){10.4}}
	\put(-0.2, 0){\line(1,0){10.4}}
	\put(0, 10){\line(0,-1){10}}
	\put(5, 10){\line(0,-1){10}}
	\put(10, 10){\line(0,-1){10}}
	\end{picture}(\mathcal{E}_a)(2))\arrow[r]&\Gamma\left(\begin{picture}(16, 10.5)(-1, 4.5)
	\put(-0.2, 15){\line(1,0){10.4}}
	\put(-0.2, 10){\line(1,0){10.4}}
	\put(-0.2, 5){\line(1,0){10.4}}
	\put(-0.2, 0){\line(1,0){5.4}}
	\put(0, 15){\line(0,-1){15}}
	\put(5, 15){\line(0,-1){15}}
	\put(10, 15){\line(0,-1){10}}
\end{picture}(\mathcal{E}_a)(2)\right)\arrow[r]&\cdots,
\end{tikzcd}	
\] 
which is the Calabi complex on Riemannian manifolds. The sequence is only a complex when $(M,[\nabla^g])$ is projectively flat. In the flat case, moreover, the sequence is exact on contractible subsets of $M$. So the image of the projective Killing operator is the kernel of the operator $\mathscr{D}_1 \colon \Gamma(\mathcal{E}{}_{(ab)}(2) )\to \Gamma\left(\begin{picture}(16, 8.0)(-1, 2.0)
	\put(-0.2, 10){\line(1,0){10.4}}
	\put(-0.2, 5){\line(1,0){10.4}}
	\put(-0.2, 0){\line(1,0){10.4}}
	\put(0, 10){\line(0,-1){10}}
	\put(5, 10){\line(0,-1){10}}
	\put(10, 10){\line(0,-1){10}}
\end{picture}(\mathcal{E}_a)(2)\right) $. For higher rank $\ell$, this easily extends to the following proposition.

\begin{proposition}\label{Prp:HR}
	Let $\ell\ge 3$. Let $k{}_{a_{1}\cdots a_\ell} \in \Gamma(\mathcal{E}_{(a_{1}\cdots a_\ell)_0}[2\ell])$ be a conformal Killing tensor and let $\sigma\in \Gamma(\mathcal{E}[1])$. Then $\sigma$ is a Killing scale for $k$ if and only if, after using $\sigma$ to trivialise densities, $\nabla^g_ak{}^a{}_{a_{2}\cdots a_\ell} $ is in the image of the rank $\ell-2$ projective Killing operator.

	When $(M,\boldsymbol{p}=[\nabla^g])$ is projectively flat, $\sigma$ is a local Killing scale if and only if $\mathscr{D}_1(\nabla^g_ak{}^a{}_{a_{2}\cdots a_\ell})=0$, where $\mathscr{D}_1 \colon \Gamma(\mathcal{E}_{(a_{2}\cdots a_\ell)}(2\ell-4))\to \Gamma(\begin{picture}(31, 8.0)(0, 2.0)
	\put(-0.2, 10){\line(1,0){25.4}}
	\put(-0.2, 5){\line(1,0){25.4}}
	\put(-0.2, 0){\line(1,0){25.4}}
	\put(0, 10){\line(0,-1){10}}
	\put(5, 10){\line(0,-1){10}}
	\put(9, 2.3){${}_{\cdots}$}
	\put(9, 7.3){${}_{\cdots}$}
	\put(20, 10){\line(0,-1){10}}
	\put(25, 10){\line(0,-1){10}}
\end{picture}(\mathcal{E}_a)(2\ell-4)) $, where each row has $\ell-1$ boxes and $\mathscr{D}{}_{1}$ is a BGG operator from the sequence including the rank $\ell-2$ projective Killing operator.
\end{proposition}

So for higher rank we have an extra fundamental reliance on conformal flatness. For a curved projective structure, local exactness of the flat BGG sequence is no longer available, so the operator $\mathscr{D}_1$ does not give a differential test for membership in the image. For an Einstein metric \(g\), projective flatness is equivalent to conformal flatness of \((M,\boldsymbol c)\). Thus, in an Einstein scale, this route to a higher-rank analogue of Theorem~\ref{Thm:LKS} requires conformal flatness.

\begin{lemma}\label{Lem:HRtf}
	Let $R{}_{A_{1}B_{1}\cdots A_\ell B_\ell}\in  \Gamma(\begin{picture}(31, 8.0)(0, 2.0)
	\put(-0.2, 10){\line(1,0){25.4}}
	\put(-0.2, 5){\line(1,0){25.4}}
	\put(-0.2, 0){\line(1,0){25.4}}
	\put(0, 10){\line(0,-1){10}}
	\put(5, 10){\line(0,-1){10}}
	\put(9, 2.3){${}_{\cdots}$}
	\put(9, 7.3){${}_{\cdots}$}
	\put(20, 10){\line(0,-1){10}}
	\put(25, 10){\line(0,-1){10}}
\end{picture}(\mathcal{E}_a))$ have top slot
\[
	r{}_{a_{1}\cdots a{}_{\ell} } = X{}^{A_{1}} Z{}^{B_{1}}{}_{a_{1}}\cdots X{}^{A_\ell} Z{}^{B_\ell}{}_{a_\ell} R_{A_{1}B_{1}\cdots A_\ell B_\ell}.
\]
Then let $K_{A_{1}B_{1}\cdots A_\ell B_\ell}$ be the trace free part of $R_{A_{1}B_{1}\cdots A_\ell B_\ell}$ with top slot
\[
	k{}_{a_{1}\cdots a{}_{\ell} } = X{}^{A_{1}} Z{}^{B_{1}}{}_{a_{1}}\cdots X{}^{A_\ell} Z{}^{B_\ell}{}_{a_\ell} K_{A_{1}B_{1}\cdots A_\ell B_\ell}.
\]
Then $k{}_{a_{1}\cdots a{}_{\ell} } = r{}_{(a_{1}\cdots a_\ell)_0} $.
\end{lemma}
\begin{proof}
	With $K$ the trace-free part of $R$, every term in the difference $R-K$ involves $h{}_{AB} $. Taking the top slot of $R-K$ hence gives a tensor where every term involves $\boldsymbol{g}{}_{ab} $. Therefore, $k$ differs from $r$ only by trace terms and $k$ is trace-free. Hence $k{}_{a_{1}\cdots a{}_{\ell} } = r{}_{(a_{1}\cdots a_\ell)_0} $.
\end{proof}

\begin{proposition}\label{Prp:HRKS}
	Let $(M,\boldsymbol{c})$ be conformally flat. Let $k{}_{a_{1}\cdots a_\ell}\in \Gamma(\mathcal{E}{}_{(a_{1}\cdots a_\ell)_0}[2\ell] ) $ be a conformal Killing tensor and let $K{}_{A_{1} B_{1}\cdots A_\ell B_\ell} $ be its full prolongation. Let $I_A = D_A\sigma$ for some $\sigma\in \Gamma(\mathcal{E}[1])$ be a parallel scale tractor. 

	Then $\sigma$ is a Killing scale for $k$ if and only if $K$ is the trace-free part of a parallel tractor $R$ defined on $M\setminus\mathcal{Z}(\sigma)$, with symmetries $\begin{picture}(31, 8.0)(0, 2.0)
	\put(-0.2, 10){\line(1,0){25.4}}
	\put(-0.2, 5){\line(1,0){25.4}}
	\put(-0.2, 0){\line(1,0){25.4}}
	\put(0, 10){\line(0,-1){10}}
	\put(5, 10){\line(0,-1){10}}
	\put(9, 2.3){${}_{\cdots}$}
	\put(9, 7.3){${}_{\cdots}$}
	\put(20, 10){\line(0,-1){10}}
	\put(25, 10){\line(0,-1){10}}
\end{picture}(I^{\perp})$.
\end{proposition}
\begin{proof}
	Let $g = \sigma^{-2}\boldsymbol{g}$ be the Einstein metric defined on $M\setminus\mathcal{Z}(\sigma)$. Since $g$ is conformally flat and Einstein, it is projectively flat. Then, by Theorem \ref{Thm:BGGL}, rank $\ell$ Killing tensors for $g = \sigma^{-2}\boldsymbol{g}$ are in linear correspondence with parallel sections of $\begin{picture}(31, 8.0)(0, 2.0)
	\put(-0.2, 10){\line(1,0){25.4}}
	\put(-0.2, 5){\line(1,0){25.4}}
	\put(-0.2, 0){\line(1,0){25.4}}
	\put(0, 10){\line(0,-1){10}}
	\put(5, 10){\line(0,-1){10}}
	\put(9, 2.3){${}_{\cdots}$}
	\put(9, 7.3){${}_{\cdots}$}
	\put(20, 10){\line(0,-1){10}}
	\put(25, 10){\line(0,-1){10}}
\end{picture}(\mathcal{E}_\alpha)$, where $\mathcal{E}{}_{\alpha} $ is the projective tractor bundle of $(M\setminus\mathcal{Z}(\sigma),[\nabla^g])$. By Theorem \ref{Thm:Iperp}, there is a connection preserving isomorphism $\mathcal{E}_\alpha\cong I^{\perp}$. So rank $\ell$ Killing tensors for $g = \sigma^{-2}\boldsymbol{g}$ are in linear correspondence with parallel sections of $\begin{picture}(31, 8.0)(0, 2.0)
	\put(-0.2, 10){\line(1,0){25.4}}
	\put(-0.2, 5){\line(1,0){25.4}}
	\put(-0.2, 0){\line(1,0){25.4}}
	\put(0, 10){\line(0,-1){10}}
	\put(5, 10){\line(0,-1){10}}
	\put(9, 2.3){${}_{\cdots}$}
	\put(9, 7.3){${}_{\cdots}$}
	\put(20, 10){\line(0,-1){10}}
	\put(25, 10){\line(0,-1){10}}
\end{picture}(I^{\perp})$. Again, by Theorem \ref{Thm:BGGL}, there is a linear bijection between the space of conformal Killing tensors of rank $\ell$ and parallel sections of \begin{picture}(31, 8.0)(0, 2.0)
\put(-0.2, 10){\line(1,0){25.4}}
	\put(-0.2, 5){\line(1,0){25.4}}
	\put(-0.2, 0){\line(1,0){25.4}}
	\put(0, 10){\line(0,-1){10}}
	\put(5, 10){\line(0,-1){10}}
	\put(9, 2.3){${}_{\cdots}$}
	\put(9, 7.3){${}_{\cdots}$}
	\put(20, 10){\line(0,-1){10}}
	\put(25, 10){\line(0,-1){10}}
	\put(25, 0){${}_{\,0}$}
\end{picture}$(\mathcal{E}_A)$. Then, by Lemma \ref{Lem:HRtf}, the conformal Killing tensors for which $\sigma$ is a Killing scale are those whose prolongation is the trace-free part of some $R\in \begin{picture}(31, 8.0)(0, 2.0)
	\put(-0.2, 10){\line(1,0){25.4}}
	\put(-0.2, 5){\line(1,0){25.4}}
	\put(-0.2, 0){\line(1,0){25.4}}
	\put(0, 10){\line(0,-1){10}}
	\put(5, 10){\line(0,-1){10}}
	\put(9, 2.3){${}_{\cdots}$}
	\put(9, 7.3){${}_{\cdots}$}
	\put(20, 10){\line(0,-1){10}}
	\put(25, 10){\line(0,-1){10}}
\end{picture}(I^{\perp})$.
\end{proof}

Identifying the image of 
	\begin{equation}\label{Eqn:tfIpp}
\operatorname{trace-free}\colon \Gamma(\begin{picture}(31, 8.0)(0, 2.0)
	\put(-0.2, 10){\line(1,0){25.4}}
	\put(-0.2, 5){\line(1,0){25.4}}
	\put(-0.2, 0){\line(1,0){25.4}}
	\put(0, 10){\line(0,-1){10}}
	\put(5, 10){\line(0,-1){10}}
	\put(9, 2.3){${}_{\cdots}$}
	\put(9, 7.3){${}_{\cdots}$}
	\put(20, 10){\line(0,-1){10}}
	\put(25, 10){\line(0,-1){10}}
\end{picture}(I^{\perp}))\to  \Gamma(\begin{picture}(31, 8.0)(0, 2.0)
\put(-0.2, 10){\line(1,0){25.4}}
	\put(-0.2, 5){\line(1,0){25.4}}
	\put(-0.2, 0){\line(1,0){25.4}}
	\put(0, 10){\line(0,-1){10}}
	\put(5, 10){\line(0,-1){10}}
	\put(9, 2.3){${}_{\cdots}$}
	\put(9, 7.3){${}_{\cdots}$}
	\put(20, 10){\line(0,-1){10}}
	\put(25, 10){\line(0,-1){10}}
	\put(25, 0){${}_{\,0}$}
\end{picture}(\mathcal{E}_A)).
\end{equation}
 is a matter for multilinear algebra. If $U = U_{1}\oplus U_{2}$ is a direct sum of vector spaces, then
 \begin{equation}\label{Eqn:WDec}
	\Lambda^k(U_{1}\oplus U_{2}) =  \bigoplus_{i=0}^k (\Lambda^i U_1)\wedge (\Lambda^{k-i} U_2)
\end{equation} 
and
\begin{equation}\label{Eqn:SDec}
	S^k(U_{1}\oplus U_{2}) =  \bigoplus_{i=0}^k (S^iU_{1})\vee  (S^{k-i}U_{2}).
\end{equation} 
If $U$ is equipped with a non-degenerate bilinear form and the decomposition $U = U_{1}\oplus  U_{2}$ is orthogonal, then the direct sums above are also orthogonal.

Fix a point $p\in M$ and let $V = \mathcal{E}_A|_p$. Let $I\in V$ be non-null. Then we have $V = W\oplus I $, where $I$ is shorthand for the subspace spanned by $I$ and $W = I^{\perp}$. Then $\Lambda^2(V) = \Lambda^2W\oplus W\wedge I $ since $I\wedge I=0$. We then apply \eqref{Eqn:SDec} to obtain
\begin{equation}\label{Eqn:Sdecomp}
	S^\ell(\Lambda^2(V)) = \bigoplus_{i=0}^\ell S^i(W\wedge I)\vee S^{\ell-i}(\Lambda^2 W)
\end{equation} 
as an orthogonal direct sum. So if $K\in \Gamma(\begin{picture}(31, 8.0)(0, 2.0)
\put(-0.2, 10){\line(1,0){25.4}}
	\put(-0.2, 5){\line(1,0){25.4}}
	\put(-0.2, 0){\line(1,0){25.4}}
	\put(0, 10){\line(0,-1){10}}
	\put(5, 10){\line(0,-1){10}}
	\put(9, 2.3){${}_{\cdots}$}
	\put(9, 7.3){${}_{\cdots}$}
	\put(20, 10){\line(0,-1){10}}
	\put(25, 10){\line(0,-1){10}}
	\put(25, 0){${}_{\,0}$}
\end{picture}(\mathcal{E}_A|_p))\subset \Gamma(S^\ell(\Lambda^2(V)))$, then we can write $K = K^{0}+K^{1}+\cdots K^\ell$, according to the direct sum \eqref{Eqn:Sdecomp}. We define $\pi^0 \colon S^\ell(\Lambda^2(V))\to S^\ell(\Lambda^2(W)) $ by $\pi^0(K) = K^0$. Since \eqref{Eqn:Sdecomp} is an orthogonal direct sum, the adjoint map to $\pi{}^{0} $ is the inclusion $\iota_0\colon S^\ell(\Lambda^2(W))\to S^\ell(\Lambda^2(V))$. As we show in the following lemma, a trace-free $K$ is determined by $K{}^{0} $ and $K{}^{1} $ alone.

\begin{lemma}\label{Lem:K0K1}
	Let $V$ be a vector space (with lower indices $A,B,C,\dots $) equipped with a non-degenerate bilinear form $h_{AB}$. Let $K_{A_{1}B_{1}\cdots A_\ell B_\ell}\in \Gamma(\begin{picture}(31, 8.0)(0, 2.0)
\put(-0.2, 10){\line(1,0){25.4}}
	\put(-0.2, 5){\line(1,0){25.4}}
	\put(-0.2, 0){\line(1,0){25.4}}
	\put(0, 10){\line(0,-1){10}}
	\put(5, 10){\line(0,-1){10}}
	\put(9, 2.3){${}_{\cdots}$}
	\put(9, 7.3){${}_{\cdots}$}
	\put(20, 10){\line(0,-1){10}}
	\put(25, 10){\line(0,-1){10}}
	\put(25, 0){${}_{\,0}$}
\end{picture}(V))$. Let $K = K^{0}+K^{1}+\cdots K^\ell$, according to the orthogonal direct sum \eqref{Eqn:Sdecomp}. Then $K=0$ if and only if $K^{0}=K^{1}=0$. Moreover,
\begin{equation}\label{Eqn:YDDec}
	\begin{picture}(31, 8.0)(0, 2.0)
	\put(-0.2, 10){\line(1,0){25.4}}
	\put(-0.2, 5){\line(1,0){25.4}}
	\put(-0.2, 0){\line(1,0){25.4}}
	\put(0, 10){\line(0,-1){10}}
	\put(5, 10){\line(0,-1){10}}
	\put(9, 2.3){${}_{\cdots}$}
	\put(9, 7.3){${}_{\cdots}$}
	\put(20, 10){\line(0,-1){10}}
	\put(25, 10){\line(0,-1){10}}
	\put(25, 0){${}_{\,0}$}
\end{picture}(\mathcal{E}_A) = \{K\in \begin{picture}(31, 8.0)(0, 2.0)
\put(-0.2, 10){\line(1,0){25.4}}
	\put(-0.2, 5){\line(1,0){25.4}}
	\put(-0.2, 0){\line(1,0){25.4}}
	\put(0, 10){\line(0,-1){10}}
	\put(5, 10){\line(0,-1){10}}
	\put(9, 2.3){${}_{\cdots}$}
	\put(9, 7.3){${}_{\cdots}$}
	\put(20, 10){\line(0,-1){10}}
	\put(25, 10){\line(0,-1){10}}
	\put(25, 0){${}_{\,0}$}
\end{picture}(\mathcal{E}_A):K^0=0\}\oplus \{K\in \begin{picture}(31, 8.0)(0, 2.0)
\put(-0.2, 10){\line(1,0){25.4}}
	\put(-0.2, 5){\line(1,0){25.4}}
	\put(-0.2, 0){\line(1,0){25.4}}
	\put(0, 10){\line(0,-1){10}}
	\put(5, 10){\line(0,-1){10}}
	\put(9, 2.3){${}_{\cdots}$}
	\put(9, 7.3){${}_{\cdots}$}
	\put(20, 10){\line(0,-1){10}}
	\put(25, 10){\line(0,-1){10}}
	\put(25, 0){${}_{\,0}$}
\end{picture}(\mathcal{E}_A):K^1=0\}  
\end{equation} 
as an orthogonal decomposition.\end{lemma}
\begin{proof}
	Let $\hat{h}{}_{AB}  = h{}_{AB} - |I|^{-2}I_AI_B$. Then $\hat{h}{}_{AB}  I^A=0$ and $\hat{h}{}_{AB} $ is the metric $h$ restricted to $W$. Since $K$ is trace-free, we have
	\begin{equation}\label{Eqn:K0K1:1}
		-\hat{h}^{A_{1}A_{2}}K_{A_{1}B_{1}A_{2}B_{2}\cdots A_\ell B_\ell} = |I|^{-2} I^{A_{1}}I^{A_{2}}K_{A_{1}B_{1}A_{2}B_{2}\cdots A_\ell B_\ell}.
	\end{equation} 
	Suppose that $K\neq 0$ and that $K^0=0$ and let $1\le j\le \ell$ be the smallest number such that $K^j\neq 0$. Let $U,V\in \Gamma(W)$. Then 
	\begin{equation}\label{Eqn:K0K1:2}
		-\hat{h}^{A_{1}A_{2}}U^{B_{1}}V{}^{B_{2}} K_{A_{1}B_{1}A_{2}B_{2}\cdots A_\ell B_\ell} = -\sum_{i=j}^{\ell-2} \hat{h}^{A_{1}A_{2}}U^{B_{1}}V{}^{B_{2}} K^i{}_{A_{1}B_{1}A_{2}B_{2}\cdots A_\ell B_\ell},
	\end{equation}
	since $K{}^{\ell}$ and $K^{\ell-1}$ both vanish after insertion of two elements of $\Lambda^2(W)$. Then, by \eqref{Eqn:K0K1:1}, the above is then equal to
	\begin{equation}\label{Eqn:K0K1:3}
		|I|^{-2}I{}^{A_{1}}I{}^{A_{2}}U{}^{B_{1}}V{}^{B_{2}}K_{A_{1}B_{1}A_{2}B_{2}\cdots A_\ell B_\ell} = |I|^{-2}\sum_{i=\max(j,2)}^{\ell} I{}^{A_{1}}I{}^{A_{2}}U{}^{B_{1}}V{}^{B_{2}}K^i_{A_{1}B_{1}A_{2}B_{2}\cdots A_\ell B_\ell}.
	\end{equation} 
	We will now fill all remaining pairs of indices $A_3B_3,\dots ,A_\ell B_\ell$. We use $\min(\ell-2,\ell-j)$ elements of $\Lambda^2W$ and then fill the remaining $\max(0,j-2)$ pairs of indices with elements of $W\wedge I$. In \eqref{Eqn:K0K1:2}, the result is zero, since the projectors only see $K^{\max(0,j-2)}$, which is zero by assumption. The result in \eqref{Eqn:K0K1:3} is $|I|^{-2}K^{\max(2,j)}$, contracted with generic elements of $W\wedge I$ and $\Lambda^2 W$. So by equality of \eqref{Eqn:K0K1:2} and \eqref{Eqn:K0K1:3}, we conclude that $K^{\max(j,2)}=0$. Since $K^j\neq 0$, we conclude that $j<2$, so $j=1$. Hence, if $K\neq 0$, then either $K^{1}$ or $K^0$ is nonzero.

	Consider the transformation $\tau\in GL(V)$ given by $\tau{}_{A}{}^{B} =\delta{}_{A}{}^{B} -2|I|^{-2}I_AI^B$. Then the $+1$ eigenspace of $\tau$ is $W$ and the $-1$ eigenspace of $\tau$ is $I$. From this it is clear that $\tau\in O(\mathcal{E}_A,h)$. So $\tau$ acts on $\Lambda^2(V)$, $S^\ell(\Lambda^2(V))$ and preserves $\begin{picture}(31, 8.0)(0, 2.0)
\put(-0.2, 10){\line(1,0){25.4}}
	\put(-0.2, 5){\line(1,0){25.4}}
	\put(-0.2, 0){\line(1,0){25.4}}
	\put(0, 10){\line(0,-1){10}}
	\put(5, 10){\line(0,-1){10}}
	\put(9, 2.3){${}_{\cdots}$}
	\put(9, 7.3){${}_{\cdots}$}
	\put(20, 10){\line(0,-1){10}}
	\put(25, 10){\line(0,-1){10}}
	\put(25, 0){${}_{\,0}$}
\end{picture}(V)$.

	Now consider the action of $\tau$ on $\Lambda^2(\mathcal{E}_A)$. The $+1$ eigenspace of $\tau$ inside $\Lambda^2(\mathcal{E}_A)$ is $W\wedge W$ and the $-1$ eigenspace is $W\wedge I$. From this it follows that
	\[
		\tau(K) = \tau(K^0+K^1+\cdots K^\ell) = K^0-K^1+K^2-\cdots = K^+-K^-,
	\] 
	that is, the $+1$ eigenspace of $\tau$ acting on $S^\ell(\Lambda^2(\mathcal{E}_A))$ is the even summands in \eqref{Eqn:Sdecomp}, and the $-1$ eigenspace of $\tau$ is the odd summands. 

	Since $\tau\in O(V,h)$, it follows that $K^+,K^-\in \begin{picture}(31, 8.0)(0, 2.0)
\put(-0.2, 10){\line(1,0){25.4}}
	\put(-0.2, 5){\line(1,0){25.4}}
	\put(-0.2, 0){\line(1,0){25.4}}
	\put(0, 10){\line(0,-1){10}}
	\put(5, 10){\line(0,-1){10}}
	\put(9, 2.3){${}_{\cdots}$}
	\put(9, 7.3){${}_{\cdots}$}
	\put(20, 10){\line(0,-1){10}}
	\put(25, 10){\line(0,-1){10}}
	\put(25, 0){${}_{\,0}$}
\end{picture}(V)$. Then it is clear that $K^-\in \{K\in \begin{picture}(31, 8.0)(0, 2.0)
\put(-0.2, 10){\line(1,0){25.4}}
	\put(-0.2, 5){\line(1,0){25.4}}
	\put(-0.2, 0){\line(1,0){25.4}}
	\put(0, 10){\line(0,-1){10}}
	\put(5, 10){\line(0,-1){10}}
	\put(9, 2.3){${}_{\cdots}$}
	\put(9, 7.3){${}_{\cdots}$}
	\put(20, 10){\line(0,-1){10}}
	\put(25, 10){\line(0,-1){10}}
	\put(25, 0){${}_{\,0}$}
\end{picture}(\mathcal{E}_A):K^0=0\}$ and that $K^+\in\{K\in \begin{picture}(31, 8.0)(0, 2.0)
\put(-0.2, 10){\line(1,0){25.4}}
	\put(-0.2, 5){\line(1,0){25.4}}
	\put(-0.2, 0){\line(1,0){25.4}}
	\put(0, 10){\line(0,-1){10}}
	\put(5, 10){\line(0,-1){10}}
	\put(9, 2.3){${}_{\cdots}$}
	\put(9, 7.3){${}_{\cdots}$}
	\put(20, 10){\line(0,-1){10}}
	\put(25, 10){\line(0,-1){10}}
	\put(25, 0){${}_{\,0}$}
\end{picture}(\mathcal{E}_A):K^0=1\}$. Hence, \eqref{Eqn:YDDec} is a direct sum, with the first summand the $-1$ eigenspace of $\tau$ and the second summand the $+1$ eigenspace of $\tau$, which is an orthogonal direct sum.
\end{proof}
The result of our algebraic analysis is the following theorem.

\begin{theorem}\label{Thm:MLA}
	Let $V$ be a vector space (with lower indices $A,B,C,\dots $) equipped with a non-degenerate bilinear form $h_{AB}$. Let $I_A\in V^*$ be non-null and let $W = I^{\perp}$. Let $K_{A_{1}B_{1}\cdots A_\ell B_\ell}\in \Gamma(\begin{picture}(31, 8.0)(0, 2.0)
\put(-0.2, 10){\line(1,0){25.4}}
	\put(-0.2, 5){\line(1,0){25.4}}
	\put(-0.2, 0){\line(1,0){25.4}}
	\put(0, 10){\line(0,-1){10}}
	\put(5, 10){\line(0,-1){10}}
	\put(9, 2.3){${}_{\cdots}$}
	\put(9, 7.3){${}_{\cdots}$}
	\put(20, 10){\line(0,-1){10}}
	\put(25, 10){\line(0,-1){10}}
	\put(25, 0){${}_{\,0}$}
\end{picture}(V))$. Let $K = K^{0}+K^{1}+\cdots K^\ell$, according to the orthogonal direct sum \eqref{Eqn:Sdecomp}. Then the following are equivalent
\begin{tabel}
	\item\label{Thm:imK:1} $K$ is in the image of 
		\begin{equation}\label{Eqn:tfWV}	\operatorname{trace-free}\colon \Gamma(\begin{picture}(26, 8.0)(0, 2.0)
	\put(-0.2, 10){\line(1,0){25.4}}
	\put(-0.2, 5){\line(1,0){25.4}}
	\put(-0.2, 0){\line(1,0){25.4}}
	\put(0, 10){\line(0,-1){10}}
	\put(5, 10){\line(0,-1){10}}
	\put(9, 2.3){${}_{\cdots}$}
	\put(9, 7.3){${}_{\cdots}$}
	\put(20, 10){\line(0,-1){10}}
	\put(25, 10){\line(0,-1){10}}
\end{picture}(W))\to  \Gamma(\begin{picture}(31, 8.0)(0, 2.0)
\put(-0.2, 10){\line(1,0){25.4}}
	\put(-0.2, 5){\line(1,0){25.4}}
	\put(-0.2, 0){\line(1,0){25.4}}
	\put(0, 10){\line(0,-1){10}}
	\put(5, 10){\line(0,-1){10}}
	\put(9, 2.3){${}_{\cdots}$}
	\put(9, 7.3){${}_{\cdots}$}
	\put(20, 10){\line(0,-1){10}}
	\put(25, 10){\line(0,-1){10}}
	\put(25, 0){${}_{\,0}$}
\end{picture}(V)).
\end{equation}
\item \label{Thm:imK:2}$K^1=0$.
\item \label{Thm:imK:3}After skewing over indices with the same subscript,
	\[
		I^{A_{1}}K_{A_{1}B_{1}A_{2}B_{2}\cdots A_\ell B_\ell} I_{C_{2}}\cdots I_{C_\ell}=0.
	\] 
	That is, each set of indices $A_iB_iC_i$ for $i\ge 2$ are completely skew.
\end{tabel}
\end{theorem}
\begin{proof}
	Recall the orthogonal decomposition 
\[
\begin{picture}(31, 8.0)(0, 2.0)
\put(-0.2, 10){\line(1,0){25.4}}
	\put(-0.2, 5){\line(1,0){25.4}}
	\put(-0.2, 0){\line(1,0){25.4}}
	\put(0, 10){\line(0,-1){10}}
	\put(5, 10){\line(0,-1){10}}
	\put(9, 2.3){${}_{\cdots}$}
	\put(9, 7.3){${}_{\cdots}$}
	\put(20, 10){\line(0,-1){10}}
	\put(25, 10){\line(0,-1){10}}
\end{picture}(V) = \begin{picture}(31, 8.0)(0, 2.0)
\put(-0.2, 10){\line(1,0){25.4}}
	\put(-0.2, 5){\line(1,0){25.4}}
	\put(-0.2, 0){\line(1,0){25.4}}
	\put(0, 10){\line(0,-1){10}}
	\put(5, 10){\line(0,-1){10}}
	\put(9, 2.3){${}_{\cdots}$}
	\put(9, 7.3){${}_{\cdots}$}
	\put(20, 10){\line(0,-1){10}}
	\put(25, 10){\line(0,-1){10}}
	\put(25, 0){${}_{\,0}$}
\end{picture}(V)\oplus \begin{picture}(31, 8.0)(0, 2.0)
\put(-0.2, 10){\line(1,0){25.4}}
	\put(-0.2, 5){\line(1,0){25.4}}
	\put(-0.2, 0){\line(1,0){25.4}}
	\put(0, 10){\line(0,-1){10}}
	\put(5, 10){\line(0,-1){10}}
	\put(9, 2.3){${}_{\cdots}$}
	\put(9, 7.3){${}_{\cdots}$}
	\put(20, 10){\line(0,-1){10}}
	\put(25, 10){\line(0,-1){10}}
	\put(25, 0){${}_{\,h}$}
\end{picture}(V)
\]
from Section \ref{Sec:MLA}. Then $\operatorname{trace-free}\colon \begin{picture}(31, 8.0)(0, 2.0)
\put(-0.2, 10){\line(1,0){25.4}}
	\put(-0.2, 5){\line(1,0){25.4}}
	\put(-0.2, 0){\line(1,0){25.4}}
	\put(0, 10){\line(0,-1){10}}
	\put(5, 10){\line(0,-1){10}}
	\put(9, 2.3){${}_{\cdots}$}
	\put(9, 7.3){${}_{\cdots}$}
	\put(20, 10){\line(0,-1){10}}
	\put(25, 10){\line(0,-1){10}}
\end{picture}(V)\to \begin{picture}(31, 8.0)(0, 2.0)
\put(-0.2, 10){\line(1,0){25.4}}
	\put(-0.2, 5){\line(1,0){25.4}}
	\put(-0.2, 0){\line(1,0){25.4}}
	\put(0, 10){\line(0,-1){10}}
	\put(5, 10){\line(0,-1){10}}
	\put(9, 2.3){${}_{\cdots}$}
	\put(9, 7.3){${}_{\cdots}$}
	\put(20, 10){\line(0,-1){10}}
	\put(25, 10){\line(0,-1){10}}
	\put(25, 0){${}_{\,0}$}
\end{picture}(V) $ is the projection onto the first summand. Let $\iota\colon \begin{picture}(31, 8.0)(0, 2.0)
\put(-0.2, 10){\line(1,0){25.4}}
	\put(-0.2, 5){\line(1,0){25.4}}
	\put(-0.2, 0){\line(1,0){25.4}}
	\put(0, 10){\line(0,-1){10}}
	\put(5, 10){\line(0,-1){10}}
	\put(9, 2.3){${}_{\cdots}$}
	\put(9, 7.3){${}_{\cdots}$}
	\put(20, 10){\line(0,-1){10}}
	\put(25, 10){\line(0,-1){10}}
	\put(25, 0){${}_{\,0}$}
\end{picture}(V)\to \begin{picture}(31, 8.0)(0, 2.0)
\put(-0.2, 10){\line(1,0){25.4}}
	\put(-0.2, 5){\line(1,0){25.4}}
	\put(-0.2, 0){\line(1,0){25.4}}
	\put(0, 10){\line(0,-1){10}}
	\put(5, 10){\line(0,-1){10}}
	\put(9, 2.3){${}_{\cdots}$}
	\put(9, 7.3){${}_{\cdots}$}
	\put(20, 10){\line(0,-1){10}}
	\put(25, 10){\line(0,-1){10}}
\end{picture}(V) $ be the inclusion, so $\operatorname{trace-free}$ and $\iota$ are adjoint maps. Then the image of \eqref{Eqn:tfWV} is the image of $\operatorname{trace-free}\circ \iota^0$, where $\iota^0 \colon \begin{picture}(31, 8.0)(0, 2.0)
\put(-0.2, 10){\line(1,0){25.4}}
	\put(-0.2, 5){\line(1,0){25.4}}
	\put(-0.2, 0){\line(1,0){25.4}}
	\put(0, 10){\line(0,-1){10}}
	\put(5, 10){\line(0,-1){10}}
	\put(9, 2.3){${}_{\cdots}$}
	\put(9, 7.3){${}_{\cdots}$}
	\put(20, 10){\line(0,-1){10}}
	\put(25, 10){\line(0,-1){10}}
\end{picture}(W) \to \begin{picture}(31, 8.0)(0, 2.0)
\put(-0.2, 10){\line(1,0){25.4}}
	\put(-0.2, 5){\line(1,0){25.4}}
	\put(-0.2, 0){\line(1,0){25.4}}
	\put(0, 10){\line(0,-1){10}}
	\put(5, 10){\line(0,-1){10}}
	\put(9, 2.3){${}_{\cdots}$}
	\put(9, 7.3){${}_{\cdots}$}
	\put(20, 10){\line(0,-1){10}}
	\put(25, 10){\line(0,-1){10}}
\end{picture}(V)  $ is the inclusion with adjoint $\pi^0$. So
\[
	\begin{aligned}
		\operatorname{im}(\operatorname{trace-free}\circ\iota^0) = \left( \ker(\pi^0\circ \iota) \right)^{\perp}& = \{K\in \begin{picture}(31, 8.0)(0, 2.0)
				\put(-0.2, 10){\line(1,0){25.4}}
				\put(-0.2, 5){\line(1,0){25.4}}
				\put(-0.2, 0){\line(1,0){25.4}}
				\put(0, 10){\line(0,-1){10}}
				\put(5, 10){\line(0,-1){10}}
				\put(9, 2.3){${}_{\cdots}$}
				\put(9, 7.3){${}_{\cdots}$}
				\put(20, 10){\line(0,-1){10}}
				\put(25, 10){\line(0,-1){10}}
				\put(25, 0){${}_{\,0}$}
				\end{picture}(\mathcal{E}_A) : K^0=0\}^{\perp} \\
				&= \{K\in \begin{picture}(31, 8.0)(0, 2.0)
				\put(-0.2, 10){\line(1,0){25.4}}
				\put(-0.2, 5){\line(1,0){25.4}}
				\put(-0.2, 0){\line(1,0){25.4}}
				\put(0, 10){\line(0,-1){10}}
				\put(5, 10){\line(0,-1){10}}
				\put(9, 2.3){${}_{\cdots}$}
				\put(9, 7.3){${}_{\cdots}$}
				\put(20, 10){\line(0,-1){10}}
				\put(25, 10){\line(0,-1){10}}
				\put(25, 0){${}_{\,0}$}
		\end{picture}(\mathcal{E}_A) : K^1=0\},
	\end{aligned}
\]
by Lemma \ref{Lem:K0K1}. Showing \ref{Thm:imK:1} $\iff$ \ref{Thm:imK:2}.

Finally we show \ref{Thm:imK:2} $\iff $ \ref{Thm:imK:3}. Observe that, on $\Lambda^2(V) = W\wedge W\oplus W\wedge I$, the linear map $\wedge I$, wedging with $I$, is injective on $W\wedge W$ and is zero on $W\wedge I$. On the other hand, $\lrcorner I$ is invertible on $W\wedge I$ and is zero on $W\wedge W$. So the left hand side of \eqref{Eqn:IKIII} is zero on $K^0,K^2,K^3,\dots ,K^\ell$ and never zero on $K^1$, so \eqref{Eqn:IKIII} holds if and only if $K^1=0$.
\end{proof}

Now we return to tractors and prolongation. Let $k_{a_{1}\cdots a_\ell}$ be a conformal Killing tensor with full prolongation $K_{A_{1}B_{1}\cdots A_\ell B_\ell}$ and let $I{}_{A}$ be a parallel non-null tractor. If $(M,\boldsymbol{c})$ is conformally flat, then by Proposition \ref{Prp:HRKS}, $\sigma$ is a local Killing scale if and only if any of the equivalent conditions in Theorem \ref{Thm:MLA} holds for $K$. For our main theorem, we strengthen this result by removing the reliance on local scales.
\begin{theorem}\label{Thm:HO1}
	Let $(M,\boldsymbol{c})$ be conformally flat. Let $k{}_{a_{1}\cdots a_\ell}\in \Gamma(\mathcal{E}{}_{(a_{1}\cdots a_\ell)_0}[2\ell] ) $ be a conformal Killing tensor and let $K{}_{A_{1} B_{1}\cdots A_\ell B_\ell} $ be its full prolongation. Let $I_A = D_A\sigma$ be a non-null parallel scale tractor. 

	Then $\sigma$ is a Killing scale for $k$ if and only if the tractor
	\begin{equation}\label{Eqn:KSl}
		I{}^{A_{1}}K{}_{A_{1}B_{1}A_{2}B_{2}\cdots A_\ell B_\ell}
	\end{equation} 
	vanishes after skewing each pair $A_iB_i$ with $I_{C_i}$ for each $i\ge 2$; that is, with an abuse of notation,
	 \begin{equation}\label{Eqn:IKIII}
		 I{}^{A_{1}}K{}_{A_1B_{1}[A_{2}B_{2}}I_{C_{2}]}\cdots {}_{[A_\ell B_\ell}I_{C_{\ell}]}=0,
	 \end{equation} 
	where all of the indices $A_i B_i$ belong to $K$.
\end{theorem}

\begin{proof}
	If $\sigma$ is a Killing scale, then it is a local Killing scale, so by Proposition \ref{Prp:HRKS}, the tractor $K$ satisfies the condition \ref{Thm:imK:1} of Theorem \ref{Thm:MLA}, so it satisfies condition \ref{Thm:imK:3} of that Theorem, which is \eqref{Eqn:KSl}.

	Since $I$ is parallel, the bundle $\begin{picture}(31, 8.0)(0, 2.0)
 	\put(-0.2, 10){\line(1,0){25.4}}
 	\put(-0.2, 5){\line(1,0){25.4}}
 	\put(-0.2, 0){\line(1,0){25.4}}
 	\put(0, 10){\line(0,-1){10}}
 	\put(5, 10){\line(0,-1){10}}
 	\put(9, 2.3){${}_{\cdots}$}
 	\put(9, 7.3){${}_{\cdots}$}
 	\put(20, 10){\line(0,-1){10}}
 	\put(25, 10){\line(0,-1){10}}
 \end{picture}(W)\to M$ is a parallel subbundle of $\begin{picture}(31, 8.0)(0, 2.0)
 	\put(-0.2, 10){\line(1,0){25.4}}
 	\put(-0.2, 5){\line(1,0){25.4}}
 	\put(-0.2, 0){\line(1,0){25.4}}
 	\put(0, 10){\line(0,-1){10}}
 	\put(5, 10){\line(0,-1){10}}
 	\put(9, 2.3){${}_{\cdots}$}
 	\put(9, 7.3){${}_{\cdots}$}
 	\put(20, 10){\line(0,-1){10}}
 	\put(25, 10){\line(0,-1){10}}
 \end{picture}(\mathcal{E}_A)$, so the tractor connection restricts to a connection on $\begin{picture}(31, 8.0)(0, 2.0)
 	\put(-0.2, 10){\line(1,0){25.4}}
 	\put(-0.2, 5){\line(1,0){25.4}}
 	\put(-0.2, 0){\line(1,0){25.4}}
 	\put(0, 10){\line(0,-1){10}}
 	\put(5, 10){\line(0,-1){10}}
 	\put(9, 2.3){${}_{\cdots}$}
 	\put(9, 7.3){${}_{\cdots}$}
 	\put(20, 10){\line(0,-1){10}}
 	\put(25, 10){\line(0,-1){10}}
\end{picture}(W)$ and the inclusion $\iota^0$ is parallel. The bundle map $\operatorname{trace-free}\colon \begin{picture}(31, 8.0)(0, 2.0)
	\put(-0.2, 10){\line(1,0){25.4}}
	\put(-0.2, 5){\line(1,0){25.4}}
	\put(-0.2, 0){\line(1,0){25.4}}
	\put(0, 10){\line(0,-1){10}}
	\put(5, 10){\line(0,-1){10}}
	\put(9, 2.3){${}_{\cdots}$}
	\put(9, 7.3){${}_{\cdots}$}
	\put(20, 10){\line(0,-1){10}}
	\put(25, 10){\line(0,-1){10}}
\end{picture}(\mathcal{E}_A)\to  \begin{picture}(31, 8.0)(0, 2.0)
\put(-0.2, 10){\line(1,0){25.4}}
	\put(-0.2, 5){\line(1,0){25.4}}
	\put(-0.2, 0){\line(1,0){25.4}}
	\put(0, 10){\line(0,-1){10}}
	\put(5, 10){\line(0,-1){10}}
	\put(9, 2.3){${}_{\cdots}$}
	\put(9, 7.3){${}_{\cdots}$}
	\put(20, 10){\line(0,-1){10}}
	\put(25, 10){\line(0,-1){10}}
	\put(25, 0){${}_{\,0}$}
\end{picture}(\mathcal{E}_A)$ is also a parallel bundle morphism. Hence, the composition $\phi := \operatorname{trace-free}\circ \iota^0$ is a parallel bundle map. It then follows that $\phi^{*} = \pi^0\circ \iota$ is also parallel. The bundles $\begin{picture}(31, 8.0)(0, 2.0)
 	\put(-0.2, 10){\line(1,0){25.4}}
 	\put(-0.2, 5){\line(1,0){25.4}}
 	\put(-0.2, 0){\line(1,0){25.4}}
 	\put(0, 10){\line(0,-1){10}}
 	\put(5, 10){\line(0,-1){10}}
 	\put(9, 2.3){${}_{\cdots}$}
 	\put(9, 7.3){${}_{\cdots}$}
 	\put(20, 10){\line(0,-1){10}}
 	\put(25, 10){\line(0,-1){10}}
\end{picture}(W)$, $\begin{picture}(31, 8.0)(0, 2.0)
	\put(-0.2, 10){\line(1,0){25.4}}
	\put(-0.2, 5){\line(1,0){25.4}}
	\put(-0.2, 0){\line(1,0){25.4}}
	\put(0, 10){\line(0,-1){10}}
	\put(5, 10){\line(0,-1){10}}
	\put(9, 2.3){${}_{\cdots}$}
	\put(9, 7.3){${}_{\cdots}$}
	\put(20, 10){\line(0,-1){10}}
	\put(25, 10){\line(0,-1){10}}
\end{picture}(\mathcal{E}_A)$ and $\begin{picture}(31, 8.0)(0, 2.0)
\put(-0.2, 10){\line(1,0){25.4}}
	\put(-0.2, 5){\line(1,0){25.4}}
	\put(-0.2, 0){\line(1,0){25.4}}
	\put(0, 10){\line(0,-1){10}}
	\put(5, 10){\line(0,-1){10}}
	\put(9, 2.3){${}_{\cdots}$}
	\put(9, 7.3){${}_{\cdots}$}
	\put(20, 10){\line(0,-1){10}}
	\put(25, 10){\line(0,-1){10}}
	\put(25, 0){${}_{\,0}$}
\end{picture}(\mathcal{E}_A)$, naturally carry actions of $O(W)$, $GL(\mathcal{E}_A)$ and $O(\mathcal{E}_A)$, respectively, so each has an action of $O(W)$ by restriction to the natural subgroup of $O(\mathcal{E}_A)$ isomorphic to $O(W)$. The map $\iota^0$ is $O(W)$-equivariant and $\operatorname{trace-free}$ is $O(\mathcal{E}_A)$-equivariant, so it is $O(W)$-equivariant. Hence, the self-adjoint composition $\phi^*\circ \phi$ is a parallel $O(W)$-equivariant endomorphism of $\begin{picture}(31, 8.0)(0, 2.0)
 	\put(-0.2, 10){\line(1,0){25.4}}
 	\put(-0.2, 5){\line(1,0){25.4}}
 	\put(-0.2, 0){\line(1,0){25.4}}
 	\put(0, 10){\line(0,-1){10}}
 	\put(5, 10){\line(0,-1){10}}
 	\put(9, 2.3){${}_{\cdots}$}
 	\put(9, 7.3){${}_{\cdots}$}
 	\put(20, 10){\line(0,-1){10}}
 	\put(25, 10){\line(0,-1){10}}
\end{picture}(W)$. Then it follows that there is a natural decomposition into subbundles, given by Littlewood's branching rules (see \cite[\S 25.3]{fulton_representation_2004} or \cite{littlewood_theory_1940}),
\begin{equation}\label{Eqn:Wdecomp}
\begin{picture}(31, 8.0)(0, 2.0)
 	\put(-0.2, 10){\line(1,0){25.4}}
 	\put(-0.2, 5){\line(1,0){25.4}}
 	\put(-0.2, 0){\line(1,0){25.4}}
 	\put(0, 10){\line(0,-1){10}}
 	\put(5, 10){\line(0,-1){10}}
 	\put(9, 2.3){${}_{\cdots}$}
 	\put(9, 7.3){${}_{\cdots}$}
 	\put(20, 10){\line(0,-1){10}}
 	\put(25, 10){\line(0,-1){10}}
\end{picture}(W) = \bigoplus_{i=1}^m \mathbb{V}_i,
\end{equation} 
where each $\mathbb{V}_i$ is an absolutely irreducible $O(W)$-equivariant subbundle of $\begin{picture}(31, 8.0)(0, 2.0)
 	\put(-0.2, 10){\line(1,0){25.4}}
 	\put(-0.2, 5){\line(1,0){25.4}}
 	\put(-0.2, 0){\line(1,0){25.4}}
 	\put(0, 10){\line(0,-1){10}}
 	\put(5, 10){\line(0,-1){10}}
 	\put(9, 2.3){${}_{\cdots}$}
 	\put(9, 7.3){${}_{\cdots}$}
 	\put(20, 10){\line(0,-1){10}}
 	\put(25, 10){\line(0,-1){10}}
	\end{picture}(W)$ and the decomposition is multiplicity free, that is the $\mathbb{V}_i$ are pairwise non-isomorphic as $O(W)$-modules. Moreover, each $\mathbb{V}_i$ is a parallel subbundle, since the connection preserves $I$ and $h$, so its holonomy is a subgroup of $O(W)$. Then by Schur's Lemma there are constants $c_i$ such that $(\phi^*\circ \phi)|_{\mathbb{V}_i} = c_i \operatorname{id}|_{\mathbb{V}_i}$, some of the $c_i$ may be zero. Suppose that $K^1=0$. Then by Theorem \ref{Thm:MLA}, this is equivalent to $K\in \operatorname{im}(\phi)$, so $K = \phi(R)$ for some $R\in \Gamma(\begin{picture}(31, 8.0)(0, 2.0)
 	\put(-0.2, 10){\line(1,0){25.4}}
 	\put(-0.2, 5){\line(1,0){25.4}}
 	\put(-0.2, 0){\line(1,0){25.4}}
 	\put(0, 10){\line(0,-1){10}}
 	\put(5, 10){\line(0,-1){10}}
 	\put(9, 2.3){${}_{\cdots}$}
 	\put(9, 7.3){${}_{\cdots}$}
 	\put(20, 10){\line(0,-1){10}}
 	\put(25, 10){\line(0,-1){10}}
\end{picture}(W))$. Then let $T = K^0 = \pi^0\circ \iota(K) =\phi^*(K)= \phi^*\circ \phi(R)$. Suppose that $T = \sum_{i=1}^{m} T_i$, where $T_i\in \Gamma(\mathbb{V}_i)$. Then let
\[
	\hat{R} = \sum_{c_i\neq 0} \frac{1}{c_i}T_i = \sum_{c_i\neq 0} R_i,
\]
where $R = \sum_{i=1}^{m} R_i$ with $R_i\in \mathbb{V}_i$. For a map between spaces equipped with bilinear forms $\ker(\phi^*\circ \phi) / \ker(\phi)$ is naturally isomorphic to $\operatorname{im}(\phi)\cap (\operatorname{im}(\phi)^{\perp}$. Inside $\begin{picture}(31, 8.0)(0, 2.0)
		\put(-0.2, 10){\line(1,0){25.4}}
		\put(-0.2, 5){\line(1,0){25.4}}
		\put(-0.2, 0){\line(1,0){25.4}}
		\put(0, 10){\line(0,-1){10}}
		\put(5, 10){\line(0,-1){10}}
		\put(9, 2.3){${}_{\cdots}$}
		\put(9, 7.3){${}_{\cdots}$}
		\put(20, 10){\line(0,-1){10}}
		\put(25, 10){\line(0,-1){10}}
		\put(25, 0){${}_{\,0}$}
\end{picture}(\mathcal{E}_A)$ the image of $\phi$ is $\{K\in \begin{picture}(31, 8.0)(0, 2.0)
		\put(-0.2, 10){\line(1,0){25.4}}
		\put(-0.2, 5){\line(1,0){25.4}}
		\put(-0.2, 0){\line(1,0){25.4}}
		\put(0, 10){\line(0,-1){10}}
		\put(5, 10){\line(0,-1){10}}
		\put(9, 2.3){${}_{\cdots}$}
		\put(9, 7.3){${}_{\cdots}$}
		\put(20, 10){\line(0,-1){10}}
		\put(25, 10){\line(0,-1){10}}
		\put(25, 0){${}_{\,0}$}
\end{picture}(\mathcal{E}_A) : K^1=0\}$, which has orthogonal complement $\{K\in \begin{picture}(31, 8.0)(0, 2.0)
		\put(-0.2, 10){\line(1,0){25.4}}
		\put(-0.2, 5){\line(1,0){25.4}}
		\put(-0.2, 0){\line(1,0){25.4}}
		\put(0, 10){\line(0,-1){10}}
		\put(5, 10){\line(0,-1){10}}
		\put(9, 2.3){${}_{\cdots}$}
		\put(9, 7.3){${}_{\cdots}$}
		\put(20, 10){\line(0,-1){10}}
		\put(25, 10){\line(0,-1){10}}
		\put(25, 0){${}_{\,0}$}
\end{picture}(\mathcal{E}_A) : K^0=0\}$, and these have trivial intersection, both by Lemma \ref{Lem:K0K1}. So, $\ker(\phi) = \ker(\phi^*\circ \phi)$. Then
\[
\phi(\hat{R}) = \sum_{c_i\neq 0} \phi(R_i) = \sum_{i=1}^{m}  \phi(R_i) = \phi(R) = K,
\] 
Hence, there is a formula, which can be applied globally to give a parallel tractor $\hat{R}$ orthogonal to $I$ with trace-free part $K$, hence the top slot of $\hat{R}$ is a Killing tensor that is a trace-adjustment of $k$.
\end{proof}
It is desirable to have an extension of the above theorem when $(M,\boldsymbol{c})$ is not assumed conformally flat, but $K$ and $I$ are parallel. As we will see in the following section, there are families of examples of non-conformally flat manifolds with plenty of parallel tractors.

\begin{theorem}\label{Thm:HO2}
	Let $k_{a_{1}\cdots a_\ell}\in \Gamma(\mathcal{E}{}_{(a_{1}\cdots a_\ell)_0}[2\ell] )$ be a normal solution of the conformal Killing equation with parallel full prolongation $K_{A_{1}B_{1}\cdots A_\ell B_\ell}$. Let $I_A = D_A\sigma$ be a non-null parallel scale tractor. If \eqref{Eqn:IKIII} holds, then $\sigma$ is a Killing scale for $k$.
\end{theorem}
\begin{proof}
	The decomposition \eqref{Eqn:Wdecomp} is well defined and parallel as long as $I$ is parallel, even if $(M,\boldsymbol{c})$ is not conformally flat. The map taking $K\mapsto \hat{R}$ in the proof of Theorem \ref{Thm:HO1} therefore takes parallel $K$ to parallel $\hat{R}$, which is orthogonal to $I$. Hence the top slot of $\hat{R}$ is a Killing tensor that is a trace-adjustment of $k$.
\end{proof}
Recall the formula for $\operatorname{trace-free}$ from Section \ref{Sec:KT},
\[
R{}_{ABCD} = K{}_{ABCD} +(P\owedge h)_{ABCD},
\] 
where 
\[
	P{}_{BD}  = \frac{1}{n}\left(R{}_{AB}{}^{A}{}_{D} - \frac{1}{2(n+1)}R{}_{AC}{}^{AC} h{}_{BD}\right).
\]
From this we can compute the values $c_i$ from the above proof when $\ell=2$. We decompose 
\[
\begin{picture}(22, 11.0)(0, 5.0)
	\put(-0.2, 16){\line(1,0){16.4}}
	\put(-0.2, 8){\line(1,0){16.4}}
	\put(-0.2, 0){\line(1,0){16.4}}
	\put(0, 16){\line(0,-1){16}}
	\put(8, 16){\line(0,-1){16}}
	\put(16, 16){\line(0,-1){16}}
\end{picture} (W) = \begin{picture}(22, 11.0)(0, 5.0)
	\put(-0.2, 16){\line(1,0){16.4}}
	\put(-0.2, 8){\line(1,0){16.4}}
	\put(-0.2, 0){\line(1,0){16.4}}
	\put(0, 16){\line(0,-1){16}}
	\put(8, 16){\line(0,-1){16}}
	\put(16, 16){\line(0,-1){16}}
	\put(16, 0){${}_{\,0}$}
\end{picture}(W)\oplus \begin{picture}(22, 7.0)(0, 1.0)
	\put(-0.2, 8){\line(1,0){16.4}}
	\put(-0.2, 0){\line(1,0){16.4}}
	\put(0, 8){\line(0,-1){8}}
	\put(8, 8){\line(0,-1){8}}
	\put(16, 8){\line(0,-1){8}}
	\put(16, 0){${}_{\,0}$}
\end{picture}(W)\oplus \mathbb{R} = \mathbb{V}_1\oplus \mathbb{V}_2\oplus \mathbb{V}_3,
\] 
where the trace factors $\begin{picture}(22, 11.0)(0, 5.0)
	\put(-0.2, 16){\line(1,0){16.4}}
	\put(-0.2, 8){\line(1,0){16.4}}
	\put(-0.2, 0){\line(1,0){16.4}}
	\put(0, 16){\line(0,-1){16}}
	\put(8, 16){\line(0,-1){16}}
	\put(16, 16){\line(0,-1){16}}
	\put(16, 0){${}_{\,g}$}
\end{picture}(W) =  \begin{picture}(22, 7.0)(0, 1.0)
	\put(-0.2, 8){\line(1,0){16.4}}
	\put(-0.2, 0){\line(1,0){16.4}}
	\put(0, 8){\line(0,-1){8}}
	\put(8, 8){\line(0,-1){8}}
	\put(16, 8){\line(0,-1){8}}
	\put(16, 0){${}_{\,0}$}
\end{picture}(W)\oplus \mathbb{R}$ consists of tractors $\hat{P}\owedge \hat{h}$, where $\hat{P}\in S^2(W)$. We have
\[
	\begin{aligned}
		K{}_{ABCD} &= R{}_{ABCD}  - (P\owedge h)_{ABCD} \\
			   &= R{}_{ABCD} - \frac{1}{n}(R{}_{E*}{}^{E}{}_{*} \owedge h)_{ABCD}+\frac{1}{2n(n+1)}R{}_{AC}{}^{AC} (h\owedge h)_{ABCD}
	\end{aligned}
\] 
So projecting onto $S^2(\Lambda^2(W))$ gives
\[
	K{}^{0}{}_{ABCD} = R{}_{ABCD} -\frac{1}{n}(R{}_{E*}{}^{E}{}_{*} \owedge \hat{h})_{ABCD}+\frac{1}{2n(n+1)}R{}_{AC}{}^{AC} (\hat{h}\owedge \hat{h})_{ABCD}
\] 
So the trace-free part of $K{}^{0}$, with respect to $\hat{h}$ is the trace-free part of $R$ and $c_1=1$. We can compute the trace of $K{}^{0}$ directly, giving
\[
K{}^{0}{}_{AB}{}^{A}{}_{D} = \frac{1}{n}R{}_{AB}{}^{A}{}_{D} -\frac{1}{n(n+1)}R{}_{AC}{}^{AC} \hat{h}{}_{BD}.
\] 
Giving $c_2=\frac{1}{n}$ and $c_3=0$.

\section{Curved Examples}\label{Sec:Ex}

The main results developed in this paper apply universally to conformally flat manifolds, which are locally conformally equivalent to open subsets of the $(p,q)$-sphere. Some results only require normal solutions to the conformal Killing equations, such as Theorems \ref{Thm:IWI} and \ref{Thm:HO2}. Normal solutions correspond exactly with prolongations which are parallel for the standard tractor connection. It is reasonable to ask to what degree such solutions can exist on non-conformally flat manifolds.

Let $n_{1},n_{2}\ge 4$. Let $M_{1}$ be the model space $\mathbb{S}^{n_{1}} = G / P$, for $G = SO(n_{1}+1,1)$ and $P$ the parabolic subgroup stabilising a null line. Then $M_{1}$ has a full space of parallel tractors, of any constant length. We will focus on the round sphere metric, which corresponds to a parallel tractor $I^1$ of constant negative length. After suitable normalisation and trivialisation of densities, we may take $J_{1} / n_{1} = +1$. In order to construct a special Einstein product, we need an Einstein manifold $(M_{2},g_{2})$ of negative scalar curvature, which crucially does not need to be conformally flat. Some examples include: the hyperbolic spaces $\mathbb{C}H^m$, $\mathbb{H}H^m$ and the Cayley hyperbolic plane $\mathbb{O}\mathbb{H}^2$, suitable quotients of these by lattices $\Gamma$, a product $\mathbb{H}^a\times\mathbb{H}^b$ with their Einstein constants scaled to be equal is negative Einstein, but not conformally flat. We could also take $M_{1}$ negatively curved hyperbolic space $\mathbb{H}^{n_{1}}$ and $M_{2}$ positively curved, such as $\mathbb{C}P^n, \mathbb{H}P^n, \mathbb{O}P{}^{2}$, or $\mathbb{S}^k\times \mathbb{S}^k$. Finally, we may also take $M_{1}$ flat $\mathbb{R}^n$ and $M_{2}$ Ricci flat, such as a Calabi-Yau manifold--although, in this case, the below dimension counts may differ.

Now we fix $M_{1} = \mathbb{S}^{n_{1}}$ and $(M_{2},g_{2})$ to be an Einstein manifold such that $J_{2} / n_{2} = -1$ and consider the conformal manifold $(M = M_{1}\times M_{2}, \boldsymbol{c} = [g_{1}\times g_{2}])$. By Theorem \ref{Thm:SEP}, the tractor bundle of $M$ is $\mathcal{T}\cong (I^{1})^{\perp}\oplus (I^{2})^{\perp}\subset \mathcal{T}_1\oplus \mathcal{T}_2$. Denote by $\pi^i \colon \mathcal{T}\to (I^i)^{\perp}\subset \mathcal{T}_i $. The tractor curvature on $\mathcal{T}$ then splits as $\kappa = \kappa^{1}\oplus \kappa^{2}$, with $\kappa^{i}$ the curvature of $\mathcal{T}_i$ restricted to $(I^i)^{\perp}$. By assumption, $\kappa^1=0$ and $\kappa^2$ captures the full tractor curvature of $M_{2}$, since $I{}^{2}$ in $\mathcal{T}_2$ is orthogonal to the tractor curvature. Hence, $M$ is not conformally flat whenever $M_{2}$ is not conformally flat.

Let $I$ be a parallel tractor on $M$, such as any constant section of the subbundle of $\mathcal{T}$ isomorphic to $(I^1)^{\perp}$. Suppose that $|I|^2\neq 0$. Then the parallel subbundle $I^{\perp}\cap (I^1)^{\perp}$ is non-degenerate and has dimension $n_{1}$ when $\pi^1(I)\neq 0$ and dimension $n_{1}+1$ when $\pi^1(I)=0$. So we conclude that the space of conformal Killing tensors of rank $\ell$ for which $\sigma = X^AI_A$ is a strong Killing scale is of dimension at least
\[
	\frac{(n_{1}+2\ell-4)(n_{1}+2\ell-3)(n_{1}+2\ell-2)(n_{1}+\ell-4)!(n_{1}+\ell-5)!}{(\ell+1)!\ell!(n_{1}-2)!(n_{1}-4)!},
\] 
which comes from the Weyl dimension formula, as the dimension of $ \begin{picture}(31, 8.0)(0, 2.0)
\put(-0.2, 10){\line(1,0){25.4}}
	\put(-0.2, 5){\line(1,0){25.4}}
	\put(-0.2, 0){\line(1,0){25.4}}
	\put(0, 10){\line(0,-1){10}}
	\put(5, 10){\line(0,-1){10}}
	\put(9, 2.3){${}_{\cdots}$}
	\put(9, 7.3){${}_{\cdots}$}
	\put(20, 10){\line(0,-1){10}}
	\put(25, 10){\line(0,-1){10}}
	\put(25, 0){${}_{\,0}$}
\end{picture}(U) $, where $U$ has dimension $n_{1}\ge 4$. 

The space of conformal Killing tensors for which $\sigma$ is a Killing scale contains the top slots of the trace-free parts $K$ of all parallel tractors $R$ with indices in $I^{\perp}\cap (I^1)^{\perp}$. In general $K$ need not lie in $(I^1)^{\perp}$, but it is still parallel.

Now suppose that all indices of \(K\) lie in \((I^1)^\perp\) and that \(\pi^2(I)\ne0\). Then $K$ must satisfy \eqref{Eqn:IKIII}. However, if $I\not\in (I^1)^{\perp}$ (which requires $M_{2}$ to have another Einstein scale, which does not hold for any of the examples we listed above), then the condition \eqref{Eqn:IKIII} becomes
\[
	I^{A_{1}}K_{A_{1}B_{1}[A_{2}B_{2}}I_{C_{2}]}\cdots{}_{[A_\ell B_\ell}I_{C_\ell]} = \pi^1(I)^{A_{1}}K_{A_{1}B_{1}[A_{2}B_{2}}I_{C_{2}]}\cdots{}_{[A_\ell B_\ell}I_{C_\ell]}=0,
\] 
which can only hold if $\pi^1(I)^{A_{1}}K_{A_{1}B_{1}\cdots A_\ell B_\ell}=0$, since $\wedge \pi^2(I)$ is injective on $\Lambda^2((I^1)^{\perp})$. Hence, \eqref{Eqn:IKIII} reduces to $K\in I^{\perp}$, or that $\sigma $ is a strong Killing scale.

\bibliographystyle{amsplain}
\bibliography{DG}

\end{document}